\documentclass[11pt,oneside,reqno]{article}
\usepackage[width=15.7cm,height=24cm,centering]{geometry}
\usepackage{amsmath,amsthm,amssymb,color,bm}
\usepackage[authoryear]{natbib}
\usepackage{booktabs}

\def\d{\,\mathrm{d}}
\usepackage{tikz}
\RequirePackage{amsthm,amsmath,amsfonts,amssymb,mathrsfs,dsfont,mathtools,thmtools}
\RequirePackage[colorlinks,citecolor=blue,urlcolor=blue,linkcolor=blue]{hyperref}
\usepackage[capitalize]{cleveref}
\RequirePackage{graphicx}

\crefname{equation}{}{}
\crefformat{equation}{\textup{(#2#1#3)}}
\crefrangeformat{equation}{\textup{(#3#1#4)--(#5#2#6)}}
\crefdefaultlabelformat{#2\textup{#1}#3}
\crefname{lemma}{Lemma}{Lemmas}
\crefname{page}{p.}{pp.}
\usepackage[normalem]{ulem}

\usepackage{bbm}
\usepackage{mathrsfs}
\usepackage{graphicx}
\usepackage{tikz}
\usepackage{enumitem}
\usetikzlibrary{arrows, automata}
\usetikzlibrary{calc}
\usetikzlibrary{positioning}

\numberwithin{equation}{section}
\allowdisplaybreaks[4]

\theoremstyle{plain}
\newtheorem{theorem}{Theorem}[section]
\newtheorem{proposition}{Proposition}[section]
\newtheorem{lemma}{Lemma}[section]
\newtheorem{corollary}{Corollary}[section]

\theoremstyle{definition}
\newtheorem{definition}{Definition}[section]
\newtheorem{example}{Example}[section]

\makeatletter

\newcount\minute
\newcount\hour
\newcount\hourMins
\def\now{%
\minute=\time%
\hour=\time \divide \hour by 60%
\hourMins=\hour \multiply\hourMins by 60%
\advance\minute by -\hourMins%
\zeroPadTwo{\the\hour}:\zeroPadTwo{\the\minute}%
}
\def\zeroPadTwo#1{\ifnum #1<10 0\fi#1}

\renewcommand{\cite}{\citet}

\def\^#1{\ifmmode {\mathaccent"705E #1} \else {\accent94 #1} \fi}
\def\~#1{\ifmmode {\mathaccent"707E #1} \else {\accent"7E #1} \fi}

\def\*#1{#1^\ast}
\edef\-#1{\noexpand\ifmmode {\noexpand\bar{#1}} \noexpand\else \-#1\noexpand\fi}
\def\>#1{\vec{#1}}
\def\.#1{\dot{#1}}

\def\atop{\@@atop}
\def\*#1{\mathscr{#1}}

\renewcommand{\le}{\leqslant}

\renewcommand{\ge}{\geqslant}

\newcommand{\IE}{\mathbb{E}}

\newcommand{\Var}{\mathop{\mathrm{Var}}\nolimits}
\newcommand{\Cov}{\mathop{\mathrm{Cov}}}

\let\PlanarStandardMid\mid
\def\mid{\vert}

\def\beqn#1\eeqn{\begin{align}#1\end{align}}
\def\beq#1\eeq{\begin{align*}#1\end{align*}}

\usepackage{graphicx}
\usepackage{latexsym}
\usepackage{amsmath,amsthm,amssymb,amscd}
\usepackage{epsf,amsmath}

\def\E{{\IE}}

\def\d{\,\mathrm{d}}

\renewcommand\section{\@startsection {section}{1}{\z@}%
{-2.8ex \@plus -.8ex \@minus -.2ex}%
{1.0ex \@plus.2ex}%
{\center\small\sc\mathversion{bold}}}

\def\subsection#1{\@startsection {subsection}{2}{0pt}%
{-2.8ex \@plus -.8ex \@minus -.2ex}%
{.8ex \@plus.2ex}%
{\bf\mathversion{bold}}{#1}}

\def\subsubsection#1{\@startsection{subsubsection}{3}{0pt}%
{\medskipamount}%
{-10pt}%
{\normalsize\itshape}{\texorpdfstring{\kern-2.2ex. #1.}{#1}}}

\def\blfootnote{\xdef\@thefnmark{}\@footnotetext}

\makeatother

\let\mid\PlanarStandardMid
\providecommand{\R}{\mathbb R}
\providecommand{\Pp}{\mathbb P}
\providecommand{\1}{\mathbf 1}
\DeclareMathOperator{\Law}{Law}
\providecommand{\dBL}{d_{\mathrm{BL}}}
\hypersetup{pdftitle={A Boundary Comparison Method for Gaussian Fluctuations in Planar Euclidean Optimization},
  pdfauthor={Zi-Yao Su},
  pdfsubject={Probability; geometric optimization; central limit theorems}}

\begin{document}
\title{A Boundary Comparison Method for Gaussian Fluctuations in Planar Euclidean Optimization}
\author{Zi-Yao Su}
\date{\it The Chinese University of Hong Kong}
\maketitle

\noindent{\bf Abstract:}
We develop a boundary comparison method for Gaussian fluctuations of optimization costs on planar random points, including traveling-salesman length, minimum perfect-matching length, and the Held--Karp relaxation. A local change can alter an optimal solution far from the modified points, making a direct decomposition into local contributions difficult. Our criterion requires the difference between the changes in the full and local costs caused by the same replacement to be controlled at the boundary scale, allowing a smaller-order $L^2$ excess. Under this condition, the centered difference between the total cost and the sum of independent cell costs vanishes in $L^2$. The multiscale argument shows that a nonvanishing variance of the cutting error would produce Gaussian tails incompatible with the boundary control. We obtain joint Gaussian limits for the three centered costs, using a conditional version for matching, and prove that the limiting covariance is positive definite in the Euclidean case. Convergence holds in quadratic Wasserstein distance for Poisson and fixed sample sizes. The results extend jointly to finitely many planar norms and to fixed sampling laws whose absolutely continuous component has a density positive almost everywhere on the square, with an arbitrary singular component.



\section{Introduction}\label{sec:introduction}

Traveling-salesman length, minimum perfect-matching length, and the Held--Karp lower bound are global functions of a point configuration. A local change can alter the optimizer far from the modified points, so their fluctuations are not readily reduced to sums of independent local contributions. We develop a boundary comparison method that works directly with the optimal values.

Classical partitions and boundary repairs compare a global optimum with the sum of the optima in smaller squares. Such comparisons give first-order limits for costs of order $\sqrt n$ \citep{BeardwoodHaltonHammersley1959,Steele1981,Yukich1998}. A cutting error proportional to the total length of the cell boundaries can be negligible relative to $\sqrt n$ and still be too large on the order-one fluctuation scale. The question addressed here is whether the random part of that error can be made small.

Let $D=[0,1]^2$, let $\mathcal P_\lambda$ be a homogeneous Poisson process on $D$, and let $\mathcal C_k$ partition $D$ into $k^2$ equal squares. Consider square costs $H_Q$ that are translation invariant and homogeneous of degree one. Write the cutting error as
\begin{equation}\label{eq:intro-remainder}
 R_{\lambda,k}=H_D(\mathcal P_\lambda)
       -\sum_{Q\in\mathcal C_k}H_Q(\mathcal P_\lambda\cap Q).
\end{equation}
Centering this decomposition gives
\[
\begin{aligned}
 H_D(\mathcal P_\lambda)-\E H_D(\mathcal P_\lambda)
 & =\sum_{Q\in\mathcal C_k}
 \bigl[H_Q(\mathcal P_\lambda\cap Q)-\E H_Q(\mathcal P_\lambda\cap Q)\bigr]\\
 &\quad+R_{\lambda,k}-\E R_{\lambda,k}.
\end{aligned}
\]
The cell terms are independent under Poisson sampling. The cutting error itself need not tend to zero: its mean is removed by centering, and only $R_{\lambda,k}-\E R_{\lambda,k}$ remains to be controlled. Proving that this remainder vanishes in $L^2$ is the main step of the method.

To control its variance, replace the points in a square $B\subset Q$ by an independent configuration with the same law. The difference between the resulting changes in $H_D$ and $H_Q$, denoted by $\Delta_{B,Q}H$, is exactly the change in $R_{\lambda,k}$. Thus the same replacement is compared in the full square and in the cell containing $B$. If $b$ is the side length of $B$, our condition is
\begin{equation}\label{eq:intro-excess}
 \E\bigl[(|\Delta_{B,Q}H|-Cb)_+^2\bigr]
       \le b^2\eta(\lambda b^2),
 \qquad \eta(u)\longrightarrow0,
\end{equation}
uniformly over the relevant nested squares. The part of the response exceeding $Cb$ therefore has $L^2$ norm $o(b)$ as the expected number $\lambda b^2$ of points in $B$ grows; exceptional configurations are allowed. For a finite vector of costs, the condition is imposed coordinatewise. We prove that it yields
\begin{equation}\label{eq:intro-cut-limit}
 \|R_{\lambda,k}-\E R_{\lambda,k}\|_2\longrightarrow0
 \qquad\text{as }\lambda/k^2\longrightarrow\infty.
\end{equation}
uniformly in $k$.

Dimension two is critical for this estimate. At one scale, a response bound of order $b$ gives a squared response of order $b^2$, while there are of order $b^{-2}$ blocks. Variance tensorization \citep{EfronStein1981} therefore gives only a constant bound, even when the excess in \eqref{eq:intro-excess} is zero. Refining that single partition does not make the estimate vanish.

Comparing several scales supplies the missing information. An orthogonal decomposition bounds the sum of the squared interaction errors across levels, so the interaction error is small on average at some level. At that level, the cutting error, centered conditionally on a block's exterior and divided by the block side length, admits two $L^2$ approximations, with errors averaged over blocks according to their areas. The boundary comparison makes these fluctuations close to random variables in a fixed bounded interval. The orthogonal decomposition gives a second approximation by sums of terms that are independent once the exterior is fixed; truncation makes the individual terms small after this normalization. If the cutting variance stayed bounded away from zero, a conditional normal approximation would give a nonvanishing averaged tail probability outside the bounded interval. This contradicts the first approximation. Theorem~\ref{thm:multiscale} gives the resulting variance bound, and Corollary~\ref{cor:average} yields \eqref{eq:intro-cut-limit}.

To obtain the Gaussian limit, we first fix the expected number of points in each cell and let the number of cells grow. Translation invariance and degree-one homogeneity turn the centered cell sum into a normalized sum of iid costs, to which the ordinary central limit theorem applies. We then increase the expected count per cell. The uniform cutting estimate removes the approximation error and makes the Gaussian covariance matrices converge. The argument also includes the standardized Poisson count; Theorem~\ref{thm:boundary} states the joint limit. Fixed sample sizes are treated by a further conditioning argument.

\begin{samepage}
\subsection{Joint limits and changes of sampling law}

For uniform samples in the unit square, the method gives the following joint limit for the three Euclidean costs.

\begin{corollary}[Uniform Euclidean samples]\label{cor:intro-uniform}
Let $X_1,X_2,\ldots$ be independent uniform points in $D$. For even $n$, let $T_n$ be the shortest tour length, $M_n$ the minimum perfect-matching length, and $L_n$ the Held--Karp lower bound, all computed on $\{X_1,\ldots,X_n\}$ with Euclidean distance. There is a positive-definite $3\times3$ matrix $S_1$ such that
\begin{equation}\label{eq:intro-euclidean-limit}
 \begin{pmatrix}
 T_n-\E T_n\\
 M_n-\E M_n\\
 L_n-\E L_n
 \end{pmatrix}
 \xrightarrow[\substack{n\to\infty\\ n\ {\rm even}}]{W_2}
 \mathcal N_3(0,S_1).
\end{equation}
The covariance matrices also converge to $S_1$.
\end{corollary}
\end{samepage}

The fluctuations have order one and are centered by the exact expectations. Positive definiteness means that every nonzero linear combination of the three centered costs has a nondegenerate Gaussian limit; positivity of the three marginal variances alone would not give this conclusion. In particular, the additive difference $T_n-L_n$ has a nondegenerate centered Gaussian limit. The Gaussian convergence questions for tour and matching lengths were listed as open in \citet[Section~5]{Chatterjee2014Survey}.

Theorem~\ref{thm:applications} extends the joint limit to any fixed finite collection of these costs computed in planar norms. For odd sample sizes, matching is completed by adding a fixed corner point. The sampling law may also be any fixed Borel probability measure
\[
 \mu=g\,dx+\tau\quad\text{on }D,
 \qquad g>0\quad\text{almost everywhere},\qquad \tau\perp dx.
\]
Here $g$ is the density of the absolutely continuous component and $\tau$ is an arbitrary singular measure. The density need not be bounded above, bounded away from zero, or square integrable. For a chosen cost vector, write $H_{n,\mu}$ for its value on $n$ independent points with law $\mu$, using the corner completion for matching. Let $m$ be the uniform law on $D$ and set
\[
 \beta=\lim_{n\to\infty}n^{-1/2}\E H_{n,m}.
\]
If $S_1$ denotes the limiting covariance for uniform samples of fixed size and $a=\int_D\sqrt g\,dx$, then
\begin{equation}\label{eq:intro-general-law}
 H_{n,\mu}-\E H_{n,\mu}
 \xrightarrow{W_2}
 \mathcal N\!\left(0,S_1+\frac{1-a^2}{4}\beta\beta^\top\right).
\end{equation}
Thus changing the sampling law adds a positive semidefinite covariance correction of rank at most one. For Poisson input, the limiting covariance of the costs is $S_1+\beta\beta^\top/4$, independent of $\mu$ in this class. The limiting cross-covariance between the centered costs and the standardized count is $a\beta/2$, and the first-order mean coefficient is $a\beta$.

There is also a deterministic characterization of degeneracy. Under the deterministic boundary conditions used for these applications, a linear combination of the costs has zero limiting variance if and only if it is bounded over all admissible finite labeled configurations in $D$. This characterization, proved in Section~\ref{sec:nondegeneracy}, yields the positive definiteness in Corollary~\ref{cor:intro-uniform}.

\subsection{Geometric comparisons and related work}

The boundary comparisons take different forms for the three problems. For tours, cutting and reconnecting along a square boundary controls the change in length. For matching, individual cells can contain odd numbers of points even when the full sample has even size, so the local configurations are completed at the boundary. With Poisson input, the cell configurations remain independent after conditioning on their count parities. We control the conditional variance of the cutting error at fixed parities and then the variance of its conditional mean. For Held--Karp, weighted cut constraints allow boundary comparisons through the addition and subsequent elimination of auxiliary labels.

Partitioning and boundary repair have long been used to study the first-order size of Euclidean optimization costs \citep{BeardwoodHaltonHammersley1959,Karp1977,Steele1981,Jaillet1992}. The boundary-functional approach develops these comparisons \citep{RedmondYukich1994,RedmondYukich1996,Yukich1998}; quantitative estimates and broader accounts are given in \citet{KooLee2007,Yukich2013Survey}. Here geometric comparisons bound the response of the cutting error to a local replacement. The multiscale argument then controls the fluctuations of that error.

The variance argument has two precedents in random-field models. \citet{AizenmanWehr1990Rounding} use the incompatibility between a surface bound and a nondegenerate Gaussian fluctuation at the critical dimension. \citet{Chatterjee2018RFIMDecay} uses a Gaussian orthogonal expansion to find a scale with small interaction energy. We apply the first principle to the cutting error and use an orthogonal decomposition over independent coordinates to identify a suitable scale. The resulting variance estimate allows a smaller-order $L^2$ excess over the boundary bound.

Stabilization gives geometric central limit theorems from local determination or convergence of add-one responses, together with moment assumptions \citep{PenroseYukich2001Graphs,Penrose2005MultivariateSpatial,BaryshnikovYukich2005,Trinh2019}. For minimum spanning trees, the connectivity structure provides the required local approximation \citep{KestenLee1996MST,ChatterjeeSen2017MST}. The two-scale method of \citet{LachiezeReyPeccatiYang2022} compares full-domain and local add-one costs and obtains quantitative normal approximation when their discrepancy decays. Our criterion permits a block replacement discrepancy of boundary order and uses comparisons across scales to make the centered cutting error vanish.

The effect of a local modification on a distant response also enters perturbative normal approximation \citep{Chatterjee2008NormalApproximation}. Local-dependence methods apply when a suitable local description can be established \citep{AvramBertsimas1993}. For planar matching, \citet{Feng2026MatchingVariance} obtains bounded variance using parity and neighboring-cell arguments. Our conditional comparison addresses the further question needed for independent-cell approximation: whether the centered cutting error vanishes.

For fixed sample sizes, we pass from conditioning on a count window to conditioning on one count by comparing cost laws at nearby sample sizes. The density-dependent mean coefficient and the correction between Poisson and fixed-count covariances also arise in stabilization theory \citep{BaryshnikovYukich2005,Penrose2007}. In the planar degree-one setting, their form predicts a density-independent Poisson cost covariance. We derive the formulas from centered cutting and the count comparison, and extend them to the sampling laws above, including arbitrary singular components.

\subsection{Organization and notation}

Section~\ref{sec:method} states the boundary criterion, including the centered cutting estimate. Section~\ref{sec:probability} proves the multiscale variance bound and the Gaussian limits for Poisson and fixed sample sizes. Sections~\ref{sec:regular} and \ref{sec:parity} give geometric conditions sufficient for the required comparisons. Section~\ref{sec:fixed-law} treats the sampling laws $\mu=g\,dx+\tau$ above. Section~\ref{sec:nondegeneracy} characterizes the covariance null space, and Section~\ref{sec:applications} verifies the three optimization applications. Section~\ref{sec:scope} compares the probabilistic and deterministic assumptions and explains why the conditional mean must be controlled. The auxiliary constructions used in the parity moment estimates are proved in the appendix.

Output dimensions are fixed. We use Euclidean vector norms and set $\|X\|_p=(\E|X|^p)^{1/p}$. We write $\Law(X)$ for the law of $X$ and $\xrightarrow{W_2}$ for quadratic Wasserstein convergence; Gaussian limits may be singular. Labels remain distinct at coincident locations unless repetition invariance is proved. A fixed consistent Borel rule assigns partition faces. Constants may depend on the fixed costs and parameters, but not on sample size or intensity.

\paragraph{Statement of AI use.}
GPT-6 Pro was used to assist with the development of mathematical arguments and the preparation of this manuscript. The authors take full responsibility for the content of the paper.

\section{The boundary criterion}\label{sec:method}

Put $D=[0,1]^2$. For each closed axis-parallel square $Q$, let $H_Q$ be a measurable, finite-valued, permutation-invariant function from finite labeled configurations in $Q$ to $\R^r$. Labels at the same location are retained. We assume that
\begin{equation}\label{eq:similarity}
 H_{x+bD}(x+bS)=bH_D(S),\qquad x\in\R^2,\quad b>0,\quad S\subset D.
\end{equation}
The square is part of the argument; in particular, $H_Q$ may incorporate a boundary convention satisfying \eqref{eq:similarity}. Assign grid faces by a fixed Borel rule.

Let $\mathcal P_\lambda$ be a homogeneous Poisson process of intensity $\lambda$ on $D$, and write
\[
 P_\lambda=H_D(\mathcal P_\lambda),\qquad N_\lambda=|\mathcal P_\lambda|.
\]
For the partition $\mathcal C_k$ of $D$ into $k^2$ congruent squares, define
\begin{equation}\label{eq:cutting-remainder}
 R_{\lambda,k}=P_\lambda-\sum_{Q\in\mathcal C_k}H_Q(\mathcal P_\lambda\cap Q).
\end{equation}
If $B\subset Q\in\mathcal C_k$ is a square of side $b$, replace the process in $B$ by an independent Poisson process of the same intensity, leaving its exterior unchanged. Denote the resulting configuration by $\mathcal P_\lambda^{(B)}$, and put
\begin{equation}\label{eq:mixed-response}
 \Delta_{B,Q}H
 =H_D(\mathcal P_\lambda)-H_D(\mathcal P_\lambda^{(B)})
  -H_Q(\mathcal P_\lambda\cap Q)+H_Q(\mathcal P_\lambda^{(B)}\cap Q).
\end{equation}
Thus $\Delta_{B,Q}H=R_{\lambda,k}-R_{\lambda,k}^{(B)}$. The comparison is between the effects of the same inner replacement in the full square and in its containing cell.

\begin{theorem}[Boundary criterion]\label{thm:boundary}
Suppose that there are $C>0$ and $\eta:[1,\infty)\to[0,\infty)$, with $\eta(u)\to0$, such that
\begin{equation}\label{eq:boundary-excess}
 \E\bigl[(|\Delta_{B,Q}H_a|-Cb)_+^2\bigr]
 \le b^2\eta(\lambda b^2),\qquad \lambda b^2\ge1,
\end{equation}
uniformly over coordinates $a\le r$, intensities $\lambda$, integers $k\ge1$, cells $Q\in\mathcal C_k$, and squares $B\subset Q$ of side $b$. Then $P_\lambda\in L^2$ for every $\lambda>0$, and the following conclusions hold.
\begin{enumerate}[label=(\roman*),leftmargin=*]
\item There are $s_0<\infty$ and a nonnegative function $\rho(s)\to0$ such that
\begin{equation}\label{eq:boundary-centered-cut}
 \|R_{\lambda,k}-\E R_{\lambda,k}\|_2
 \le\rho(\lambda/k^2),\qquad \lambda/k^2\ge s_0,
\end{equation}
uniformly in $k$. Thus the centered total cost is approximated in $L^2$ by the sum of the centered cell costs.
\item There are a positive semidefinite matrix $\Sigma_P$ and a vector $\beta\in\R^r$ such that
\begin{equation}\label{eq:boundary-joint-limit}
 \left(P_\lambda-\E P_\lambda,\frac{N_\lambda-\lambda}{\sqrt\lambda}\right)
 \xrightarrow[\lambda\to\infty]{W_2}
 \mathcal N_{r+1}\!\left(0,
  \begin{pmatrix}\Sigma_P&\beta/2\\ \beta^\top/2&1\end{pmatrix}\right),
 \qquad \frac{\E P_\lambda}{\sqrt\lambda}\longrightarrow\beta.
\end{equation}
All joint covariance entries converge. The Gaussian law may be singular.
\end{enumerate}
\end{theorem}

\begin{proof}
Corollary~\ref{cor:average}, based on the multiscale variance bound, proves \eqref{eq:boundary-centered-cut} and square integrability of the remainder for sufficiently large $\lambda/k^2$. Theorem~\ref{thm:cut-clt} then gives square integrability of $P_\lambda$ at every positive intensity and the limits in \eqref{eq:boundary-joint-limit}.
\end{proof}

Square integrability, the covariance limit, and the first-order mean asymptotic all follow from the boundary comparison. It is enough to impose \eqref{eq:boundary-excess} on the nested squares used in Corollary~\ref{cor:average}, which also gives a version with the excess averaged over blocks and scales.

\subsection{Fixed sample sizes}

Let $X_1,X_2,\ldots$ be iid uniform points in $D$, and set $H_n=H_D(\{X_1,\ldots,X_n\})$. For probability laws on a Euclidean space, let $\dBL$ be the supremum of the differences of integrals over real functions bounded by one and with Lipschitz constant at most one.

\begin{corollary}\label{cor:boundary-fixed}
Under the assumptions of Theorem~\ref{thm:boundary}, suppose that $\{|H_n-\E H_n|^2:n\ge1\}$ is uniformly integrable and that, for some integer $q\ge1$ and $K<\infty$,
\begin{equation}\label{eq:step-count}
 \dBL\bigl(\Law(H_{j+q}),\Law(H_j)\bigr)\le K(j+1)^{-1/2},\qquad j\ge0.
\end{equation}
Then
\begin{equation}\label{eq:boundary-fixed-limit}
 H_n-\E H_n\xrightarrow[n\to\infty]{W_2}
 \mathcal N_r\!\left(0,\Sigma_P-\frac14\beta\beta^\top\right),
 \qquad \E H_n-\E P_n\longrightarrow0.
\end{equation}
The same limit holds with center $\E P_n$, and the covariance matrices converge.
\end{corollary}

\begin{proof}
Theorem~\ref{thm:boundary} gives the centered cutting estimate required by Corollary~\ref{cor:fixed-uniform}. Also, $H_n\in L^2$ because $P_\lambda\in L^2$ and $\Pp(N_\lambda=n)>0$. That corollary now gives the result.
\end{proof}

Lemma~\ref{lem:count-window} gives a weaker window condition in place of \eqref{eq:step-count}. Lemma~\ref{lem:exact-count} transfers sufficiently high Poisson moments to exact-count moments. Neither statement requires monotonicity in the sample size.

\subsection{Conditional replacements}

Constraints such as matching parity can make it useful to fix a state $\pi$ before replacing the points in a cell. Theorem~\ref{thm:multiscale} applies when the inputs are independent conditionally on $\pi$; it then bounds the expected conditional variance. For a square-integrable cutting remainder,
\begin{equation}\label{eq:variance-split}
 \E|R_{\lambda,k}-\E R_{\lambda,k}|^2
 =\E|R_{\lambda,k}-\E[R_{\lambda,k}\mid\pi]|^2
  +\E|\E[R_{\lambda,k}\mid\pi]-\E R_{\lambda,k}|^2.
\end{equation}
The multiscale argument controls the first term. The second is the variance of the conditional mean and needs a separate estimate. Once both terms vanish, Theorem~\ref{thm:cut-clt} and the fixed-count argument apply. For matching, $\pi$ consists of the cell-count parities; Section~\ref{sec:parity} bounds the conditional mean in this case. Example~\ref{ex:state-fluctuation} shows that this term need not vanish merely because the first one does.

\section{The probabilistic argument}\label{sec:probability}

We first turn replacement estimates at several scales into a variance bound. We then prove that a vanishing centered cutting error implies square integrability and a joint Gaussian limit for the centered cost and the standardized point count. The final step passes from Poisson to fixed sample sizes.

\subsection{A multiscale variance bound}

In the geometric application, a node will represent a square, its children will form a partition of that square, and its mass will be its area. We state the variance estimate for a general finite rooted forest. Let $I$ be its leaf set, with every leaf at depth $H\ge1$. Identify each node with its set of descendant leaves. Give the nodes deterministic masses $m(B)>0$, additive over children, with total root mass one, and put
\[
 \delta=\max_{A\text{ child of }B}\sqrt{\frac{m(A)}{m(B)}}.
\]
Suppose that $U=(U_i)_{i\in I}$ has a regular conditional product law given a state $\pi$. For a real $F=F(\pi,U)\in L^1$, set
\[
 \mathcal G_B=\sigma(\pi,U_{B^c}),\qquad
 X_B=F-\E[F\mid\mathcal G_B],\qquad
 \mathcal E_C=\sum_{B\text{ nonroot}}\E\bigl[(|X_B|-C\sqrt{m(B)})_+^2\bigr].
\]
Let $F^{(B,\pi)}$ be obtained by replacing $U_B$ with a conditionally independent copy at fixed $\pi$. Define $\mathcal E_C^{\rm res}$ by replacing $X_B$ in the preceding sum with $F-F^{(B,\pi)}$.

The proof uses the orthogonal decomposition for independent coordinates \citep{EfronStein1981,WehrAizenman1990Extensive}. As in the scale-selection argument of \citet{Chatterjee2018RFIMDecay}, a bound on the sum of the squared interaction norms identifies a suitable level. At that level, a quantitative normal approximation conflicts with the boundary bound if the variance is too large, following the principle of \citet{AizenmanWehr1990Rounding}.

\begin{theorem}[Multiscale replacement bound]\label{thm:multiscale}
For $C>0$, put
\[
 \varepsilon=\frac{2+4\mathcal E_C/C^2}{H},\qquad
 R=\min\{\delta^{-1},\varepsilon^{-1}\}.
\]
If $R\ge2$, then $F-\E[F\mid\pi]\in L^2$ and
\begin{equation}\label{eq:multiscale-bound}
 \E\Var(F\mid\pi)\le\frac{32C^2}{\log R},\qquad
 \Var F\le\frac{32C^2}{\log R}+\Var(\E[F\mid\pi]).
\end{equation}
The second inequality is understood in the extended sense. Both conclusions hold with $\mathcal E_C^{\rm res}$ in place of $\mathcal E_C$. In the unconditional resampling version, $F$ need only be measurable and finite valued: $\mathcal E_C^{\rm res}<\infty$ already implies $F\in L^2$.
\end{theorem}

\begin{proof}
Order the leaves and let $\mathcal F_i=\sigma(\pi,U_1,\ldots,U_i)$. Conditional independence gives
\[
 \E[F\mid\mathcal F_i]-\E[F\mid\mathcal F_{i-1}]
 =\E[X_{\{i\}}\mid\mathcal F_i].
\]
Finiteness of $\mathcal E_C$ makes each $X_{\{i\}}$ square integrable. Summing these finitely many martingale differences proves that $F-\E[F\mid\pi]\in L^2$. Subtract this conditional mean and divide by $C$. It remains to prove the assertion when $C=1$, $\E[F\mid\pi]=0$, and
\[
 d=\E\Var(F\mid\pi)=\E F^2.
\]

Let $P_i$ integrate out $U_i$ with its conditional marginal at fixed $\pi$. These are commuting orthogonal projections. Expanding the product of $P_i+(\operatorname{Id}-P_i)$ gives
\[
 F=\sum_{\varnothing\ne S\subset I}f_S,
 \qquad
 f_S=\prod_{i\in S}(\operatorname{Id}-P_i)\prod_{i\notin S}P_iF.
\]
The component indexed by the empty set is zero. For two distinct subsets, integrating in any index belonging to just one of them shows that the corresponding components are conditionally orthogonal. Thus, with $v_S=\E f_S^2$,
\[
 \sum_{\varnothing\ne S\subset I}v_S=d.
\]
For a child $A$ of $B$, write
\[
 Y_{B,A}=\E[F\mid\mathcal G_B,U_A]-\E[F\mid\mathcal G_B],\qquad
 D_B=\E X_B^2,\qquad
 e_B=\E\left(X_B-\sum_{A\text{ child of }B}Y_{B,A}\right)^2.
\]
If $c_B(S)$ counts the children of $B$ met by $S$, then
\begin{equation}\label{eq:interaction-budget-local}
 D_B=\sum_{S:S\cap B\ne\varnothing}v_S,\qquad
 e_B=\sum_{S:c_B(S)\ge2}v_S
 \le\sum_{A\text{ child of }B}D_A-D_B.
\end{equation}
Indeed, $Y_{B,A}$ retains exactly the components whose nonempty intersection with $B$ is contained in $A$; the rightmost difference gives $v_S$ the weight $(c_B(S)-1)_+$. Summing over internal nodes telescopes, and every nonempty $S$ meets a root. Consequently
\begin{equation}\label{eq:interaction-budget}
 \sum_{B\text{ internal}}e_B
 \le\sum_{i\text{ leaf}}D_i-\sum_{B\text{ root}}D_B
 \le\sum_{i\text{ leaf}}D_i-d.
\end{equation}

Let $t_j$ be the excess summed over nodes at depth $j$, and let $e_j$ be the sum of $e_B$ at an internal depth $j$. Put $E=\sum_{j=1}^Ht_j$. Since
\[
 x^2\le2a^2+2(|x|-a)_+^2,
\]
we have $\sum_{i\text{ leaf}}D_i\le2+2t_H$. If $d\le\varepsilon$, the conclusion follows from $\varepsilon\le R^{-1}<32/\log R$. Otherwise \eqref{eq:interaction-budget} yields
\[
 \sum_{j=1}^{H-1}(e_j+t_j+t_{j+1})
 \le2+4E-d<(H-1)\varepsilon.
\]
Here $R\ge2$ implies $H\ge4$. There is therefore a deterministic $j\in\{1,\ldots,H-1\}$ such that
\begin{equation}\label{eq:selected-depth}
 e_j+t_j+t_{j+1}\le\varepsilon.
\end{equation}
The sums in \eqref{eq:selected-depth} already average over $\pi$, so the selected depth is deterministic.

For each node $B$ at this depth and each child $A$, define
\[
 K_A=\bigl(\E[\,\cdot\mid\mathcal G_B,U_A]-\E[\,\cdot\mid\mathcal G_B]\bigr)
       \operatorname{clip}(X_A,\sqrt{m(A)}),
 \qquad
 \operatorname{clip}(x,a)=\max\{-a,\min\{x,a\}\}.
\]
Given $\mathcal G_B$, the $K_A$ are independent and centered. They satisfy $|K_A|\le2\sqrt{m(A)}$. Their conditional ranges have length at most $2\sqrt{m(A)}$, so
\[
 \Var(K_A\mid\mathcal G_B)\le m(A).
\]
The same projection applied to $X_A$ gives $Y_{B,A}$: conditional independence implies
\[
 \E[\E[F\mid\mathcal G_A]\mid\mathcal G_B,U_A]=\E[F\mid\mathcal G_B].
\]
The projection is an $L^2$ contraction. The errors $Y_{B,A}-K_A$ for different children are mutually orthogonal, and each is orthogonal to the residual defining $e_B$. Hence
\begin{equation}\label{eq:clipped-approximation}
 \E\left(X_B-\sum_AK_A\right)^2
 \le e_B+\sum_A\E\bigl[(|X_A|-\sqrt{m(A)})_+^2\bigr].
\end{equation}

Select a depth-$j$ node $B$, independently of all inputs, with probability $m(B)$. In the rest of this argument set
\[
 U_*=\frac{X_B}{\sqrt{m(B)}},\qquad
 V_* =\frac{\sum_AK_A}{\sqrt{m(B)}},\qquad
 W_* =\operatorname{clip}(U_*,1).
\]
Equations \eqref{eq:selected-depth}--\eqref{eq:clipped-approximation} give
\begin{align*}
 \E U_*^2&=\sum_{B\text{ at depth }j}D_B\ge d,\qquad
 \E(U_*-V_*)^2\le e_j+t_{j+1},\\
 \E(U_*-W_*)^2&=t_j,\qquad
 \E(V_*-W_*)^2\le2\varepsilon.
\end{align*}
Let $v$ be the variance of $V_*$ conditionally on the selected node, its exterior, and $\pi$. Then $v\le1$. If $d\ge16\varepsilon$, the $L^2$ triangle inequality gives
\[
 \E v=\E V_*^2\ge(\sqrt d-\sqrt\varepsilon)^2\ge d/2.
\]
It follows that $d\le2$ and $\Pp(v\ge d/4)\ge d/4$. On this event the conditionally independent summands of $V_*$ have absolute value at most $2\delta$, and their absolute third moments sum to at most $2\delta v$. The Berry--Esseen inequality with constant $10$ \citep[Theorem~3.6]{ChenGoldsteinShao2011} therefore gives
\begin{equation}\label{eq:gaussian-tail-contradiction}
 \frac d4\left[1-\Phi(4/\sqrt d)-\frac{40\delta}{\sqrt d}\right]_+
 \le\Pp(V_*>2)\le\E(V_*-W_*)^2\le2\varepsilon,
\end{equation}
where $\Phi$ is the standard normal distribution function.

Suppose now that $d\ge32/\log R$. Since $\varepsilon\le R^{-1}$, this implies $d\ge16\varepsilon$. The resulting bound $d\le2$ gives $\log R\ge16$. Integrating the standard normal density from $x$ to $x+1/x$, with $x=4/\sqrt d$, yields
\[
 1-\Phi(4/\sqrt d)\ge c_0\sqrt d\,e^{-8/d},
 \qquad c_0=\frac{e^{-3/2}}{4\sqrt{2\pi}}>\frac1{60}.
\]
The ratio of the Berry--Esseen error to this Gaussian tail lower bound is at most
\[
 75(\log R)R^{-3/4}\le1200e^{-12}<\frac12.
\]
Thus the left side of \eqref{eq:gaussian-tail-contradiction} is at least
\[
 \frac{c_0}{8}\left(\frac{32}{\log R}\right)^{3/2}R^{-1/4}
 >\frac2R\ge2\varepsilon.
\]
The strict inequality holds at $\log R=16$, and its left side divided by $2/R$ increases thereafter. This contradiction proves the first assertion of \eqref{eq:multiscale-bound}. Restore $C$ and use total variance for the second.

For the resampling version, conditional Jensen gives
\[
 X_B=\E[F-F^{(B,\pi)}\mid\pi,U],\qquad
 \mathcal E_C\le\mathcal E_C^{\rm res}.
\]
Finally, suppose $\pi$ is constant and $F$ is merely finite and measurable. Replace the coordinates of $U$ successively by those of an independent copy $U'$. Every increment has the law of a one-leaf replacement difference and belongs to $L^2$ when $\mathcal E_C^{\rm res}<\infty$. Thus $F(U)-F(U')\in L^2$. Fubini gives a fixed $u'$ with $F(u')$ finite and $F(U)-F(u')\in L^2$, proving $F\in L^2$ and justifying the preceding argument.
\end{proof}

\subsection{From centered cutting to a joint Gaussian limit}

Return to the square costs in Section~\ref{sec:method}. For laws on a finite-dimensional Euclidean space, $W_2$ is the infimum of the $L^2$ distance over couplings. We use the equivalence between $W_2$ convergence and weak convergence together with convergence of second moments \citep[Theorem~5.11]{Santambrogio2015Transport}.

\begin{theorem}[Gaussian limit from centered cutting]\label{thm:cut-clt}
Suppose that there are $s_0<\infty$ and a nonnegative function $\rho(s)\to0$ such that, whenever $s=\lambda/k^2\ge s_0$, the remainder $R_{\lambda,k}$ is integrable and
\begin{equation}\label{eq:centered-cut}
 \|R_{\lambda,k}-\E R_{\lambda,k}\|_2\le\rho(s).
\end{equation}
Then $P_\lambda\in L^2$ for every $\lambda>0$. Set
\[
 m(\lambda)=\E P_\lambda,\qquad
 V_\lambda=\left(P_\lambda-m(\lambda),\frac{N_\lambda-\lambda}{\sqrt\lambda}\right),
 \qquad \Gamma(s)=\Cov(V_s).
\]
There are $\Sigma_P$ and $\beta$ such that
\begin{equation}\label{eq:cut-joint-limit}
 V_\lambda\xrightarrow[\lambda\to\infty]{W_2}\mathcal N_{r+1}(0,\Gamma),
 \qquad
 \Gamma=\begin{pmatrix}\Sigma_P&\beta/2\\\beta^\top/2&1\end{pmatrix},
 \qquad \frac{m(\lambda)}{\sqrt\lambda}\longrightarrow\beta.
\end{equation}
The covariance matrices converge, and the limiting Gaussian may be singular. For all sufficiently large $s$,
\begin{equation}\label{eq:covariance-rate}
 W_2\bigl(\mathcal N_{r+1}(0,\Gamma(s)),\mathcal N_{r+1}(0,\Gamma)\bigr)\le\rho(s),
 \qquad \|\Gamma(s)-\Gamma\|_{\max}\le K\rho(s),
\end{equation}
where $K$ may depend on the limiting covariance.
\end{theorem}

\begin{proof}
We first obtain square integrability. Fix $\lambda$ with $\lambda/4\ge s_0$, and put
\[
 h(S)=H_D(S)-H_D(\varnothing).
\]
For a measurable function $u$ on configurations with $u(\varnothing)=0$, let $T_ju$ be $2^{-j}$ times the sum of $u$ evaluated on the $4^j$ rescaled dyadic cells at depth $j$. Thus $T_0=\operatorname{Id}$ and $T_jT_1=T_{j+1}$. By \eqref{eq:similarity},
\[
 g=h-T_1h=R_{\lambda,2}+H_D(\varnothing)
\]
when evaluated at $\mathcal P_\lambda$. Hence $g\in L^2(\mathcal P_\lambda)$ and $g(\varnothing)=0$.

Write $\mu=4^{-j}\lambda$, and denote expectation at intensity $\mu$ by $\E_\mu$. On nonempty configurations the likelihood ratio of the intensity-$\mu$ law relative to the intensity-$\lambda$ law is
\[
 e^{\lambda-\mu}4^{-jN}\le e^\lambda4^{-j}.
\]
Also
\[
 |\E_\mu g|^2\le\Pp_\mu(N>0)\E_\mu|g|^2\le\mu\E_\mu|g|^2.
\]
Independence of the dyadic cells therefore gives
\begin{equation}\label{eq:automatic-l2-series}
 \E_\lambda|T_jg|^2
 =\E_\mu|g|^2+(4^j-1)|\E_\mu g|^2
 \le(1+\lambda)e^\lambda4^{-j}\E_\lambda|g|^2.
\end{equation}
Thus $\sum_{j\ge0}T_jg$ converges in $L^2$. Telescoping gives the pointwise identity
\[
 h=\sum_{j=0}^{m-1}T_jg+T_mh.
\]
Conditionally on $N_\lambda=n$, the probability that two points lie in the same depth-$m$ cell is at most $n^2 4^{-m}$. Given distinct occupied cell labels, the $n$ rescaled locations are independent and uniform on $D$. On that event $T_mh$ has the law of
\[
 2^{-m}\sum_{i=1}^nh(\{U_i\}),
\]
where the $U_i$ are iid uniform. For each fixed $n$, this tends to zero in probability because all summands are finite valued. Averaging over $N_\lambda$ gives $T_mh\to0$ in probability. It follows that $h$ equals the $L^2$ limit of the series. This proves square integrability at every sufficiently large intensity. A smaller positive intensity has a bounded likelihood ratio relative to a larger one, so $P_\lambda\in L^2$ for all $\lambda>0$.

Next, $u\mapsto\Law(V_u)$ is continuous in $W_2$ at each $s>0$. On a compact interval of intensities, dominate all laws by a Poisson law of a strictly larger intensity. The likelihood ratios are uniformly bounded and converge pointwise as $u\to s$. Dominated convergence gives weak convergence and convergence of the first and second moments of $(P_u,N_u)$. Centering and rescaling prove the asserted continuity.

Put $\bar\rho(s)=\sup_{u\ge s}\rho(u)$, which is finite for large $s$ and tends to zero. Fix such an $s$, let $\lambda\to\infty$, and choose
\[
 k=\left\lfloor\sqrt{\lambda/s}\right\rfloor,\qquad u=\lambda/k^2.
\]
Then $u\ge s$ and $u\to s$. Independence, homogeneity and additivity of the count give
\begin{equation}\label{eq:joint-cell-decomposition}
 V_\lambda=\frac1k\sum_{j=1}^{k^2}V_u^{(j)}
                  +(R_{\lambda,k}-\E R_{\lambda,k},0),
\end{equation}
where the $V_u^{(j)}$ are independent copies of $V_u$. Couple each copy to a copy of $V_s$, using independent optimal couplings. The differences are centered, so the squared $L^2$ distance between the two normalized sums equals $W_2(V_u,V_s)^2$, which tends to zero. At fixed $s$, the iid finite-variance central limit theorem holds in $W_2$, since the covariance of the normalized sum is exactly $\Gamma(s)$. Consequently
\begin{equation}\label{eq:gaussian-approx-fixed-scale}
 \limsup_{\lambda\to\infty}
 W_2\bigl(V_\lambda,\mathcal N_{r+1}(0,\Gamma(s))\bigr)\le\bar\rho(s).
\end{equation}
The triangle inequality, applied through the same $V_\lambda$, gives
\[
 W_2\bigl(\mathcal N_{r+1}(0,\Gamma(s)),\mathcal N_{r+1}(0,\Gamma(t))\bigr)
 \le\bar\rho(s)+\bar\rho(t).
\]
These laws form a Cauchy family in the complete $W_2$ space. Their covariance entries converge, and their characteristic functions identify the limit as a centered Gaussian. Letting $s\to\infty$ in \eqref{eq:gaussian-approx-fixed-scale} proves the joint $W_2$ limit. This independent-cell approximation is a joint $W_2$ version of the scalar criterion in \citet[Lemma~2.2]{Trinh2019}, with the approximation error given by \eqref{eq:centered-cut}.

For fixed sufficiently large $s$, take $\lambda=k^2s$ exactly in \eqref{eq:joint-cell-decomposition} and let $k\to\infty$. This proves the first inequality in \eqref{eq:covariance-rate}. For centered coupled vectors $X,Y$, entrywise Cauchy--Schwarz gives
\[
 \|\Cov X-\Cov Y\|_{\max}
 \le(\|X\|_2+\|Y\|_2)\|X-Y\|_2.
\]
The Gaussian second moments are bounded for large $s$, giving the second inequality.

Finally, the Poisson series for $m$ can be differentiated locally. Integrability at a larger intensity provides a summable bound: the exponential decay of the likelihood ratio dominates the factor $N$ introduced by differentiation. We obtain
\begin{equation}\label{eq:poisson-mean-derivative}
 m'(s)=s^{-1}\Cov(P_s,N_s)=s^{-1/2}\Gamma(s)_{1:r,r+1}.
\end{equation}
Let $d=\lim_{s\to\infty}\Gamma(s)_{1:r,r+1}$. For fixed $M>0$,
\[
 \frac1{\sqrt\lambda}
 \left|\int_M^\lambda s^{-1/2}\bigl(\Gamma(s)_{1:r,r+1}-d\bigr)\,ds\right|
 \le2\sup_{s\ge M}|\Gamma(s)_{1:r,r+1}-d|.
\]
Let $\lambda\to\infty$ and then $M\to\infty$. This gives $m(\lambda)/\sqrt\lambda\to2d$. Set $\beta=2d$; the count coordinate has variance one at every intensity, completing \eqref{eq:cut-joint-limit}.
\end{proof}

\subsection{Averaging the replacement excess}

For $s=\lambda/k^2$ sufficiently large, put
\[
 h=\left\lfloor\frac12\sqrt{\log_2s}\right\rfloor.
\]
Starting from the $k^2$ root cells, divide each square into $2^{2h}$ equal children at each of $h$ generations. Let $\mathcal B_j$ be the depth-$j$ squares and let $Q(B)$ be the root containing $B$. Define
\begin{equation}\label{eq:averaged-excess}
 e_{\lambda,k}(C)=\max_{1\le a\le r}\frac1h\sum_{j=1}^h\sum_{B\in\mathcal B_j}
 \E\bigl[(|\Delta_{B,Q(B)}H_a|-C\operatorname{side}(B))_+^2\bigr].
\end{equation}

\begin{corollary}[Average-scale criterion]\label{cor:average}
Suppose that $e_{\lambda,k}(C)\to0$ uniformly as $\lambda/k^2\to\infty$, for some fixed $C>0$. Then the centered cutting estimate \eqref{eq:centered-cut} holds, and Theorem~\ref{thm:cut-clt} gives the joint Gaussian limit. In particular, \eqref{eq:boundary-excess} suffices. Whenever $\mathfrak R_{\lambda,k}\ge2$,
\begin{equation}\label{eq:averaged-cut-bound}
 \E|R_{\lambda,k}-\E R_{\lambda,k}|^2
 \le\frac{32rC^2}{\log\mathfrak R_{\lambda,k}},\qquad
 \mathfrak R_{\lambda,k}
 =\min\left\{2^h,\left(\frac2h+\frac{4e_{\lambda,k}(C)}{C^2}\right)^{-1}\right\}.
\end{equation}
If $e_{\lambda,k}(C)=O(h^{-1})$ uniformly, the cutting variance is $O((\log\log s)^{-1})$, and the homogeneous covariance error in \eqref{eq:covariance-rate} is $O((\log\log s)^{-1/2})$.
\end{corollary}

\begin{proof}
The leaf configurations are independent. Give each node its area as mass, so the root masses sum to one and $\delta=2^{-h}$. Replacing $B$ changes $R_{\lambda,k}$ by $\Delta_{B,Q(B)}H$. The replacement excess for each coordinate is at most $he_{\lambda,k}(C)$. The unconditional resampling form of Theorem~\ref{thm:multiscale} gives both square integrability of the remainder and \eqref{eq:averaged-cut-bound}.

For every node, with $b=\operatorname{side}(B)$,
\[
 \lambda b^2\ge s2^{-2h^2}\ge\sqrt s.
\]
Since areas sum to one at every depth, \eqref{eq:boundary-excess} implies
\[
 e_{\lambda,k}(C)\le\sup_{v\ge\sqrt s}\eta(v)\longrightarrow0.
\]
Under the averaged excess assumption, $h\to\infty$ and $\mathfrak R_{\lambda,k}\to\infty$ uniformly. Taking the square root of the bound in \eqref{eq:averaged-cut-bound}, and then its supremum over all larger cell intensities, gives a function $\rho$ satisfying \eqref{eq:centered-cut}. If $e_{\lambda,k}(C)=O(h^{-1})$, then $\mathfrak R_{\lambda,k}\ge ch$ for large $s$. Since $h$ has order $\sqrt{\log s}$, the variance rate follows, and \eqref{eq:covariance-rate} gives the covariance rate.
\end{proof}

The same argument applies when the leaf inputs are independent conditionally on $\pi$. Use replacements at fixed $\pi$ in \eqref{eq:averaged-excess}, and average over both $\pi$ and the conditional inputs. Denote the resulting quantities by $e_{\lambda,k}^{\pi}(C)$ and $\mathfrak R_{\lambda,k}^{\pi}$, and assume $R_{\lambda,k}\in L^1$. Theorem~\ref{thm:multiscale} gives, whenever $\mathfrak R_{\lambda,k}^{\pi}\ge2$,
\begin{equation}\label{eq:conditional-cut}
 \E|R_{\lambda,k}-\E R_{\lambda,k}|^2
 \le\frac{32rC^2}{\log\mathfrak R_{\lambda,k}^{\pi}}
       +\E|\E[R_{\lambda,k}\mid\pi]-\E R_{\lambda,k}|^2.
\end{equation}
The inequality remains valid if the right side is infinite. If the variance of the conditional mean vanishes uniformly and $e_{\lambda,k}^{\pi}(C)\to0$, then \eqref{eq:centered-cut} follows and Theorem~\ref{thm:cut-clt} applies.

\subsection{From a count window to an exact count}

The event that a Poisson variable equals its mean has probability tending to zero, so the joint limit alone does not justify conditioning on that event. We first condition on a window of width proportional to $\sqrt n$, then compare the cost laws at the individual counts in that window. The Gaussian regression term is the one appearing in conditional limit theorems; compare \citet[Theorems~2.2--2.3]{Janson2007Conditioning}. Here continuity across the count window replaces stochastic monotonicity.

Let $N_\ell\sim\operatorname{Poi}(n_\ell)$, where the integers $n_\ell\to\infty$, and let $P_\ell\in L^1(\R^r)$ be associated costs. Write
\[
 b_\ell=\E P_\ell,\qquad
 \nu_{j,\ell}=\Law(P_\ell\mid N_\ell=j),\qquad
 \Law(Y_\ell)=\nu_{n_\ell,\ell}.
\]
For a fixed integer $q\ge1$, set
\[
 D_\ell=\{N_\ell\equiv n_\ell\pmod q\},\qquad
 Z_\ell=\frac{N_\ell-n_\ell}{\sqrt{n_\ell}}.
\]

\begin{lemma}[Count-window comparison]\label{lem:count-window}
Suppose that, conditionally on $D_\ell$, the pair $(P_\ell-b_\ell,Z_\ell)$ converges weakly to a centered Gaussian pair $(G,Z)$ with covariance
\[
 \begin{pmatrix}A&d\\d^\top&1\end{pmatrix}.
\]
Let $M_{\ell,\delta}$ be the conditional law of $P_\ell-b_\ell$ given $D_\ell$ and $|Z_\ell|\le\delta$. Then
\begin{equation}\label{eq:count-window-comparison}
 Y_\ell-b_\ell\ \Longrightarrow\ G-dZ
 \quad\Longleftrightarrow\quad
 \lim_{\delta\downarrow0}\limsup_{\ell\to\infty}
 \dBL\bigl(\Law(Y_\ell-b_\ell),M_{\ell,\delta}\bigr)=0.
\end{equation}
A sufficient condition is
\begin{equation}\label{eq:local-count-continuity}
 \lim_{\delta\downarrow0}\limsup_{\ell\to\infty}
 \sup_{\substack{|j-n_\ell|\le\delta\sqrt{n_\ell}\\j\equiv n_\ell\pmod q}}
 \dBL(\nu_{j,\ell},\nu_{n_\ell,\ell})=0.
\end{equation}
In particular, the neighboring-count bounds
\[
 \dBL(\nu_{j+q,\ell},\nu_{j,\ell})\le K(j+1)^{-1/2}
\]
uniformly throughout these windows suffice. If the weak limit in \eqref{eq:count-window-comparison} holds and $\{|Y_\ell-\E Y_\ell|^2\}_\ell$ is uniformly integrable, then both centerings by $b_\ell$ and by $\E Y_\ell$ converge in $W_2$ to $\mathcal N_r(0,A-dd^\top)$. Moreover,
\[
 \E Y_\ell-b_\ell\longrightarrow0,
\]
and the covariance matrices converge.
\end{lemma}

\begin{proof}
For fixed $\delta>0$, the limiting window has positive probability and its boundary has probability zero. Hence
\[
 M_{\ell,\delta}\Longrightarrow\Law(G\mid |Z|\le\delta).
\]
The Gaussian vector $G_0=G-dZ$ is independent of $Z$ and has covariance $A-dd^\top$. It can be coupled to the limiting window law with distance at most $|d|\delta$. Thus this window law tends to $\Law(G_0)$ as $\delta\downarrow0$. The triangle inequality in both directions proves \eqref{eq:count-window-comparison}.

The finite-window law is a mixture over the counts in \eqref{eq:local-count-continuity}. Convexity and translation invariance of $\dBL$ bound its distance to $\Law(Y_\ell-b_\ell)$ by the displayed supremum. The step bound telescopes over $|j-n_\ell|/q$ increments. All intervening counts are at least $n_\ell/2$ for large $\ell$, so this gives a bound $C\delta$, uniformly in the window.

Uniform integrability makes $Y_\ell-\E Y_\ell$ tight, and weak convergence makes $Y_\ell-b_\ell$ tight. Their deterministic difference is bounded: choose a ball containing each with probability greater than $3/4$ and intersect the two events. Bounded deterministic shifts preserve uniform integrability of squared norms. It follows that $\{|Y_\ell-b_\ell|^2\}_\ell$ is uniformly integrable, so its first and second moments converge to those of $G_0$. This gives $\E Y_\ell-b_\ell\to0$, the covariance convergence and $W_2$ convergence under both centerings.
\end{proof}

\begin{lemma}[Poisson tails at an exact count]\label{lem:exact-count}
Let $N\sim\operatorname{Poi}(n)$, where $n\ge1$ is an integer, let $P$ be an associated cost in $\R^r$, and let $Y$ have its conditional law given $N=n$. Suppose that consecutive counts in the residue class of $n$, both within distance $\sqrt n$ of $n$, have cost laws at $\dBL$ distance at most $K/\sqrt n$ per step of length $q$. Then, for every $b\in\R^r$ and $t\ge0$,
\begin{equation}\label{eq:exact-count-tail}
 \Pp(|Y-b|>t+1)\le C_{q,K}\sqrt{\Pp(|P-b|>t)}.
\end{equation}
Consequently, a uniform $p$th moment bound for the Poisson mixture about $b$ gives uniform exact-count moments about $b$ of every order $p'<p/2$. For $p>4$, the squared exact-count norms about $b$ are uniformly integrable.
\end{lemma}

\begin{proof}
Put $S(t)=\Pp(|P-b|>t)$. For $n^{-1/2}\le\delta\le1$, let
\[
 W_\delta=\{|N-n|\le\delta\sqrt n,\ N\equiv n\pmod q\}.
\]
Stirling's formula gives a lower bound $c/\sqrt n$ for each Poisson mass in this window. There are at least $c_q\delta\sqrt n$ eligible integers, including $n$, so
\[
 \Pp(W_\delta)\ge c_q\delta.
\]
Use the bounded Lipschitz function $f(x)=\min\{1,(|x-b|-t)_+\}$. Its expectation at count $n$ is at least $\Pp(|Y-b|>t+1)$. Telescoping the step comparisons from any eligible count to $n$, and then averaging, yields
\[
 \Pp(|Y-b|>t+1)
 \le\frac{S(t)}{\Pp(W_\delta)}+C_{q,K}\delta
 \le C_{q,K}\left(\frac{S(t)}\delta+\delta\right).
\]
Take $\delta=\sqrt{S(t)}$ if it lies in $[n^{-1/2},1]$. If $0<S(t)<n^{-1}$, use the single atom $N=n$, whose mass is at least $c/\sqrt n$, to obtain
\[
 \Pp(|Y-b|>t+1)\le C\sqrt n\,S(t)\le C\sqrt{S(t)}.
\]
The same atom argument handles $S(t)=0$. This proves \eqref{eq:exact-count-tail}.

A uniform bound on $\E|P-b|^p$, together with Markov's inequality and \eqref{eq:exact-count-tail}, gives a uniform tail bound of order $u^{-p/2}$ for $|Y-b|$ at large $u$. Integrating the tail gives every moment of order $p'<p/2$. If $p>4$, choose $2<p'<p/2$.
\end{proof}

\begin{corollary}[Homogeneous fixed-sample limit]\label{cor:fixed-uniform}
Let $H_n$ be the cost of $n$ iid uniform points in $D$. Under the assumptions of Theorem~\ref{thm:cut-clt}, suppose that the count-window comparison on the right side of \eqref{eq:count-window-comparison} holds along $n_\ell=\ell$ for some fixed $q$, with $P_\ell=P_{n_\ell}$ and $\nu_{j,\ell}=\Law(H_j)$, and that $\{|H_n-\E H_n|^2\}_n$ is uniformly integrable. In particular, it is enough to assume \eqref{eq:local-count-continuity} or \eqref{eq:step-count} together with this uniform integrability. Then
\begin{equation}\label{eq:fixed-uniform-limit}
 H_n-\E H_n\xrightarrow[n\to\infty]{W_2}
 \mathcal N_r\!\left(0,\Sigma_P-\frac14\beta\beta^\top\right),
 \qquad \E H_n-m(n)\longrightarrow0.
\end{equation}
The same limit holds with center $m(n)$, and the covariance matrices converge. The moment hypothesis follows, for example, from \eqref{eq:step-count} and
\[
 \sup_{s\ge1}\|P_s-m(s)\|_p<\infty
\]
for some $p>4$.
\end{corollary}

\begin{proof}
We first show that the joint limit in \eqref{eq:cut-joint-limit} persists after conditioning on any count residue modulo $q$. For $q>1$, fix a nontrivial $q$th root of unity $z$ and $t\in\R^{r+1}$. In \eqref{eq:joint-cell-decomposition}, with $k$ fixed and $\lambda$ large, each centered cost--count cell term has $L^2$ norm at most $C/k$, because $\Gamma(\lambda/k^2)$ is bounded. Omit one cell term and the cutting remainder from the exponent. The count in the omitted cell remains in $z^{N_\lambda}$; it is independent of the other cells and has Poisson mean $\lambda/k^2$. Thus
\begin{equation}\label{eq:residue-fourier}
 \left|\E\left[z^{N_\lambda}e^{it\cdot V_\lambda}\right]\right|
 \le e^{-(\lambda/k^2)(1-\operatorname{Re}z)}
  +|t|\left(\frac Ck+\|R_{\lambda,k}-\E R_{\lambda,k}\|_2\right).
\end{equation}
Let $\lambda\to\infty$ at fixed $k$, and then let $k\to\infty$. The right side tends to zero. With $\omega=e^{2\pi i/q}$, Fourier inversion gives, uniformly in $a\in\{0,\ldots,q-1\}$,
\[
 \E\bigl[e^{it\cdot V_\lambda}\1_{\{N_\lambda\equiv a\pmod q\}}\bigr]
 =\frac1q\sum_{j=0}^{q-1}\omega^{-ja}
       \E\bigl[\omega^{jN_\lambda}e^{it\cdot V_\lambda}\bigr]
 \longrightarrow\frac1q e^{-t^\top\Gamma t/2}.
\]
At $t=0$ each residue has probability tending to $1/q$. Dividing by that probability gives the conditional limit. The convergence is uniform over the finitely many residues, so the selected residue may vary with $\lambda$. The case $q=1$ needs no conditioning argument.

Apply Lemma~\ref{lem:count-window} with $n_\ell=\ell$, $b_\ell=m(\ell)$ and $d=\beta/2$. This proves \eqref{eq:fixed-uniform-limit} and covariance convergence. Under the stated sufficient moment condition, Lemma~\ref{lem:exact-count} gives a uniform exact-count moment of some order $p'>2$ about $m(n)$. Jensen's inequality transfers this bound to centering by $\E H_n$, which gives the required uniform integrability.
\end{proof}

\section{Regular costs and boundary replacement}
\label{sec:regular}

The deterministic comparison in this section replaces a nonempty configuration inside a square by one boundary point. Applying it with two exterior configurations cancels the common interior cost and controls the effect of a block replacement on the cutting error. Empty blocks require a separate estimate. For that estimate, we assume that the absolutely continuous part of the sampling measure has a density positive almost everywhere; its singular part is unrestricted.

Throughout this section, $D=[0,1]^2$, $m$ denotes planar Lebesgue measure, and all squares have sides parallel to the coordinate axes. The notation $S\sqcup T$ joins two labeled configurations, keeping distinct labels even when their locations coincide. We begin with costs on sets and then extend them to such configurations.

\subsection{Consequences of the boundary comparison}

Tour length and the Held--Karp value satisfy the boundary comparison below, as proved in Section~\ref{sec:applications}.

\begin{definition}[Regular cost]\label{def:regular}
A real-valued function $H$ on finite subsets of $\R^2$ is called regular if
\[
 H(x+rS)=rH(S),\qquad x\in\R^2,\quad r>0,
\]
and there are a point $z_*\in\partial D$ and a finite constant $B_\partial$ such that
\begin{equation}\label{eq:regular-definition}
 \bigl|H(I\cup O)-H(I)-H(O\cup\{z_*\})\bigr|
 \le B_\partial
\end{equation}
for every nonempty finite $I\subset D^\circ$ and every finite $O\subset\R^2\setminus D$.
\end{definition}

Translation invariance and homogeneity force the empty and singleton costs to be zero. For a square of side $b$, the error in \eqref{eq:regular-definition} is at most $bB_\partial$. The absolute value is essential: the displacement and insertion estimates below use both directions of the comparison.

\begin{lemma}[Displacement, insertion, and arbitrary boundary points]
\label{lem:regular-insert}
A regular cost has a unique continuous, permutation-invariant extension to each labeled configuration space $(\R^2)^n$. Removing all but one label at a repeated location leaves the cost unchanged. The extensions satisfy
\begin{equation}\label{eq:regular-displacement}
 \bigl|H(x_1,\ldots,x_n)-H(y_1,\ldots,y_n)\bigr|
 \le K\sum_{i=1}^n|x_i-y_i|,
 \qquad K=2B_\partial.
\end{equation}
For every nonempty labeled configuration $S$,
\begin{equation}\label{eq:regular-insertion}
 \bigl|H(S\sqcup\{x\})-H(S)\bigr|
 \le K\operatorname{dist}(x,S).
\end{equation}
Moreover, if $Q$ is a square of side $b$, $z\in\partial Q$, $I$ is a nonempty labeled configuration in $Q$, and $O$ is a labeled configuration outside $Q^\circ$, then
\begin{equation}\label{eq:regular-full-boundary}
 \bigl|H(I\sqcup O)-H(I)-H(O\sqcup\{z\})\bigr|
 \le B_0b,
 \qquad B_0=(1+2\sqrt2)B_\partial.
\end{equation}
Boundary labels may be assigned to either configuration, including distinct labels at the same location.
\end{lemma}

\begin{proof}
Fix a finite background $O$ and two points $x,y\notin O$. If $x$ and $y$ lie in the interior of a square of side $b$ whose closure avoids $O$, apply \eqref{eq:regular-definition} with singleton interiors. Their interior costs are zero, and the common exterior term cancels. Hence
\[
 |H(O\cup\{x\})-H(O\cup\{y\})|\le 2B_\partial b.
\]
For sufficiently close $x,y$ away from $O$, we may let $b$ decrease to $\|x-y\|_\infty$. For arbitrary $x,y\notin O$, choose a polygonal path avoiding $O$ with length at most $|x-y|+\varepsilon$. Such a path can be obtained by choosing an intermediate point arbitrarily close to the midpoint of $xy$ and avoiding the finitely many forbidden lines through a point of $O$. The path is compact and has positive distance from $O$. Subdivide it into sufficiently short segments, apply the local estimate, and let $\varepsilon\downarrow0$. This proves the one-coordinate displacement estimate.

For two lists of distinct points, choose distinct auxiliary points arbitrarily close to the target points and disjoint from both lists. Move the labels first to the auxiliary points and then to their targets, one at a time. Summing the one-coordinate bounds and passing to the limit gives \eqref{eq:regular-displacement} on the dense set of distinct lists. The estimate therefore defines a unique continuous extension to $(\R^2)^n$. Permutation invariance, translation invariance, and homogeneity pass to this extension.

It remains to verify consistency across cardinalities. An $M$-tuple concentrated at one point has zero cost, by translation invariance and homogeneity. Fix a background $O$ with no label at $z$. In a square of side $b$ shrinking to $z$, place $M$ distinct interior labels and apply the boundary comparison. As $b\downarrow0$, the interior cost tends to zero, the prescribed boundary point tends to $z$, and continuity gives
\begin{equation}\label{eq:regular-repetition}
 H(O\sqcup\{z\}^{[M]})=H(O\sqcup\{z\}),
 \qquad M\ge1.
\end{equation}
Here $\{z\}^{[M]}$ denotes $M$ labels at $z$. Perturbing the background and using continuity removes the restriction on $O$. Thus adding further labels at an occupied site does not change the cost. Moving the inserted point to a nearest old label in \eqref{eq:regular-displacement} proves \eqref{eq:regular-insertion}.

For a square of side $b$, moving its prescribed boundary point to $z$ costs at most $K\sqrt2b$. This proves \eqref{eq:regular-full-boundary} first for strictly interior and strictly exterior labels. Approximate any boundary labels by distinct points on their designated sides and use continuity. The resulting error is
\[
 B_\partial b+K\sqrt2b=(1+2\sqrt2)B_\partial b,
\]
as claimed.
\end{proof}

Finite linear combinations of regular costs are regular. Indeed, use \eqref{eq:regular-full-boundary} at a common boundary point and sum the absolute errors.

\begin{lemma}[Short paths]\label{lem:short-path}
The $\ell\times\ell$ grid of cells of side $b/\ell$ has a Hamiltonian path, using only shared sides, that ends at any prescribed corner cell. Joining successive cell centers along this path gives a curve of length at most $b\ell$. If $\ell$ is even, any ordered pair of adjacent corner cells can be prescribed as its endpoints. Furthermore, any $n\ge1$ labeled sites in a square of side $b$ can be ordered so that
\begin{equation}\label{eq:short-site-path}
 \sum_{i=2}^n|x_i-x_{i-1}|\le 4b\sqrt n.
\end{equation}
Consequently, a regular cost with insertion constant $K$ satisfies
\begin{equation}\label{eq:regular-growth}
 |H(S)|\le 4Kb\sqrt{|S|}
\end{equation}
for every labeled configuration $S$ in that square.
\end{lemma}

\begin{proof}
Traverse the columns in turn, alternating between upward and downward motion. Starting in the lower-left cell, this path ends in the lower-right cell when $\ell$ is even. The symmetries of the square therefore give every ordered pair of adjacent corner cells in that case; for arbitrary $\ell$, they give any prescribed terminal corner cell. The path through the centers has $\ell^2-1$ edges of length $b/\ell$.

For the second assertion, take $\ell=\lceil\sqrt n\rceil$. At the center of each visited cell, make a round trip to every label in that cell. The total length is at most
\[
 b\ell+\frac{\sqrt2bn}{\ell}\le4b\sqrt n.
\]
List the sites in their order of visits and shortcut the intervening walk. This gives \eqref{eq:short-site-path}. Insert the labels in this order, starting with the zero singleton cost. Each insertion distance is at most the distance to the preceding label, so \eqref{eq:regular-insertion} gives \eqref{eq:regular-growth}.
\end{proof}

\subsection{Moments uniform over the sampling law}

Deleting one point changes the cost by at most its distance to the nearest remaining point, multiplied by a fixed constant. A packing argument bounds the sum of the squares of these distances. Variance tensorization, applied with the costs obtained after deletion, then gives the moment estimates below; see \citet{EfronStein1981,BoucheronBousquetLugosiMassart2005}. The argument applies to every sampling law, including laws with atoms.

\begin{lemma}[Distribution-free moments]\label{lem:regular-moments}
Let $H$ be a fixed finite vector of regular costs, and let $H_n$ be its value on the first $n$ points of an iid sequence with any Borel probability law $\mu$ on $D$. For an independent $N_s\sim\operatorname{Poi}(s)$, put $P_s=H_{N_s}$. Uniformly over $\mu$,
\begin{equation}\label{eq:regular-count-increment}
 \|H_{n+1}-H_n\|_p\le C(n+1)^{-1/2},
 \qquad 1\le p\le2.
\end{equation}
For every finite $p\ge1$,
\begin{equation}\label{eq:regular-moments}
 \sup_{n\ge0}\|H_n-\E H_n\|_p
 +\sup_{s\ge1}\|P_s-\E P_s\|_p
 +\sup_{n\ge1}\|H_n-\E P_n\|_p
 \le C_p.
\end{equation}
If $\mu\ge cm$ for some $c>0$, then \eqref{eq:regular-count-increment} holds for every finite $p\ge1$, with right-hand side $C_pc^{-1/2}(n+1)^{-1/2}$.
\end{lemma}

\begin{proof}
We first consider one coordinate. For a configuration of $N\ge2$ labels in $D$, let $r_i$ be the distance from label $i$ to its nearest other label. The open discs of radii $r_i/2$ are disjoint. Indeed, an intersection would place two centers at distance less than $(r_i+r_j)/2\le\max(r_i,r_j)$, contradicting the definition of the larger nearest-neighbor distance. A zero-radius disc is empty, so the assertion also covers repeated locations. These discs lie in a square of side $1+\sqrt2$, and hence
\begin{equation}\label{eq:nearest-square-sum}
 \sum_{i=1}^N r_i^2\le\frac{4(1+\sqrt2)^2}{\pi}.
\end{equation}
Exchangeability gives $\E r_N^2\le C/N$. Applying \eqref{eq:regular-insertion} to the last label proves \eqref{eq:regular-count-increment}. The increment at $n=0$ is zero.

Write $F=H_N$, and let $F_i$ be the cost after deleting label $i$. The same argument gives the pointwise estimate
\begin{equation}\label{eq:regular-deletion-square-sum}
 \sum_{i=1}^N|F-F_i|^2\le V
\end{equation}
with $V$ depending only on $H$. For $N\le1$, all differences are zero. Since $F_i$ does not depend on the $i$th input, variance tensorization gives $\Var F\le V$.

Set $U=F-\E F$ and $U_i=F_i-\E F$. Apply the same tensorization inequality to $|U|^j$, using $|U_i|^j$ as the functions independent of the $i$th input. For each integer $j\ge2$, this gives
\begin{align}
 \E|U|^{2j}
 &\le (\E|U|^j)^2
   +\sum_i\E\bigl(|U|^j-|U_i|^j\bigr)^2 \notag\\
 &\le (\E|U|^j)^2
   +j^2V\E(|U|+\sqrt V)^{2j-2}.
 \label{eq:regular-moment-recursion}
\end{align}
The deterministic growth bound ensures that these moments exist at every fixed $N$. Starting with the second moment, induction in $j$ and interpolation give all fixed centered moments uniformly in $N$ and $\mu$.

Let $a_n=\E H_n$. Summing the first-moment increment estimate gives
\begin{equation}\label{eq:regular-mean-modulus}
 |a_j-a_\ell|\le2C|\sqrt j-\sqrt\ell|.
\end{equation}
For $s\ge1$, put $\ell=\lfloor s\rfloor$. The elementary Poisson moment estimate $\|N_s-s\|_p\le C_p\sqrt s$ implies
\[
 \|\sqrt{N_s}-\sqrt\ell\|_p
 \le\frac{\|N_s-\ell\|_p}{\sqrt\ell}\le C_p.
\]
For completeness, the Poisson estimate follows for even integer moments by expanding
$\exp\{s(e^t-1-t)\}$, whose coefficient of $t^{2j}$ is a polynomial in $s$ of degree at most $j$, and then by interpolation. Conditioning on $N_s$, first center $H_{N_s}$ at $a_{N_s}$ and then use \eqref{eq:regular-mean-modulus}. This proves the Poisson moment bound. At $s=n$ it also gives $|\E P_n-a_n|\le C$, proving the final bound in \eqref{eq:regular-moments}.

Finally, a square contained in $B(y,t)\cap D$ shows that, for $y\in D$ and $0<t\le\sqrt2$,
\[
 m(B(y,t)\cap D)\ge t^2/8.
\]
When $\mu\ge cm$, the nearest-neighbor distance from an independent query point to $n\ge1$ old labels therefore has tail at most $\exp(-cnt^2/8)$. Integrating this tail and applying \eqref{eq:regular-insertion} gives the stated $L^p$ increment bound. Sum the scalar estimates over the fixed number of coordinates to obtain the vector assertions.
\end{proof}

\subsection{Empty blocks and the centered cutting error}

Let $X$ be a Poisson process with finite intensity measure $\lambda\mu$ on $D$. Partition $D$ into $k^2$ equal squares and define
\[
 R_k(X)=H(X)-\sum_{Q\in\mathcal C_k}H(X\cap Q).
\]
Whenever a square is subdivided, use a consistent Borel assignment of its boundary points to its children. Thus the assigned cells are disjoint, while their closures are the usual geometric squares.

For a block $B\subset Q$ of side $b$, replace its assigned points by an independent Poisson configuration with the same intensity measure. Denote the result by $X^{(B)}$ and write
\[
 \Delta_B=R_k(X)-R_k(X^{(B)}).
\]
Let $I,I'$ be the original and replacement configurations in $B$. For $U=D,Q$, let $O_U$ consist of the points assigned to $U$ but not to $B$. Choose $z\in\partial B$ and set
\[
 F(I)=H(I\sqcup O_D)-H(I\sqcup O_Q),
 \qquad T=H(O_D\sqcup\{z\})-H(O_Q\sqcup\{z\}).
\]
Then $\Delta_B=F(I)-F(I')$. If $I$ is nonempty, \eqref{eq:regular-full-boundary} cancels its interior cost and gives
\[
 |F(I)-T|\le2B_0b.
\]
If $I$ is empty, \eqref{eq:regular-insertion} gives
\[
 |F(I)-T|\le K(d_D+d_Q),
 \qquad
 d_U=\begin{cases}
 \operatorname{dist}(z,O_U),&O_U\ne\varnothing,\\
 0,&O_U=\varnothing.
 \end{cases}
\]
Consequently, if $\mathcal E_B$ is the event that exactly one of $I,I'$ is empty, then
\begin{equation}\label{eq:regular-excess-pointwise}
 \bigl(|\Delta_B|-(1+4B_0)b\bigr)_+^2
 \le C\1_{\mathcal E_B}(b^2+d_D^2+d_Q^2).
\end{equation}
Lemma~\ref{lem:regular-insert} allows labels on the block boundary, so the same bound holds for the chosen assignment of boundary points.

\begin{lemma}[Regular cutting]\label{lem:regular-cut}
Suppose that $\mu=g\,m+\tau$ is a fixed finite measure on $D$, where $\tau\perp m$ and $g>0$ $m$-almost everywhere. For any fixed finite vector of regular costs,
\begin{equation}\label{eq:regular-cut-limit}
 \|R_k-\E R_k\|_2\longrightarrow0
 \qquad\text{whenever }\lambda/k^2\longrightarrow\infty.
\end{equation}
The singular part of $\mu$ and the mass assigned to grid faces are unrestricted.
\end{lemma}

\begin{proof}
It suffices to prove the assertion for one coordinate. Fix one generation of blocks, all of side $b$, and realize all independent replacements using a second Poisson process $Y$ with intensity $\lambda\mu$. Mark a block if exactly one of $X(B),Y(B)$ is empty. Since the exterior of $B$ is independent of both block configurations and the two interiors have the same law, each expected term
$\E[\1_{\mathcal E_B}d_U^2]$ equals twice its contribution from marked blocks with $X(B)>0$.

Consider first one ambient square $U$, either $D$ or one of the roots $Q$. For an occupied block with $d_U>200b$, associate the set
\[
 A_B=(B(z,d_U/8)\cap U)\setminus B(z,3b).
\]
The area estimate for a ball intersected with a square, used in Lemma~\ref{lem:regular-moments}, gives
\begin{equation}\label{eq:regular-hole-area}
 m(A_B)\ge d_U^2/512-9\pi b^2\ge c_0d_U^2.
\end{equation}
Let $U^\sharp$ denote the assigned ambient cell (and let $U^\sharp=D$ when $U=D$). For $y\in A_B$, the query region
\[
 V_y=B(y,b/8)\cap U^\sharp
\]
avoids the closure of $B$ and contains no point of $X$. Indeed, $|y-z|>3b$ separates it from $B$, whereas $|y-z|<d_U/8$ and $d_U>200b$ place it strictly within distance $d_U$ of $z$.
Using the assigned cell $U^\sharp$ makes $V_y$ count precisely the labels in $O_U$, even when the sampling measure charges a partition face. The face has zero area, so this choice leaves the area and density estimates unchanged.

The regions $A_B$ belonging to distinct occupied blocks in the same ambient square are disjoint. If $A_B$ and $A_{B'}$ intersect, choose the larger of their two radii, say $d_U(B)\ge d_U(B')$. The chosen boundary points are then at distance less than $d_U(B)/4$. Any point of the occupied block $B'$ lies within $\sqrt2b$ of its boundary point and is exterior to $B$. It would therefore lie at distance less than
\[
 d_U(B)/4+\sqrt2b<d_U(B)
\]
from the chosen boundary point of $B$, a contradiction. For different ambient roots, the associated regions lie in disjoint roots up to area-zero boundaries. By \eqref{eq:regular-hole-area}, the total expected large-radius contribution is bounded by
\begin{equation}\label{eq:regular-hole-integral}
 C\int_D\Pp\{X(V_y)=0\}\,dm(y).
\end{equation}
This estimate applies separately to $U=D$ and to the family $U=Q$. Radii at most $200b$ contribute at most $Cb^2$ per marked block.

For $t>0$, put $G_t=\{g<t\}$. There is an absolute $\kappa>0$ such that $m(V_y)\ge\kappa b^2$, uniformly over the ambient squares. Fubini's theorem gives
\[
 \int_D m(V_y\cap G_t)\,dm(y)\le Cb^2m(G_t).
\]
Outside a set of area at most $C'm(G_t)$, it follows that
$m(V_y\setminus G_t)\ge\kappa b^2/2$. Thus the mean count in $V_y$ is at least $\kappa\lambda tb^2/2$ there. The integral in \eqref{eq:regular-hole-integral} is consequently bounded by
\[
 C\{m(G_t)+\exp(-c\lambda tb^2)\}.
\]
For the short-radius contribution, the total area of blocks satisfying
$m(B\cap G_t)>b^2/2$ is at most $2m(G_t)$. Every other block has empty probability at most $\exp(-\lambda tb^2/2)$. Combining these estimates with \eqref{eq:regular-excess-pointwise} yields, for every generation,
\begin{equation}\label{eq:regular-generation-excess}
 \sum_B\E\bigl(|\Delta_B|-(1+4B_0)b\bigr)_+^2
 \le C\{m(G_t)+\exp(-c\lambda tb^2)\}.
\end{equation}
The lower bounds on mean counts use only $g\,m$ and hence remain valid in the presence of the singular measure $\tau$.

We now apply Theorem~\ref{thm:multiscale}. Fix an integer $J\ge8$. Starting from the $k^2$ roots, construct a forest of height $J$ by dividing each node into $2^J\times2^J$ children. Give each node its Euclidean area as mass. The root masses sum to one, the square-root child-to-parent ratio is $2^{-J}$, and at every nonroot depth,
\[
 \lambda b^2\ge (\lambda/k^2)2^{-2J^2}.
\]
For fixed $J$ and $t>0$, sum \eqref{eq:regular-generation-excess} over the finitely many depths and let $\lambda/k^2\to\infty$. The limiting upper bound is at most $CJm(G_t)$. Letting $t\downarrow0$ makes it zero, since $g>0$ almost everywhere. Thus the resampling excess in Theorem~\ref{thm:multiscale} tends to zero for this fixed forest height. With its notation, $\delta=2^{-J}$ and $\varepsilon\to2/J$, so
\[
 \limsup_{\lambda/k^2\to\infty}\Var R_k
 \le \frac{C}{\log J},
\]
where $C$ is independent of $J$. Integrability follows from \eqref{eq:regular-growth} and the finiteness of all Poisson count moments. Let $J\to\infty$.

If partition faces carry positive $\mu$-mass, the assigned cell counts are still independent, the area estimates are unchanged, and the boundary comparison allows labels assigned to either side. Equivalently, for each fixed finite configuration one may perturb every label into the interior of its assigned finest cell and pass to the limit using Lemma~\ref{lem:regular-insert}. Applying the scalar conclusion to the finitely many coordinates proves \eqref{eq:regular-cut-limit}.
\end{proof}

\section{Boundary comparisons for parity-constrained costs}\label{sec:parity}

Perfect matching requires an even number of labels. After a sample is
partitioned, some cells may contain an odd number, so their local matching
costs require one additional label. We place that label at a prescribed
boundary point.

For a Poisson input, conditioning on the number of points in each cell modulo two preserves
independence of the cell configurations. Replacements that keep these
parities fixed can then be treated by the multiscale variance theorem.
To recover the unconditional variance, we must also control the variation
of the conditional mean as the parities change. These are the two steps
of the cutting argument in this section.

For \(a\in\{0,1\}\), let \(\{z\}^{[a]}\) denote \(a\) new labels at \(z\).
All unions of labeled configurations below are denoted by \(\sqcup\);
coincident labels are retained.  Distances are Euclidean.

\begin{definition}[Parity cost]\label{def:parity}
A parity cost is a real-valued functional \(H\) on finite even subsets of
\(\R^2\), including the empty set, satisfying
\[
 H(x+rS)=rH(S),\qquad x\in\R^2,\quad r>0.
\]
There exist \(z_*\in\partial D\) and \(\kappa_\partial<\infty\) such that,
whenever \(I\subset D^\circ\) and \(O\subset\R^2\setminus D\) are finite
sets of distinct sites, \(|I|+|O|\) is even, and \(a=|I|\bmod 2\),
\begin{equation}\label{eq:parity-interface}
 \left|
 H(I\cup\{z_*\}^{[a]})+
 H(O\cup\{z_*\}^{[a]})-H(I\cup O)
 \right|\le \kappa_\partial.
\end{equation}
\end{definition}

\subsection{Continuity, paired insertion, and boundary comparison}

The boundary comparison implies continuity in the point locations. We
first establish this fact so that the later arguments apply to atoms
and to labels on partition faces.

\begin{lemma}\label{lem:parity-insert}
A parity cost has a unique continuous, permutation-invariant extension
to \( (\R^2)^n \) for each even \(n\).  Put \(L=3\kappa_\partial\).
For corresponding labeled configurations of the same even cardinality,
\begin{equation}\label{eq:parity-displacement}
 |H(x_1,\ldots,x_n)-H(y_1,\ldots,y_n)|
 \le L\sum_{i=1}^n|x_i-y_i|.
\end{equation}
For every even labeled configuration \(S\),
\begin{equation}\label{eq:parity-pair-insertion}
 H(S\sqcup\{x,x\})=H(S),\qquad
 |H(S\sqcup\{x,y\})-H(S)|\le L|x-y|.
\end{equation}
Let \(Q\) be a square of side \(b\), let \(z\in\partial Q\), and split an
even labeled configuration into labels \(I\) assigned to \(Q\) and labels
\(O\) assigned outside \(Q^\circ\).  Labels on \(\partial Q\) may be
assigned to either side.  With \(a=|I|\bmod2\) and
\(\kappa=\kappa_\partial+2\sqrt2L\),
\begin{equation}\label{eq:parity-full-boundary}
 \left|
 H(I\sqcup\{z\}^{[a]})+
 H(O\sqcup\{z\}^{[a]})-H(I\sqcup O)
 \right|\le \kappa b.
\end{equation}
Finally, if \(S,S'\) are even labeled configurations and their shared
labels are identified, then
\begin{equation}\label{eq:parity-symmetric-difference}
 |H(S)-H(S')|\le L\,\mathsf M_{\mathrm E}(S\mathbin{\triangle}S'),
\end{equation}
where \(\mathsf M_{\mathrm E}\) is the minimum Euclidean perfect-matching
length on the remaining labels.
\end{lemma}

\begin{proof}
Write \(K=D-z_*\) and \(p(v)=H(\{0,v\})\).  For \(v\in bK^\circ\), let
\(r=\inf\{t>0:v\in tK\}\).  If \(r<t<2r\), the square \(tK\), whose
prescribed anchor is \(0\), contains \(v\) in its interior and excludes
\(2v\).  Apply \eqref{eq:parity-interface} with \(I=\{v\}\) and
\(O=\{2v\}\).  Translation invariance and homogeneity give
\[
 |p(v)+p(2v)-H(\{v,2v\})|=|2p(v)|\le\kappa_\partial t.
\]
Letting \(t\downarrow r\) shows that
\begin{equation}\label{eq:parity-two-site}
 |p(v)|\le \tfrac12\kappa_\partial b,\qquad v\in bK^\circ.
\end{equation}

Fix an odd set \(F\).  Suppose \(x,y\) lie in a square of side \(b\) whose
closure avoids \(F\), and compare \(H(F\cup\{x\})\) and
\(H(F\cup\{y\})\) through that square's prescribed anchor.  The completed
exterior is the same in both comparisons.  Each completed interior cost
has absolute value at most \(\kappa_\partial b/2\), by
\eqref{eq:parity-two-site}.  Thus the difference is at most
\(3\kappa_\partial b\).  Locally off \(F\), the side \(b\) can decrease
to \(\|x-y\|_\infty\).

Any two points outside the finite set \(F\) can be joined by a polygonal
path avoiding \(F\), of length at most \(|x-y|+\varepsilon\).  Indeed,
choose an intermediate point arbitrarily close to their midpoint and
outside the finitely many lines that would force a segment through a
point of \(F\).  Subdivide the resulting path into sufficiently short
pieces and apply the local bound.  Letting \(\varepsilon\downarrow0\)
proves the one-label Lipschitz estimate off \(F\).  To compare two lists
of distinct sites, first move the labels to distinct auxiliary sites
arbitrarily close to their targets and outside both lists, and then to
their targets.  Changing one coordinate at a time gives
\eqref{eq:parity-displacement} on distinct lists.  Density gives the
unique continuous extension to every even Cartesian power; permutation
invariance passes to the extension.

Homogeneity forces the empty cost, and the cost of an even list
supported at one point, to be zero.  To extend the boundary comparison,
approximate the designated inner labels from \(Q^\circ\) and the
designated outer labels from the complement of \(\overline Q\), making
all sites distinct.  The continuity just proved passes the comparison to
the limit.  This includes coincident labels on a shared face that have
been assigned to different sides.  Moving the two possible copies of
the prescribed anchor to \(z\) changes the sum by at most
\(2\sqrt2Lb\), proving \eqref{eq:parity-full-boundary}.

For \(x\) outside the support of \(S\), enclose the two new labels at
\(x\) in squares shrinking to \(x\).  In
\eqref{eq:parity-full-boundary} the interior has even cardinality and
zero cost, so the comparison error tends to zero.  Hence
\(H(S\sqcup\{x,x\})=H(S)\).  Continuity removes the restriction on \(x\).
Moving one of the two new labels proves the second assertion in
\eqref{eq:parity-pair-insertion}.

Choose a perfect matching of \(S\mathbin{\triangle}S'\).  Delete every
old--old pair, insert every new--new pair, and move the old label to the
new site in each mixed pair.  Every intermediate configuration has even
cardinality.  Equations \eqref{eq:parity-displacement} and
\eqref{eq:parity-pair-insertion} bound the total change by \(L\) times
the matching length.  Minimizing proves
\eqref{eq:parity-symmetric-difference}.
\end{proof}

For a square \(Q\), let \(z_Q\) be its lower-left corner and define the
completed cost
\[
 A_Q(S)=H\bigl(S\sqcup\{z_Q\}^{[|S|\bmod2]}\bigr).
\]
Consecutive pairing in the visiting order of Lemma~\ref{lem:short-path}
and \eqref{eq:parity-pair-insertion} give
\begin{equation}\label{eq:parity-growth}
 |H(S)|\le4Lb\sqrt{|S|}\quad(|S|\ \text{even}),\qquad
 |A_Q(S)|\le4Lb\sqrt{|S|+1}
\end{equation}
when \(S\subset Q\) and \(Q\) has side \(b\).  The labeled extensions
and their completions are Borel measurable.  Finite linear combinations
of parity costs satisfy the comparison in \eqref{eq:parity-full-boundary} and paired
insertion bound, with constants summed using absolute coefficients.

Some costs, including matching, are defined directly on labeled
configurations. If such a cost is translation invariant, homogeneous
of degree one, permutation invariant, and satisfies
\eqref{eq:parity-full-boundary}, it agrees with the
continuous extension just constructed.  Enclose each distinct
location in a disjoint square of side \(\varepsilon\), and separate its
labels in a fixed scaled pattern.  In the comparisons before and after
separation the exterior terms cancel, whereas each completed interior
has cost \(O(\varepsilon)\) by homogeneity.  Summing over the finitely
many locations shows that the original and separated costs differ by
\(O(\varepsilon)\).  The separated configurations converge to the
original labeled configuration, so their continuous extensions converge
to its value.  Thus both formulations agree.

\subsection{Changes of sample count and of local parity}

Let \(A_n\) be a fixed finite vector evaluated at \(n\) iid points with
Borel law \(\nu\) on \(D\).  Parity coordinates use completion at one
fixed point \(z\in D\), and regular coordinates have their ordinary
values.  Let \(P_s\) be the same vector at an independent
\(\operatorname{Poi}(s)\) count.  Constants denoted by \(C_H\) below
depend only on this fixed vector and its comparison constants.

\begin{lemma}[Neighboring counts]\label{lem:neighbor-count}
For every Borel probability \(\nu\) and every \(n\ge0\),
\begin{equation}\label{eq:parity-step-two-bl}
 \dBL\bigl(\Law(A_{n+2}),\Law(A_n)\bigr)
 \le \frac{C_H}{\sqrt{n+1}}.
\end{equation}
If all coordinates are parity costs and
\(L_H=(\sum_jL_j^2)^{1/2}\), there are couplings such that
\begin{equation}\label{eq:parity-step-two-coupling}
 \Pp\left\{|A_{n+2}-A_n|>\frac{CL_H}{\sqrt{n+1}}\right\}
 \le\frac{C}{\sqrt{n+1}}.
\end{equation}
For \(n\ge4\) and even \(d=|m-n|\le n\), there is also a coupling with
\begin{equation}\label{eq:parity-long-count-coupling}
 \Pp\left\{|A_m-A_n|>\frac{CL_Hd}{\sqrt n}\right\}
 \le\frac{Cd}{\sqrt n}.
\end{equation}
If \(p=\nu(B(z,\eta))>0\), then
\begin{equation}\label{eq:parity-step-one-bl}
 \dBL\bigl(\Law(A_{n+1}),\Law(A_n)\bigr)
 \le C_H\left\{\eta+\frac1{\sqrt{(n+1)p}}\right\}.
\end{equation}
\end{lemma}

The bounded Lipschitz distance is invariant under a common translation,
so the same comparisons hold when both costs are shifted by any
fixed vector.

\begin{proof}
Put \(q=\lceil\sqrt{n+1}\rceil\), partition \(D\) into \(q^2\) squares,
and write \(p_i>0\), \(1\le i\le J\), for their positive \(\nu\)-masses.
Independently of the old sample, choose \(i\) with probability \(p_i\)
and add two independent sites with law \(\nu(\,\cdot\,|C_i)\).  After a
uniform random permutation, the augmented list has likelihood
\[
 W=\binom{n+2}{2}^{-1}\sum_{a<b}K(X_a,X_b),
 \qquad
 K(x,y)=\sum_{i=1}^J p_i^{-1}\1_{\{x,y\in C_i\}},
\]
relative to an iid \(\nu\)-sample of size \(n+2\).  The first projection
of \(K-1\) vanishes:
\[
 \E[K(X,Y)\mid X]=1,\qquad \E K(X,Y)^2=J.
\]
Thus distinct pair terms are uncorrelated, including pairs sharing one
index.  This is the variance calculation for a degenerate \(U\)-statistic
\citep{Hoeffding1948}. It gives
\begin{equation}\label{eq:parity-likelihood-variance}
 \|W-1\|_2^2=\frac{J-1}{\binom{n+2}{2}}\le\frac{C}{n+1}.
\end{equation}
The total variation distance is at most \(\|W-1\|_2/2\).
The two added sites are at distance at most \(\sqrt2/q\).  Equation
\eqref{eq:parity-pair-insertion} therefore bounds the change of the
completed parity vector by \(2L_H/q\).  Maximally couple the augmented
list to an iid list of size \(n+2\), and retain the old list through its
conditional distribution given the augmentation.  This proves
\eqref{eq:parity-step-two-coupling}.

If \(d<\sqrt n\), compose \(d/2\) neighboring couplings.  The
intermediate counts are at least \(n-\sqrt n\ge n/2\), so the accumulated
cost and the union bound give
\eqref{eq:parity-long-count-coupling}.  If \(d\ge\sqrt n\), its
right-hand side is at least one after increasing \(C\), and the
assertion is automatic.  This also covers \(m\) much smaller than \(n\).

For a regular coordinate with insertion constant \(K_0\), the
augmentation changes the cost by at most
\[
 K_0\{2r_1(x,S)+|x-y|\},\qquad n\ge1,
\]
where \(r_1(x,S)=\operatorname{dist}(x,S)\).
Indeed, insertion of a label can only decrease subsequent
nearest-neighbor distances, and
\(r_1(y,S)\le r_1(x,S)+|x-y|\).
The first added point has marginal law \(\nu\), independently of \(S\).
The packing estimate in Lemma~\ref{lem:regular-moments} gives
\(\E r_1(x,S)\le C/\sqrt{n+1}\).
Apply bounded Lipschitz tests and the preceding total variation bound.
The case \(n=0\) is absorbed in \(C_H\), proving
\eqref{eq:parity-step-two-bl} for mixed vectors.

For a one-step comparison, add one point sampled from
\(\nu(\,\cdot\,|B(z,\eta))\) and permute the labels.  Relative to an iid
list of size \(n+1\), the likelihood is
\[
 \frac{N_B}{(n+1)p},
 \qquad
 d_{\mathrm{TV}}\le
 \frac12\sqrt{\frac{1-p}{(n+1)p}}.
\]
For a parity coordinate, the added point and the changed completion
label form a paired insertion or a one-label move of cost at most
\(L\eta\).  For a regular coordinate, insertion costs at most
\(2K_0\eta\) if the old sample has a point in the ball; the complementary
event has probability \((1-p)^n\le C/\sqrt{(n+1)p}\).
Bounded Lipschitz tests prove \eqref{eq:parity-step-one-bl}.
\end{proof}

The next lemma compares costs when the parities in two neighboring
cells are reversed. The configurations outside these cells can be
arbitrary and are kept fixed. Since the total parity is unchanged,
the global completion label is also the same in both configurations.
The one-label count comparison is the Poisson instance of the
log-concave parity comparison in \citet{Pinelis2020ParityCoupling}.
A related parity-exchange involution appears in
\citet{Petrov2020ParityCoupling}; the construction below chooses the
exchanged pair by its geometric length.

\begin{lemma}[Local parity couplings]\label{lem:parity-couple}
For a Poisson process with finite, strictly positive total intensity,
the laws conditional on even and odd counts admit a coupling in which
the two clouds differ by exactly one label.

Consider two touching rectangles with side lengths at most \(h\).
Suppose their intensity measures dominate \(v\gamma_1\) and
\(v\gamma_2\), where \(v\ge1\) and \(\gamma\) is a probability measure
on pairs with marginals \(\gamma_1,\gamma_2\).  Set
\[
 F(t)=\gamma\{(x,y):|x-y|\le t\},
\]
and let \(D_0\) bound the diameter of the union of the two marginal
supports.  Under complementary prescribed cell-parity patterns, the
completed costs of the full configurations can be coupled so that their
expected absolute difference is at most
\begin{equation}\label{eq:parity-pair-coupling}
 CL\left\{\int_0^{D_0}e^{-2vF(t)}\,\mathrm dt
             +(h+D_0\sqrt v)e^{-2v}\right\}.
\end{equation}
In particular, side-adjacent squares containing a homogeneous component
of mass \(v\) in each square have bound \(CLh/v\).
Two patches of mass at least \(v\), with combined diameter at most \(d\),
have bound
\begin{equation}\label{eq:parity-two-patches}
 CL(d+he^{-2v}).
\end{equation}
In a cell \(Q\) with side lengths at most \(h\), if a patch of mass at
least \(v\) lies within distance \(d\) of a fixed completion point \(z\),
the opposite-parity comparison is
bounded by
\begin{equation}\label{eq:parity-one-patch}
 CL\{d+(h+\operatorname{dist}(z,Q))e^{-2v}\}.
\end{equation}
The last two assertions remain valid for \(0<v<1\), with an absolute
change of constants.
\end{lemma}

\begin{proof}
Write the even and odd Poisson counts of mean \(w>0\) as \(2J_e\) and
\(2J_o+1\).  The likelihood ratio of the mass of \(J_o\) at \(j\) to
that of \(J_e\) at \(j\) is a constant times \(w/(2j+1)\), and hence is
decreasing in \(j\).  This gives \(J_o\le_{\mathrm{st}}J_e\).
Likewise, on \(j\ge1\), the ratio of the mass of \(J_e\) at \(j\) to
that of \(J_o+1\) at \(j\) is a constant times \(w/(2j)\).  Including
the possible atom of \(J_e\) at zero gives
\(J_e\le_{\mathrm{st}}J_o+1\).
The common-quantile coupling therefore satisfies
\[
 J_o\le J_e\le J_o+1.
\]
Use one common iid sequence of locations and the two corresponding
prefixes.  The clouds differ by exactly one label.  Completion at a
point in the cell closure changes the cost by at most \(2Lh\).

Split off a selected process of mass \(v\) from each cell.  If the
residual process in one cell has mean \(w\), its likelihood, conditional
on the total cell parity \(p\), relative to its unconditional law is
\begin{equation}\label{eq:parity-residual-likelihood}
 \frac{1+(-1)^{p+N_{\mathrm{res}}}e^{-2v}}
      {1+(-1)^p e^{-2(v+w)}}.
\end{equation}
The two residual parity distributions under opposite target parities
are consequently at total variation distance \(O(e^{-2v})\), uniformly
in \(w\ge0\).  Couple these parity bits maximally.  Conditional on
agreement, use identical residual clouds.  Conditional on disagreement,
use the one-label coupling just proved.  Do this in both cells.

On a residual-disagreement branch, couple the selected clouds
identically if their required parities agree, and by the one-label
coupling otherwise.  In each cell exactly one of the selected and
residual parts changes parity.  Hence only one label per cell changes,
and \eqref{eq:parity-symmetric-difference} bounds the cost by \(CLh\).
The probability of this branch is \(O(e^{-2v})\).

It remains to couple the selected clouds when the residual parities
agree in both cells.  Generate a Poisson process of pairs of mean
\(2v\), with pair law \(\gamma\), and attach to each pair an independent
fair mark selecting either its first or its second endpoint.  Poisson
splitting gives independent selected marginal clouds with intensities
\(v\gamma_1\) and \(v\gamma_2\).
Write \(r=e^{-2v}\).  The probabilities of the four parity patterns
intersected with the event that the pair process is nonempty are
\[
 \begin{aligned}
 \Pp(00,\mathrm{nonempty})=\Pp(11,\mathrm{nonempty})
       &=\frac{(1-r)^2}{4},\\
 \Pp(01,\mathrm{nonempty})=\Pp(10,\mathrm{nonempty})
       &=\frac{1-r^2}{4}.
 \end{aligned}
\]
In particular these probabilities are bounded below by an absolute
positive constant for \(v\ge1\).  In the \(00\) case, the empty-process probability \(r\) has been
subtracted from the unconditional parity probability.

On a nonempty pair configuration, choose the shortest pair, breaking
ties by independent continuous tags, and flip its mark.  The choice of
pair is independent of all marks.  The flip preserves the unconditioned
marked-process law, is an involution, and toggles both cell parities.
The displayed equality of masses therefore makes it an exact coupling
of the two complementary parity laws conditional on nonemptiness.
If the selected pair has length greater than \(t\), the original
Poisson process contains no pair of length at most \(t\).  Dividing by
the above lower bound shows that this probability is at most
\(Ce^{-2vF(t)}\).  Moving the retained endpoint along that pair costs
at most its length times \(L\).

Couple the nonemptiness indicators first.  Their disagreement
probability is \(O(e^{-2v})\).  On a branch where these indicators disagree, the nonempty selected
cloud has conditional expected count \(O(v)\), because the event
specifying its parity and nonemptiness has probability bounded below
and the unconditional total count has mean \(2v\).
The selected sites lie in a square of side at most \(D_0\).
The visiting bound and
\eqref{eq:parity-symmetric-difference}, followed by Cauchy--Schwarz,
bound the expected cost on this branch by \(CLD_0\sqrt v\).
Integrating the shortest-pair tail and adding the residual and
nonemptiness exceptions proves \eqref{eq:parity-pair-coupling}.

For uniform reflected pairs across a common side,
\(F(t)=\min\{1,t/(2h)\}\).  The integral is \(O(h/v)\), and
\((h+D_0\sqrt v)e^{-2v}=O(h/v)\).
For two patches, on residual agreement use the two one-label
couplings directly: the two changed labels lie in their union and
cost at most \(Ld\).  On residual disagreement use the preceding
\(CLh\) bound, proving \eqref{eq:parity-two-patches}.
For one patch, pair the changed selected label with the changed
completion label.  The cost is at most \(Ld\) on residual agreement,
and at most \(CL(h+\operatorname{dist}(z,Q))\) otherwise.
This proves \eqref{eq:parity-one-patch}.  When \(0<v<1\),
\(e^{-2v}\ge e^{-2}\), so the direct one-label bounds imply both patch
assertions after increasing \(C\).
\end{proof}

\subsection{Distribution-free centered moments}

The moment estimate separates fluctuations at fixed bin parities from
fluctuations of the conditional mean. A space-filling curve orders the
sample while controlling distances, and a concentration bound for
cumulative parities controls the second term. We state these two tools
here; their proofs are in Appendix~\ref{app:curve-prefix}.

Space-filling orderings and their \(1/2\)-H\"older geometry are classical
\citep{BartholdiPlatzman1988Spacefilling,GaoSteele1994Spacefilling}.
Related invariant Hilbert paths appear in
\citet[Proposition~5]{HolroydLiggett2001ExtraHead}.

\begin{lemma}[A space-filling order]\label{lem:hilbert-order}
There is a continuous surjection \(\phi:[0,1]\to D\) starting and
ending at any prescribed ordered pair of adjacent corners of \(D\),
such that
\begin{equation}\label{eq:parity-hilbert}
 |\phi(u)-\phi(v)|\le C|u-v|^{1/2}.
\end{equation}
For every \(j\ge0\), the map sends the \(4^j\) equal closed parameter
intervals onto the \(4^j\) dyadic squares of side \(2^{-j}\), one
interval for each square. Its first-preimage map
\[
 b(x)=\min\phi^{-1}(\{x\})
\]
is Borel and satisfies \(\phi(b(x))=x\) for every \(x\in D\).
\end{lemma}

\begin{lemma}[Concentration for prefix parities]\label{lem:parity-prefix}
Let \(\sigma_1,\ldots,\sigma_J\) be independent bits and set
\(\tau_i=\sigma_1\oplus\cdots\oplus\sigma_i\).
Suppose that, at each step, the transition laws from opposite current
states have total variation distance at most \(\rho<1\).
Let \(f\) be a function of \((\tau_1,\ldots,\tau_J)\) with coordinate
oscillations at most \(a_i\). If \(\sum_i a_i^2=0\), then \(f\) is
constant. Otherwise,
\begin{equation}\label{eq:parity-prefix-concentration}
 \Var f\le\frac{\sum_i a_i^2}{4(1-\rho)^2},
 \qquad
 \Pp(|f-\E f|>t)
 \le2\exp\left\{-\frac{2(1-\rho)^2t^2}{\sum_i a_i^2}\right\}.
\end{equation}
\end{lemma}

Prefix coordinates for parity patterns also appear in
\citet[Lemma~2.6]{Feng2026MatchingVariance}. The independent Poisson
parities used here allow the direct coupling proof given in the appendix.

\begin{lemma}[Distribution-free parity moments]\label{lem:parity-moments}
For one parity coordinate, every Borel law \(\nu\) on \(D\), and every
\(1\le p<\infty\),
\begin{equation}\label{eq:parity-moment-bound}
 \sup_{s\ge1}\|P_s-\E P_s\|_p+
 \sup_{n\ge1}\|A_n-\E P_n\|_p\le CL\sqrt p.
\end{equation}
The same bound holds for the fixed-count variables centered at their
own means, and
\begin{equation}\label{eq:parity-mean-step}
 |\E A_{n+2}-\E A_n|\le CL(n+1)^{-1/2}.
\end{equation}
These assertions are uniform over \(\nu\), including atomic and
singular laws.  A fixed finite vector is controlled coordinatewise.
\end{lemma}

\begin{proof}
If \(L=0\), the growth bound makes all costs zero, so assume \(L>0\).
Choose the Hilbert curve to end at the completion corner.
Let \(Q_\nu\) be a quantile function for \(b(X)\), \(X\sim\nu\).
Since \(\phi(b(X))=X\), mapping a Poisson process of intensity \(s\) on
\([0,1]\) through \(\phi\circ Q_\nu\) gives a labeled Poisson process
of intensity \(s\nu\).

Put \(J=\lfloor s\rfloor\), partition the parameter interval into
\(J\) equal bins, and set
\[
 t_0=0,\qquad t_i=Q_\nu(i/J)\ (1\le i<J),\qquad
 t_J=1,\qquad d_i=t_i-t_{i-1}.
\]
The bin clouds are independent, their counts have the common mean
\(s/J\in[1,2)\), and the sites from bin \(i\) lie in
\(\phi([t_{i-1},t_i])\).  This construction remains valid for labeled
clouds when one atom occupies several bins or \(d_i=0\).

Condition on all bin parities \(\pi\).  The bin clouds remain
independent.  Reveal them in order, and let \(D_i\) be the Doob
increments of \(P_s\).  Replacing bin \(i\) by an independent cloud
with the same prescribed parity preserves the total parity and hence
the global completion.  The union of the two bin clouds is contained
in a square of side \(C\sqrt{d_i}\), by
\eqref{eq:parity-hilbert}.  Equations
\eqref{eq:parity-symmetric-difference} and
\eqref{eq:parity-growth} give the replacement bound
\[
 CL\sqrt{d_i}\sqrt{N_i+N_i'}.
\]
Averaging over the replacement and the unrevealed bins gives
\[
 |D_i|\le CL\sqrt{d_i}(\sqrt{N_i}+1).
\]
The Poisson means are bounded, and either count parity has probability
at least \((1-e^{-2})/2\).  Consequently, for an absolute \(C_1\),
\begin{equation}\label{eq:parity-conditional-exponential-square}
 \E\left[
 \left.\exp\left\{\frac{D_i^2}{C_1L^2d_i}\right\}
 \right|\,\text{past},\pi\right]\le2
\end{equation}
when \(d_i>0\); for \(d_i=0\), paired coincident-label invariance gives
\(D_i=0\).

We use the elementary implication
\[
 \E Y=0,\quad \E e^{Y^2/a^2}\le2
 \quad\Longrightarrow\quad
 \E e^{tY}\le e^{Ca^2t^2},\qquad t\in\R.
\]
For \(|t|a\le1\), apply
\(e^x-1-x\le x^2e^{|x|}/2\) and Cauchy--Schwarz, using the
exponential-square moment.  For \(|t|a>1\), use
\(tY\le Y^2/(2a^2)+t^2a^2/2\) and absorb the constant factor into
\(e^{Ca^2t^2}\).
Apply this implication conditionally in
\eqref{eq:parity-conditional-exponential-square} and iterate.
Since \(\sum_i d_i=1\), the conditional fluctuation
\(P_s-\E[P_s\mid\pi]\) is sub-Gaussian with scale \(CL\),
uniformly in \(\pi\).

It remains to bound the fluctuations of the conditional mean
\(g(\pi)=\E[P_s\mid\pi]\). Express \(g\) as a function of the
prefix parities.
Flipping an internal prefix changes precisely the parities of bins
\(i\) and \(i+1\).  Lemma~\ref{lem:parity-couple} changes one label in
each bin.  The two changed sites lie in
\(\phi([t_{i-1},t_{i+1}])\), so their comparison costs at most
\(CL\sqrt{d_i+d_{i+1}}\).
Flipping the last prefix changes the last-bin parity and the global
completion.  The changed site can be paired with \(\phi(1)\), at cost
at most \(CL\sqrt{d_J}\).
The sum of squared coordinate oscillations of \(g\) is therefore at
most \(CL^2\).
The transition laws of the prefix chain from opposite states have
total variation distance \(e^{-2s/J}\le e^{-2}\).
Equation \eqref{eq:parity-prefix-concentration}, together with the
conditional estimate, proves
\begin{equation}\label{eq:parity-poisson-subgaussian}
 \Pp(|P_s-\E P_s|>t)\le C e^{-ct^2/L^2}.
\end{equation}

Divide the cost by \(L\).
Lemma~\ref{lem:neighbor-count} gives a uniform step-two
bounded-Lipschitz comparison for this normalized cost.
Apply Lemma~\ref{lem:exact-count} to its Poisson mixture of mean \(n\),
at the center \(\E P_n/L\).
That lemma bounds the tail at a fixed count by a constant times the
square root of the corresponding Poisson tail. By
\eqref{eq:parity-poisson-subgaussian}, the fixed-count cost therefore
also has a uniform sub-Gaussian tail; the bounded shift in the tail
argument is absorbed into the constants.  Integration gives
\eqref{eq:parity-moment-bound}.
The \(p=1\) estimate bounds the change of center from \(\E P_n\)
to \(\E A_n\), giving the assertion with the fixed-count mean.

Finally, use the two-point augmentation in
Lemma~\ref{lem:neighbor-count}.
Its expected cost differs from \(\E A_n\) by at most \(2L/q\).
Its likelihood relative to the iid sample of size \(n+2\) is \(W\),
so its expected cost differs from \(\E A_{n+2}\) by at most
\[
 \left|\E[(W-1)(A_{n+2}-\E A_{n+2})]\right|
 \le \|W-1\|_2\sqrt{\Var(A_{n+2})}.
\]
The moment bound already proved and
\eqref{eq:parity-likelihood-variance} prove
\eqref{eq:parity-mean-step}.
An alternative fixed completion point in \(D\) changes the cost by at
most \(\sqrt2L\); the moment bounds persist.  The same paired
augmentation, whose completion parity is unchanged, proves the mean
increment bound for that completion as well.
\end{proof}

\subsection{Conditional cutting and the conditional mean}

Write \(A=A_D\) and, on the \(k\times k\) partition, set
\[
 R_k=A-\sum_{Q\in\mathcal C_k}A_Q.
\]
Faces are assigned by a fixed Borel rule, consistent under refinement.
Suppose a block \(B\) lies in a root square \(Q\), and replace its
assigned labels while preserving their parity.
The global and root completion labels stay fixed.  Both lie outside
\(B^\circ\); if a completion coincides with a resampled label, retain
it as a separate exterior label.
Apply \eqref{eq:parity-full-boundary} to the global and local costs,
using the same anchor on \(\partial B\).  The completed interior cost
cancels, and both completed exterior costs stay fixed before and after
replacement.  Thus
\begin{equation}\label{eq:parity-block-oscillation}
 \operatorname{osc}_B(R_k)\le4\kappa\,\operatorname{side}(B),
\end{equation}
including the case of an empty block.
The multiscale theorem applied at fixed parities now controls the
conditional variance of \(R_k\). To bound its full variance, it remains
to show that \(\E[R_k\mid\pi]\) also has vanishing variance. The next
lemma will be used for this second step.

\begin{lemma}[Variance of the conditional mean on a square grid]\label{lem:state-var}
Divide a square \(U\) into \(m^2\) leaves of side \(h\), with positive
Poisson means \(\mu_i\).  Write \(\pi_U\) for their parity vector and
\(g_U=\E[A_U\mid\pi_U]\).
For each side edge \(e\) of the leaf adjacency graph, suppose
\(Lhb_e\) bounds the oscillation of \(g_U\) when the two endpoint
parities are flipped, with all other parities fixed.
Then, with \(r_i=e^{-2\mu_i}\),
\begin{equation}\label{eq:parity-state-variance}
 \Var(g_U)\le CL^2h^2
       \left(1+\sum_i r_i+\sum_e b_e^2\right).
\end{equation}
\end{lemma}

\begin{proof}
Scale \(U\) to \(D\), choose the curve in
Lemma~\ref{lem:hilbert-order} to end at the completion corner, and
order the \(m\)-grid centers by their first preimages.
Consecutive times, and the final gap to \(1\), have length at most
\(C/m^2\).
To verify this for \(m\ge4\), choose a dyadic square side between
\(2/m\) and \(4/m\).
Every such dyadic square contains an \(m\)-grid center strictly in
its interior.  Every preimage of that center belongs to the parameter
interval of this dyadic square.  A gap between consecutive center
times therefore cannot contain a whole such interval, and has length
at most two interval lengths.  Increasing \(C\) covers the finitely many cases \(m<4\).
Equation \eqref{eq:parity-hilbert} now shows that the union of leaves
\(i,\ldots,j\) has diameter at most \(Ch\sqrt{j-i+1}\), and that the
last leaf lies within \(Ch\) of the completion corner.

Route each consecutive pair of centers by shortest horizontal and
vertical grid steps.  Each route has \(O(1)\) edges and lies in a
\(Ch\)-ball about its distinct starting center, so every grid edge
belongs to \(O(1)\) routes.
If \(a_i\) is the sum of the \(b_e\)'s along route \(i\), telescoping
successive edge flips gives a two-parity oscillation bound \(Lha_i\).
The bounded lengths and overlaps give
\begin{equation}\label{eq:parity-route-energy}
 \sum_{i<m^2}a_i^2\le C\sum_e b_e^2.
\end{equation}

Put \(N=m^2\). We construct the prefix parities from independent
random variables, one for each step.  At step \(i\), copy the preceding prefix with probability
\(r_i\), and otherwise reset to an independent fair bit.
The mod-two increment at step \(i\), conditionally on its past, is
one with probability \((1-r_i)/2\).  The increments are therefore
independent and have exactly the required Poisson parity laws.

Resample the random variables used at step \(i\), keeping those
at all other steps fixed.
The prefix parity at step \(i\) changes with probability at most \(2(1-r_i)\).
If it changes, the two prefix paths remain different until the next
reset \(T>i\), with \(T=N+1\) if there is no further reset.
Consequently only the parities at \(i\) and \(T\) change; the second
change is omitted when \(T=N+1\).
For \(i<N\) and \(T=i+1\), the cost is at most \(Lha_i\).
For \(T>i+1\), the one-label couplings and
\eqref{eq:parity-symmetric-difference} bound the change by
\[
 CLh\sqrt{T-i+1}.
\]
When \(T=N+1\), the second label in this comparison is the completion
corner; the same bound follows from the preceding geometric estimate.
Resampling the last step, \(i=N\), changes the cost by at most \(CLh\).

The initial change is independent of the future reset indicators, and
\[
 \Pp(T>t)=\prod_{i<j\le t}r_j,\qquad i<t\le N.
\]
The tail-sum identity gives
\[
 \E[(T-i+1)\1_{\{T>i+1\}}]
 \le4\sum_{t=i+1}^N\prod_{i<j\le t}r_j .
\]
Apply the Efron--Stein inequality to the independent random variables
in this construction \citep{EfronStein1981}. The contribution from an
immediate reset is bounded by \(CL^2h^2\sum_i a_i^2\), the last step
contributes at most \(CL^2h^2\),
and the delayed-reset contribution by a constant times
\[
 L^2h^2\sum_{t=2}^N\sum_{i<t}(1-r_i)
                              \prod_{i<j\le t}r_j.
\]
The exact telescoping identity
\[
 \sum_{i<t}(1-r_i)\prod_{i<j\le t}r_j
 =r_t\left(1-\prod_{i<t}r_i\right)\le r_t
\]
and \eqref{eq:parity-route-energy} prove
\eqref{eq:parity-state-variance}.
\end{proof}

\begin{lemma}[Centered parity cutting]\label{lem:parity-cut}
Let a Poisson process on \(D\) have intensity \(\lambda\mu\), where
\(\mu=g\,\mathrm dx+\tau\) is finite, \(g\ge0\), and
\(\tau\perp\mathrm dx\).
Fix \(k\ge1\).  Use \(H_0\ge4\) generations, each dividing every square
into \(2^{H_0}\times2^{H_0}\) children, and put
\[
 J=2^{H_0^2},\qquad h=(kJ)^{-1}.
\]
Inside each finest leaf, against the midpoint of each side, place an
open square patch of side \(\varepsilon h/3\), where
\(0<\varepsilon\le1\).
Call a leaf bad if at least one of its four patches has
\(\lambda\mu\)-mass below \(v\ge1\), and let \(B\) be the number of bad
leaves.  Then
\begin{equation}\label{eq:parity-cut-bound}
 \Var(R_k)\le
 \frac{C\kappa^2}{\log H_0}
 +CL^2\{J^{-2}+\varepsilon^2+e^{-2v}+h^2B\}.
\end{equation}
The first term bounds \(\Var(R_k\mid\pi)\) for every attainable finest
parity pattern \(\pi\).

In particular, if \(g>0\) area-almost everywhere on \(D\), then
\[
 \|R_k-\E R_k\|_2\longrightarrow0
 \quad\text{whenever}\quad \lambda/k^2\longrightarrow\infty,
\]
with \(k\) allowed to vary.
\end{lemma}

\begin{proof}
If \(\kappa=0\), then \(L=0\), and the growth bound gives \(R_k=0\).
We may therefore assume \(\kappa>0\).
First assume all leaf means are positive.
Conditional on \(\pi\), the leaf clouds are independent; the parity
of every block is fixed as the mod-two sum of its leaf parities.
Thus block replacement preserves that parity, and
\eqref{eq:parity-block-oscillation} applies.
Use the resampling form of Theorem~\ref{thm:multiscale} with area
masses, zero excess, and child-to-parent square-root mass ratio
\(2^{-H_0}\).
For each fixed attainable parity pattern it gives
\[
 \Var(R_k\mid\pi)\le
 \frac{C\kappa^2}{\log H_0}.
\]
The required integrability follows from \eqref{eq:parity-growth}
and the finite moments of a Poisson count. We next estimate the
variance of the conditional mean.

Write
\[
 g_D=\E[A_D\mid\pi],\qquad
 g_Q=\E[A_Q\mid\pi_Q],\qquad
 \E[R_k\mid\pi]=g_D-\sum_{Q\in\mathcal C_k}g_Q.
\]
The parity vectors in distinct root squares are independent, and
therefore so are the conditional means \(g_Q\). Hence
\begin{equation}\label{eq:parity-root-state-sum}
 \Var\!\left(\E[R_k\mid\pi]\right)
 \le2\Var(g_D)+2\sum_Q\Var(g_Q).
\end{equation}
Apply Lemma~\ref{lem:state-var} globally and in each root.
For a good leaf, \(r_i\le e^{-2v}\).
For two side-adjacent good leaves, their facing patches have combined
diameter at most \(C\varepsilon h\); hence
\eqref{eq:parity-two-patches} gives
\(b_e\le C(\varepsilon+e^{-2v})\).
For any edge, the direct one-label couplings give \(b_e\le5\), since
the changed labels lie in two side-adjacent squares.

There are \(O(h^{-2})\) grid edges, and at most \(4B\) edges have a bad
endpoint.  An edge appears at most once globally and once within its
root.  The same is true of leaf contributions.
The constant term in \eqref{eq:parity-state-variance}, summed over
the global square and the \(k^2\) roots, contributes at most
\[
 C(1+k^2)h^2\le CJ^{-2}.
\]
Consequently \eqref{eq:parity-root-state-sum} is bounded by
\(CL^2\{J^{-2}+\varepsilon^2+e^{-2v}+h^2B\}\).
The total variance identity gives \eqref{eq:parity-cut-bound}:
\[
 \Var(R_k)=\E\Var(R_k\mid\pi)
             +\Var\!\left(\E[R_k\mid\pi]\right)
\]

If some leaf has zero mean, add an independent Poisson process of
intensity \(\delta\,\mathrm dx\), keeping the grid and original
intensity fixed.  The added process is empty with probability tending
to one as \(\delta\downarrow0\).
The growth bound and Poisson moment bounds show that the perturbed
cuts converge in \(L^2\) to the original cut.
A patch that was good remains good after this perturbation, so the
number of bad leaves cannot increase.
Pass to the limit in \eqref{eq:parity-cut-bound}.
The conditional multiscale estimate itself applies directly to each
attainable parity pattern of the original process, including deterministic
empty leaves.

For the final assertion, put \(a=\varepsilon^2/9\).
If \(\lambda h^2at\ge2v\), a patch of mass below \(v\) must have at
least \(ah^2/2\) of its area in \(\{g<t\}\).  Otherwise its absolutely
continuous mass alone would be at least
\(\lambda t ah^2/2\ge v\).
Choose one deficient patch from each bad leaf.  Their interiors are
disjoint, and hence
\[
 h^2B\le2a^{-1}|\{g<t\}|.
\]
For fixed \(H_0,\varepsilon,v,t\), the premise eventually holds because
\(\lambda h^2=(\lambda/k^2)/J^2\to\infty\).
After taking this limit, send \(t\downarrow0\), then \(v\to\infty\),
then \(\varepsilon\downarrow0\), and finally \(H_0\to\infty\).
The assumption \(g>0\) almost everywhere gives
\(|\{g<t\}|\to0\), which proves the claimed centered cutting limit.
\end{proof}

\section{Fixed sampling laws and the covariance correction}\label{sec:fixed-law}

We now extend the homogeneous limit to fixed sampling laws whose
absolutely continuous component has a density positive almost everywhere
on the square. The proof compares centered costs in small cells with
their homogeneous counterparts. The law may have an arbitrary singular
component, and the density need not be bounded above, bounded away from
zero, or square integrable.

Let \(m\) be Lebesgue measure on \(D=[0,1]^2\), and write a Borel
probability as
\[
 \mu=g_\mu m+\tau_\mu,\qquad \tau_\mu\perp m.
\]
Fix a finite vector \(H\) whose coordinates are regular costs or parity
costs in Definitions~\ref{def:regular} and \ref{def:parity}.
Let \(H_{n,\mu}\) be its value at \(n\) iid points with law \(\mu\), and
let \(P_{\lambda,\mu}\) be its value at a Poisson process of intensity
\(\lambda\mu\), with total count \(N\).  Coincident labels are retained;
parity coordinates use the prescribed completion of odd samples.
Put
\begin{equation}\label{eq:fixed-law-vector}
 V_{\lambda,\mu}
 =\left(P_{\lambda,\mu}-\E P_{\lambda,\mu},
                    \frac{N-\lambda}{\sqrt\lambda}\right),
 \qquad
 a(\mu)=\int_D\sqrt{g_\mu}\,\mathrm dm.
\end{equation}
The coefficient \(a(\mu)\) is the Hellinger affinity between \(\mu\) and
the uniform law \(m\); the singular component contributes zero to the
affinity integral. Let \(\beta\) be the mean coefficient for uniform
sampling, and let \(\Sigma_P\) be the limiting covariance of the centered
costs under uniform Poisson sampling. Their existence is established
below. Define
\begin{equation}\label{eq:fixed-law-targets}
 \Gamma(a)=
 \begin{pmatrix}
  \Sigma_P&a\beta/2\\
  a\beta^\top/2&1
 \end{pmatrix},
 \qquad
 \Sigma(a)=\Sigma_P-\frac{a^2}{4}\beta\beta^\top .
\end{equation}
We write \(\mathcal N(0,B)\) for the centered Gaussian law with
covariance \(B\), including singular \(B\).

\begin{theorem}[Fixed-law transfer]\label{thm:fixed-law}
Suppose \(g_\mu>0\) \(m\)-almost everywhere on \(D\), and set
\(a=a(\mu)\).  Then
\begin{align}
 W_2\bigl(\Law(V_{\lambda,\mu}),\mathcal N(0,\Gamma(a))\bigr)
       &\longrightarrow0,
 &\Cov(V_{\lambda,\mu})&\longrightarrow\Gamma(a),
 \label{eq:fixed-law-poisson}\\
 W_2\bigl(\Law(H_{n,\mu}-\E H_{n,\mu}),
                         \mathcal N(0,\Sigma(a))\bigr)
       &\longrightarrow0,
 &\Cov(H_{n,\mu})&\longrightarrow\Sigma(a),
 \label{eq:fixed-law-binomial}\\
 \E H_{n,\mu}-\E P_{n,\mu}&\longrightarrow0.
 &&\label{eq:fixed-law-centers}
\end{align}
The fixed-count limit also holds when centered by \(\E P_{n,\mu}\).
The joint Poisson limit \eqref{eq:fixed-law-poisson}, including covariance
convergence, remains valid after conditioning on either parity of \(N\);
the centering remains the unconditional mean in
\eqref{eq:fixed-law-vector}.  Moreover,
\begin{equation}\label{eq:fixed-law-first-order}
 \frac{\E P_{\lambda,\mu}}{\sqrt\lambda}\longrightarrow a\beta,
 \qquad
 \frac{\E H_{n,\mu}}{\sqrt n}\longrightarrow a\beta,
 \qquad
 \frac{H_{n,\mu}}{\sqrt n}\longrightarrow a\beta
          \quad\text{in }L^2.
\end{equation}
No restriction is imposed on the singular measure \(\tau_\mu\).
The limiting covariance matrices need not be invertible.
\end{theorem}

\subsection{Centered cutting and the homogeneous constants}

Partition \(D\) into \(k^2\) squares \(C_j\) of side \(h=1/k\), assigning
shared faces by a fixed Borel rule, consistently under the refinements
used in the cutting arguments.  Every assigned cell is contained in
the corresponding closed geometric square and contains its interior.
Set \(p_j=\mu(C_j)\).  A local parity cost is completed at its own
square's lower-left corner.  Let \(R_k\) denote the global cost minus
the sum of these local costs.

\begin{corollary}[Centered cutting for a fixed law]\label{thm:fixed-cut}
If \(g_\mu>0\) \(m\)-almost everywhere on \(D\), then
\begin{equation}\label{eq:fixed-law-cutting}
 \|R_k-\E R_k\|_2\longrightarrow0
 \quad\text{whenever}\quad s=\lambda/k^2\longrightarrow\infty,
\end{equation}
with \(k\) allowed to vary.
\end{corollary}

\begin{proof}
For a regular coordinate this is Lemma~\ref{lem:regular-cut}.
For a parity coordinate, Lemma~\ref{lem:parity-cut} controls both the
expected conditional variance and the variance of the conditional mean
of the cutting error, given the leaf-count parities.
The squared \(L^2\)-norm of the centered vector remainder is the sum of
its coordinate variances.
\end{proof}

Under the uniform law \(m\), translation invariance and degree-one
homogeneity give the similarity relation for cell costs. Together with
the growth bounds and Corollary~\ref{thm:fixed-cut}, they verify the
assumptions of Theorem~\ref{thm:cut-clt}. Thus there are constants
\(\beta,\Sigma_P\) such that
\begin{equation}\label{eq:fixed-law-homogeneous}
 W_2\bigl(\Law(V_{s,m}),\mathcal N(0,\Gamma(1))\bigr)\longrightarrow0,
 \qquad
 \frac{\E P_{s,m}}{\sqrt s}\longrightarrow\beta .
\end{equation}
In particular the covariance matrices converge.
Lemmas~\ref{lem:regular-moments} and \ref{lem:parity-moments} also give
\(\sup_n\|H_{n,m}-\E P_{n,m}\|_2<\infty\), and hence
\begin{equation}\label{eq:fixed-law-uniform-mean}
 \frac{\E H_{n,m}}{\sqrt n}\longrightarrow\beta,\qquad
 \frac{H_{n,m}}{\sqrt n}\longrightarrow\beta
              \quad\text{in }L^2.
\end{equation}
Since \(\Gamma(1)\) is a covariance matrix, its Schur complement
satisfies \(\Sigma(1)\succeq0\). Also,
\[
 0\le a(\mu)\le\left(\int_Dg_\mu\,\mathrm dm\right)^{1/2}\le1.
\]
It follows that
\begin{equation}\label{eq:fixed-law-psd}
 \Sigma(a)=\Sigma(1)+\frac{1-a^2}{4}\beta\beta^\top\succeq0,
\end{equation}
and consequently \(\Gamma(a)\succeq0\).

\subsection{Approximation of Hellinger affinity}

Within each cell, we compare the restriction of \(\mu\) with the uniform
measure having the same total mass. For the partition above, define
\begin{equation}\label{eq:fixed-law-affinity-data}
 a_k=\sum_jh\sqrt{p_j},\qquad
 v_j=\frac{p_j}{h^2},\qquad
 u_j=\tau_\mu(C_j)+
       \int_{C_j}\bigl(\sqrt{g_\mu}-\sqrt{v_j}\bigr)^2\,\mathrm dm.
\end{equation}
The quantity \(a_k=\sum_j\sqrt{\mu(C_j)m(C_j)}\) is the Hellinger
affinity between the discrete distributions
\((\mu(C_j))_j\) and \((m(C_j))_j\). The quantity \(u_j\) is the squared
Hellinger distance between the finite measures \(\mu|_{C_j}\) and
\(v_jm|_{C_j}\), including the singular mass. If \(p_j=0\), set
\(u_j/\sqrt{v_j}=0\). Since
\(\int_{C_j}g_\mu\,\mathrm dm+\tau_\mu(C_j)=p_j\), expansion gives
\[
 u_j=2p_j-2\sqrt{v_j}\int_{C_j}\sqrt{g_\mu}\,\mathrm dm,
 \qquad 0\le u_j\le2p_j.
\]
Therefore
\begin{equation}\label{eq:fixed-law-affinity-identity}
 2\{a_k-a(\mu)\}=\sum_j\frac{u_j}{\sqrt{v_j}}\ge0.
\end{equation}
These quantities are finite for every Borel probability \(\mu\).

\begin{lemma}[Hellinger affinity under cell partitions]\label{lem:affinity}
For every fixed Borel probability \(\mu\) on \(D\),
\[
 a_k\longrightarrow a(\mu)\qquad(k\to\infty).
\]
This holds for every fixed rule assigning shared faces.
\end{lemma}

\begin{proof}
For a strictly positive continuous function \(\varphi\) on \(D\), put
\[
 A_\varphi=\frac12\left(\int_D\varphi\,\mathrm dm
                  +\int_D\varphi^{-1}\,\mathrm d\mu\right).
\]
We first prove the variational identity
\begin{equation}\label{eq:fixed-law-affinity-variational}
 a(\mu)=\inf_{\varphi\in C(D),\,\varphi>0}A_\varphi.
\end{equation}
Relative to \(\rho=m+\mu\), let \(b=\mathrm dm/\mathrm d\rho\) and
\(e=\mathrm d\mu/\mathrm d\rho\).
The pointwise infimum of \((bx+e/x)/2\) over \(x>0\) is \(\sqrt{be}\),
and \(\int\sqrt{be}\,\mathrm d\rho=a(\mu)\).
First minimize with \(x\) constrained to \([L^{-1},L]\), \(L\ge1\).
The minimum is attained by a measurable function in that interval;
the minimized values decrease to \(\sqrt{be}\) and are bounded by
\((b+e)/2\).
For each fixed \(L\), approximate the measurable minimizer in
\(L^1(\rho)\) by continuous functions and clip them to \([L^{-1},L]\).
The map \(x\mapsto x^{-1}\) is Lipschitz on this interval, so the
corresponding values of \(A_\varphi\) converge to the minimum over
measurable functions.
Dominated convergence as \(L\to\infty\) proves
\eqref{eq:fixed-law-affinity-variational}; the reverse inequality also
follows directly from the pointwise arithmetic--geometric mean
inequality.

Choose any \(z_j\) in the closed geometric square corresponding to
\(C_j\).  Arithmetic--geometric mean gives
\[
 h\sqrt{p_j}\le
 \tfrac12\{h^2\varphi(z_j)+p_j/\varphi(z_j)\}.
\]
After summing over \(j\), the first term on the right converges to
\(\int\varphi\,\mathrm dm\), and the second to
\(\int\varphi^{-1}\,\mathrm d\mu\).
For the latter, the error is bounded by the modulus of continuity of
\(\varphi^{-1}\) at \(\sqrt2h\): every point assigned to \(C_j\) lies in
its closed square. This also covers measures that charge partition
faces. Hence \(\limsup_k a_k\le A_\varphi\).
Take the infimum over \(\varphi\) and use
\eqref{eq:fixed-law-affinity-identity}.
\end{proof}

\subsection{Comparison of centered cell costs}

Write \(s\) for the average number of points per cell and \(s_j\) for
the mean count in cell \(C_j\). Define
\[
 s=\lambda h^2,\qquad s_j=\lambda p_j,\qquad
 \Delta=h^2\sum_j\min\{1,\lambda u_j\},
 \qquad
 B_M=\sum_{s_j<M}(h^2+p_j).
\]
The quantity \(\Delta\) controls the comparison of the two input laws,
while \(B_M\) bounds the contribution from cells whose mean count is
less than \(M\). The homogeneous Gaussian approximation error is
\[
 \epsilon(M)=\sup_{t\ge M}
 W_2\bigl(\Law(V_{t,m}),\mathcal N(0,\Gamma(1))\bigr).
\]
By \eqref{eq:fixed-law-homogeneous}, \(\epsilon(M)\to0\).
For a matrix, \(\|\cdot\|_{\max}\) denotes the largest absolute entry.

\begin{lemma}[Centered cell comparison]\label{lem:cell-compare}
For \(s\ge1\) and sufficiently large \(M\), set
\(r_k=\|R_k-\E R_k\|_2\). Then
\begin{align}
 W_2\bigl(\Law(V_{\lambda,\mu}),\mathcal N(0,\Gamma(a_k))\bigr)
 &\le r_k+C\{\Delta^{1/4}+B_M^{1/2}+\epsilon(M)\},
 \label{eq:fixed-law-cell-wasserstein}\\
 \|\Cov(V_{\lambda,\mu})-\Gamma(a(\mu))\|_{\max}
 &\le C\{r_k+r_k^2+\sqrt\Delta+B_M+\epsilon(M)\}.
 \label{eq:fixed-law-cell-covariance}
\end{align}
The constants depend only on the fixed vector of costs.
\end{lemma}

\begin{proof}
The first step is to compare the same vector statistic under two
input laws. Let \(X=F(\xi)\) and \(Y=F(\eta)\), where \(\xi\) and \(\eta\)
have laws \(P\) and \(Q\), respectively. Write
\[
 m_P=\bigl(\E|X-\E X|^4\bigr)^{1/4},\qquad
 m_Q=\bigl(\E|Y-\E Y|^4\bigr)^{1/4}.
\]
Let \(p,q\) be densities of \(P,Q\) relative to a common dominating
measure, and set
\(d_{\mathrm H}(P,Q)=\|\sqrt p-\sqrt q\|_2\).
There is a maximal coupling of \(P\) and \(Q\) for which
\begin{equation}\label{eq:fixed-law-centered-hellinger}
 \begin{split}
 &\E\bigl|(X-\E X)-(Y-\E Y)\bigr|^2\\
 &\hspace{2em}
 +\|\Cov X-\Cov Y\|_{\max}
 \le C(m_P^2+m_Q^2)d_{\mathrm H}(P,Q).
 \end{split}
\end{equation}

To prove the estimate, construct the maximal coupling from the common
density \(\min(p,q)\) and denote by \(B\) the event that the inputs are
drawn from the residual densities. Conditional on the first input,
\[
 \Pp(B\mid\xi)=t_P=(1-q/p)_+
 \quad\text{on }\{p>0\}.
\]
The pointwise inequality \(\sqrt p\,t_P\le2|\sqrt p-\sqrt q|\) gives
\[
 \|t_P\|_{2,P}\le2d_{\mathrm H}(P,Q),\qquad
 \Pp(B)=d_{\mathrm{TV}}(P,Q)\le d_{\mathrm H}(P,Q).
\]
For an independent copy \((X',Y',B')\), Cauchy--Schwarz gives
\[
 \begin{split}
 \E[|X-X'|^2\1_B]
 &\le2\E[|X-\E X|^2t_P]
          +2\E|X-\E X|^2\Pp(B)\\
 &\le C m_P^2d_{\mathrm H}(P,Q),
 \end{split}
\]
The same bound holds with \(B'\), and with the second input and \(m_Q\).
Outside \(B\cup B'\), we have \(X=Y\) and \(X'=Y'\).
On \(B\cup B'\), both
\(\lvert(X-X')-(Y-Y')\rvert^2\) and the absolute value of each entry of
\[
 (X-X')(X-X')^\top-(Y-Y')(Y-Y')^\top
\]
are bounded by \(C\{|X-X'|^2+|Y-Y'|^2\}\).
The independent-copy identities
\[
 \E|(X-\E X)-(Y-\E Y)|^2
 =\tfrac12\E|(X-X')-(Y-Y')|^2
\]
and
\(\Cov X=\tfrac12\E[(X-X')(X-X')^\top]\)
prove \eqref{eq:fixed-law-centered-hellinger}.

Now compare the actual Poisson input in cell \(C_j\), of intensity
\(\lambda\mu|_{C_j}\), with the homogeneous input of intensity
\(\lambda v_jm|_{C_j}\).
Both have total mean \(s_j\).
If \(p_j>0\), the Hellinger affinity between their conditional one-point
laws is
\[
 \alpha_j=\frac1{p_j}
          \int_{C_j}\sqrt{g_\mu v_j}\,\mathrm dm.
\]
Expanding over the common Poisson count gives the Hellinger affinity
between the two Poisson configuration laws:
\begin{equation}\label{eq:fixed-law-poisson-affinity}
 \sum_{n\ge0}e^{-s_j}\frac{s_j^n}{n!}\alpha_j^n
 =e^{-s_j(1-\alpha_j)}
 =e^{-\lambda u_j/2}.
\end{equation}
This computation includes the singular part through the term
\(\tau_\mu(C_j)\) in \(u_j\).
It is the equal-mass case of the Poisson Hellinger-affinity identity;
see \citet{Leskela2024} for the general formulation.
To couple the two Poisson inputs, draw a common count and, conditional
on that count, use a maximal coupling of the two configuration laws.
The count strata are disjoint and have the same probabilities under
both laws. Thus the resulting coupling is also maximal for the full
Poisson inputs, and the two counts agree almost surely.
If \(P_j\) and \(Q_j\) denote these two Poisson input laws, then
\eqref{eq:fixed-law-poisson-affinity} gives
\[
 d_{\mathrm H}(P_j,Q_j)^2
 =2\bigl(1-e^{-\lambda u_j/2}\bigr)
 \le2d_j,\qquad d_j=\min\{1,\lambda u_j\}.
\]

First retain the cells with \(s_j\ge M\), where \(M\ge1\). Let
\(Y_j\) and \(\widetilde Y_j\) be the centered cost vectors of the two
coupled Poisson inputs in \(C_j\), with the actual and homogeneous laws,
respectively. Parity coordinates are completed at the cell's prescribed
corner. After rescaling the square to \(D\), the corresponding centered
vectors are \(Y_j/h\) and \(\widetilde Y_j/h\).
Lemmas~\ref{lem:regular-moments} and \ref{lem:parity-moments} bound the
\(L^4\) norms of \(Y_j/h\) and \(\widetilde Y_j/h\) uniformly over these
cells and their sampling laws. The centered count divided by
\(\sqrt{s_j}\) also has bounded \(L^4\) norm, since its fourth moment is
\(3+s_j^{-1}\). In the global vector, the weights are \(h\) for the
unit-square costs and \(\sqrt{p_j}\) for the standardized cell count.

Take the couplings independently across cells. The paired count
coordinates agree, and the centered vectors \(Y_j-\widetilde Y_j\)
are independent. Equation~\eqref{eq:fixed-law-centered-hellinger} yields
\[
 \E\left|\sum_{s_j\ge M}(Y_j-\widetilde Y_j)\right|^2
 \le Ch^2\sum_j\sqrt{d_j}\le C\sqrt\Delta.
\]
The last inequality follows from Cauchy--Schwarz and
\(\sum_jh^2=1\). The same estimate bounds the difference between the
covariance matrices of the two sums of cost vectors.
For a cost--count entry, the corresponding bound is
\[
 C\sum_j h\sqrt{p_jd_j}
 \le C\left(\sum_jh^2d_j\right)^{1/2}
       \left(\sum_jp_j\right)^{1/2}
 =C\sqrt\Delta.
\]

For each retained homogeneous cell, couple its unit-square
cost--count vector \(V_{s_j,m}\) to a Gaussian with covariance
\(\Gamma(1)\), with \(L^2\)-error at most \(\epsilon(M)\).
Take these couplings independently and apply
\[
 D_j=\operatorname{diag}(hI,\sqrt{p_j}).
\]
Their total squared error is at most \(2\epsilon(M)^2\), because
\[
 \sum_j\|D_j\|_{\mathrm{op}}^2
 \le\sum_j(h^2+p_j)=2.
\]
The covariance error is \(O(\epsilon(M))\), by Cauchy--Schwarz and
the bounded second moments of the centered sums.
The full weighted Gaussian sum, including all cells, has covariance
\begin{equation}\label{eq:fixed-law-gaussian-sum}
 \sum_jD_j\Gamma(1)D_j^\top
 =
 \begin{pmatrix}
  \Sigma_P&(\sum_jh\sqrt{p_j})\beta/2\\
  (\sum_jh\sqrt{p_j})\beta^\top/2&1
 \end{pmatrix}
 =\Gamma(a_k).
\end{equation}

Omitting the centered costs and counts from cells with \(s_j<M\),
or omitting their Gaussian counterparts, changes the sum by
\(O(\sqrt{B_M})\) in \(L^2\) and its covariance by \(O(B_M)\).
Indeed, the rescaled centered costs have bounded second moments at
intensities at least one by the moment lemmas.
Below intensity one, the growth bounds and finite Poisson moments give
the same boundedness directly.
The weighted count variance is \(p_j\), and a cost--count covariance
is bounded by \(C h\sqrt{p_j}\le C(h^2+p_j)\).
Zero-mass cells have zero actual cost and count; their normalized count
is defined to be zero.
Their Gaussian cost contribution is still covered by the \(h^2\)
part of \(B_M\).

The sum of the actual cell-count coordinates is exactly
\((N-\lambda)/\sqrt\lambda\).
Adding the centered cutting remainder changes the vector by \(r_k\)
in \(L^2\), and its covariance by \(O(r_k+r_k^2)\).
These comparisons prove
\eqref{eq:fixed-law-cell-wasserstein} and the covariance assertion with
\(\Gamma(a_k)\) in place of \(\Gamma(a(\mu))\).

It remains to compare the covariance targets \(\Gamma(a_k)\) and
\(\Gamma(a(\mu))\).
On cells with \(\lambda u_j\ge1\), the total area is at most \(\Delta\).
Since \(u_j\le2p_j\), their contribution to
\(a_k-a(\mu)=\tfrac12\sum_j u_j/\sqrt{v_j}\) is at most
\[
 \sum_{\lambda u_j\ge1}h\sqrt{p_j}\le\sqrt\Delta.
\]
On the remaining cells,
\(u_j/\sqrt{v_j}\le\sqrt2h\sqrt{u_j}\), and
\[
 \sum_{\lambda u_j<1}u_j
 \le\frac{\Delta}{s}.
\]
Another application of Cauchy--Schwarz gives
\begin{equation}\label{eq:fixed-law-target-affinity-error}
 0\le a_k-a(\mu)
 \le\sqrt\Delta+\sqrt{\frac{\Delta}{2s}}.
\end{equation}
Since \(s\ge1\), this difference is \(O(\sqrt\Delta)\).
Equation \eqref{eq:fixed-law-targets} proves
\eqref{eq:fixed-law-cell-covariance}.
\end{proof}

\subsection{Proof of the fixed-law transfer}

\begin{proof}[Proof of Theorem~\ref{thm:fixed-law}]
We first choose a partition for which the centered cutting error and
the comparison of cell laws both vanish. Write
\[
 w(t)=m\{g_\mu<t\},\qquad L_k=a_k-a(\mu).
\]
Then \(w(t)\to0\) as \(t\downarrow0\), by positivity of the density,
and \(L_k\to0\), by Lemma~\ref{lem:affinity}.
Choose an increasing sequence of integers \(K_i\ge i\) such that
\(L_k\le i^{-3}\) for all \(k\ge K_i\), and make the thresholds
\(4iK_i^2\) strictly increasing.
For
\[
 4iK_i^2\le\lambda<4(i+1)K_{i+1}^2,
 \qquad
 k=\left\lfloor\sqrt{\lambda/i}\right\rfloor,
\]
we have \(k\ge K_i\), \(k\to\infty\), and
\[
 i\le s=\lambda/k^2\le4i,\qquad
 sL_k\le4i^{-2}\longrightarrow0.
\]
Thus the cells become small while their average Poisson count tends
to infinity, and \(sL_k\to0\) holds without a prescribed rate for
\(L_k\to0\).

Put \(b_j=u_j/h^2\).
Since \(\min(1,x)\le x^{2/3}\), H\"older's inequality gives
\begin{equation}\label{eq:fixed-law-delta-vanishing}
 \begin{split}
 \Delta
 &\le s^{2/3}\sum_jh^2b_j^{2/3}\\
 &\le s^{2/3}
       \left(\sum_{p_j>0}h^2\frac{b_j}{\sqrt{v_j}}\right)^{2/3}
       \left(\sum_jh^2v_j\right)^{1/3}\\
 &=\{2sL_k\}^{2/3}\longrightarrow0.
 \end{split}
\end{equation}
The terms with \(p_j=0\) have \(b_j=0\).
Hence convergence of the affinities suffices to control \(\Delta\);
no \(L^2\) approximation of \(g_\mu\) is needed.

For fixed \(M\) and \(t>0\), suppose \(st\ge2M\).
If \(s_j<M\), at least half the area of \(C_j\) must lie in
\(\{g_\mu<t\}\); otherwise
\[
 s_j=\lambda p_j
 \ge\lambda\int_{C_j}g_\mu\,\mathrm dm
 \ge st/2\ge M.
\]
These cells have total area at most \(2w(t)\).
Their total mass is at most \(k^2M/\lambda=M/s\).
Thus
\[
 B_M\le2w(t)+M/s\longrightarrow0
\]
after first sending \(s\to\infty\) and then \(t\downarrow0\).
Corollary~\ref{thm:fixed-cut} gives \(r_k\to0\).
Lemma~\ref{lem:cell-compare}, first with \(M\) fixed and then as
\(M\to\infty\), reduces the joint Poisson limit to convergence of the
Gaussian targets.
Since \(a_k\to a(\mu)\), the positive semidefinite square roots
\(\Gamma(a_k)^{1/2}\) converge to \(\Gamma(a(\mu))^{1/2}\).
Applying these matrices to the same standard Gaussian vector proves
convergence in \(W_2\), also for singular limiting covariances.
This proves \eqref{eq:fixed-law-poisson} and its covariance assertion.

We next show that this joint limit persists after conditioning on
either count parity.
Let \(z\) be the common completion corner, choosing the lower-left
corner also when there is no parity coordinate.
For every \(\eta>0\), the ball \(B(z,\eta)\) has positive \(\mu\)-mass:
its intersection with \(D\) has positive area and \(g_\mu>0\) almost
everywhere there. Lemma~\ref{lem:neighbor-count} therefore gives
\begin{equation}\label{eq:fixed-law-adjacent-vanishing}
 \sup_{\lambda/2\le j\le2\lambda}
 \dBL\bigl(\Law(H_{j+1,\mu}),
          \Law(H_{j,\mu})\bigr)\longrightarrow0.
\end{equation}
Indeed, first fix \(\eta>0\), let \(\lambda\to\infty\) in
\eqref{eq:parity-step-one-bl}, and then let \(\eta\downarrow0\).
By translation invariance of \(\dBL\), the same estimate holds after
subtracting any common deterministic center.

Let \(w_j=e^{-\lambda}\lambda^j/j!\).
Unimodality and the maximal-atom estimate give
\[
 \sum_{j\ge0}|w_{j+1}-w_j|=O(\lambda^{-1/2}).
\]
For fixed \(t\in\R^r\) and \(u\in\R\), set
\[
 \chi_j=\E\exp\left\{
  i t^\top(H_{j,\mu}-\E P_{\lambda,\mu})
       + i u\,\frac{j-\lambda}{\sqrt\lambda}\right\}.
\]
Equation \eqref{eq:fixed-law-adjacent-vanishing} and the deterministic
count increment \(\lambda^{-1/2}\) give
\(\sup_{\lambda/2\le j\le2\lambda}|\chi_{j+1}-\chi_j|\to0\).
Pairing consecutive terms yields
\[
 \begin{split}
 \E[(-1)^N e^{i(t,u)\cdot V_{\lambda,\mu}}]
 &=\sum_{\ell\ge0}
       w_{2\ell}(\chi_{2\ell}-\chi_{2\ell+1})\\
 &\quad+\sum_{\ell\ge0}
       (w_{2\ell}-w_{2\ell+1})\chi_{2\ell+1}
 \longrightarrow0.
 \end{split}
\]
The Poisson probability outside \([\lambda/2,2\lambda]\) tends to
zero, and the second sum is bounded by the displayed total variation
of the weights.
Since
\[
 \Pp(N\ \text{even})=\tfrac12(1+e^{-2\lambda}),
 \qquad
 \Pp(N\ \text{odd})=\tfrac12(1-e^{-2\lambda}),
\]
both conditional characteristic functions converge to the same
Gaussian characteristic function as the unconditional one.
The moment lemmas give uniformly bounded centered cost moments of
every fixed order, and normalized Poisson counts have bounded moments
of every fixed order for \(\lambda\ge1\).
Both parity probabilities stay bounded away from zero. The same
moment bounds therefore hold after conditioning, and the conditional
squared norms are uniformly integrable. This upgrades both conditional
limits to \(W_2\) convergence and gives convergence of their means and
covariances.

Now take \(\lambda=n\) and condition on
\[
 D_n=\{N\equiv n\pmod2\}.
\]
The conditional joint Gaussian limit just proved has cost--count
covariance \(d=a\beta/2\).
For counts in the same parity class and
\(|j-n|\le\delta\sqrt n\), telescoping
\eqref{eq:parity-step-two-bl} gives
\begin{equation}\label{eq:fixed-law-count-window}
 \sup_{\substack{|j-n|\le\delta\sqrt n\\j\equiv n\pmod2}}
 \dBL\bigl(\Law(H_{j,\mu}-\E P_{n,\mu}),
          \Law(H_{n,\mu}-\E P_{n,\mu})\bigr)
 \le C_H\delta+o(1).
\end{equation}
The intermediate counts are at least \(n/2\) for large \(n\).
Apply Lemma~\ref{lem:count-window} with \(q=2\).
The resulting fixed-count weak limit is Gaussian with covariance
\[
 \Sigma_P-dd^\top
   =\Sigma_P-\frac{a^2}{4}\beta\beta^\top=\Sigma(a).
\]
Lemmas~\ref{lem:regular-moments} and \ref{lem:parity-moments} provide a
uniform moment of order \(p>2\), both at the fixed-count mean and at
\(\E P_{n,\mu}\).
The final assertion of Lemma~\ref{lem:count-window} now yields
\(W_2\) convergence for both centerings, convergence of the covariance
matrices, and
\(\E H_{n,\mu}-\E P_{n,\mu}\to0\).
This proves \eqref{eq:fixed-law-binomial}--\eqref{eq:fixed-law-centers}.
Since \(D_n\) uses the parity of \(n\), this proves the limit along
both even and odd sample sizes.

Finally, we recover the first-order mean from the covariance with the
Poisson count. Set
\[
 \gamma(t)=
 \Cov\left(P_{t,\mu},\frac{N-t}{\sqrt t}\right).
\]
The deterministic growth bounds justify differentiating the Poisson
series term by term on compact subintervals of \((0,\infty)\):
\begin{equation}\label{eq:fixed-law-mean-derivative}
 \frac{\mathrm d}{\mathrm dt}\E P_{t,\mu}
 =\frac1t\,\E[(N-t)P_{t,\mu}]
 =t^{-1/2}\gamma(t).
\end{equation}
By \eqref{eq:fixed-law-poisson},
\(\gamma(t)\to a\beta/2\).
Integrate \eqref{eq:fixed-law-mean-derivative} from any fixed
\(M\ge1\) to \(\lambda\) and divide by \(\sqrt\lambda\).
Since \(\int_M^\lambda t^{-1/2}\,\mathrm dt
=2(\sqrt\lambda-\sqrt M)\), the limit is \(a\beta\).
Equation \eqref{eq:fixed-law-centers} gives the fixed-count mean
limit.  The uniform centered second moment bound then gives
\(H_{n,\mu}/\sqrt n\to a\beta\) in \(L^2\).
\end{proof}

Changing the sampling law adds a positive semidefinite correction of
rank at most one, proportional to \(\beta\beta^\top\):
\[
 \Sigma(a)-\Sigma(1)
   =\frac{1-a^2}{4}\beta\beta^\top\succeq0.
\]
Thus the uniform law minimizes the fixed-count limiting covariance,
in the positive semidefinite order, among the fixed laws of
Theorem~\ref{thm:fixed-law}.  The deterministic null-space criterion
below identifies which of these Gaussian limits are nondegenerate.

\section{A deterministic characterization of covariance null directions}
\label{sec:nondegeneracy}

We now characterize the null space of the limiting covariance in terms
of the values of the costs on finite configurations. The proof starts
with a lower bound obtained from two local configurations. Comparing a
local cost difference with a boundary error and then summing orthogonal
projections is the principle used by
\citet[Proposition~3.1 and Lemma~4.1]{WehrAizenman1990Extensive}.
We use configurations of the same cardinality, which allows us to
condition on the cell counts before applying this argument.

Let $H=(H_1,\ldots,H_r)$ be a fixed vector of regular and parity costs as in Definitions~\ref{def:regular} and \ref{def:parity}. Let $\mathcal A$ be the class of all finite labeled configurations in $D$ if every coordinate is regular, and the class of all even labeled configurations in $D$ otherwise. On $\mathcal A$ we always use the original costs, without a completion label. For each regular coordinate, let $B_j$ be the constant $B_0$ in
\eqref{eq:regular-full-boundary}; for each parity coordinate, let it be
the constant $\kappa$ in \eqref{eq:parity-full-boundary}. Thus $B_j$ is
the error bound for that coordinate on the unit square. For
$t\in\R^r$, write
\[
 H_t=t^\top H,
 \qquad B(t)=\sum_{j=1}^r|t_j|B_j.
\]
The square side lengths in these comparisons are Euclidean, even when
the costs are defined using other norms.

\begin{proposition}[Two local configurations]\label{prop:two-config}
Suppose that $S_-,S_+\in\mathcal A$ have the same cardinality $q\ge1$ and satisfy
\begin{equation}\label{eq:two-config-gap}
 H_t(S_+)-H_t(S_-)>2B(t).
\end{equation}
For every $c>0$ there are $v>0$ and $n_0<\infty$ such that
\begin{equation}\label{eq:two-config-variance}
 \Var(t^\top H_{n,\mu})\ge v
 \qquad(n\ge n_0,\quad\mu\ge cm),
\end{equation}
where $\mu$ is a probability on $D$. At odd counts, parity coordinates use the prescribed corner completion.
\end{proposition}

\begin{proof}
By continuity at fixed cardinality, there are relative open neighborhoods $U_-,U_+\subset D^q$ and $\delta>0$ such that
\[
 H_t(s_+)-H_t(s_-)\ge2B(t)+\delta,
 \qquad s_-\in U_-,\quad s_+\in U_+.
\]
Their probabilities under the uniform product measure, denoted by
$\ell_-$ and $\ell_+$, are positive. This also holds when $S_-$ or $S_+$
has coincident locations or labels on the boundary of $D$.

Put $k=\lfloor\sqrt{n/q}\rfloor$ and $h=1/k$, and divide $D$ into $k^2$ assigned cells $C_i$. Let $N_i$ be their sample counts and $I_i$ their within-cell configurations. Conditional on the count vector $\mathbf N$, the $I_i$ are independent, with iid locations in each cell according to the corresponding conditional one-point law. Write $Y=t^\top H_{n,\mu}$.

On $\{N_i=q\}$, homogeneity and the boundary comparisons
\eqref{eq:regular-full-boundary} and \eqref{eq:parity-full-boundary} give
\begin{equation}\label{eq:two-config-projection}
 \E[Y\mid\mathbf N,I_i]
 =hH_t(\widetilde I_i)+b_i(\mathbf N)+e_i(\mathbf N,I_i),
 \qquad |e_i|\le hB(t),
\end{equation}
where $\widetilde I_i$ is the configuration $I_i$ rescaled to $D$.
For this identity, apply the boundary comparison to each coordinate
and then integrate over the exterior. For a regular cost, add the
cell's boundary point to the exterior configuration. If a parity cost
occurs, $q$ is even, so the interior needs no completion. The global
completion label, when present, lies outside the open cell. Conditional
on $\mathbf N$, the resulting exterior costs are independent of the
locations in $I_i$ and contribute to $b_i(\mathbf N)$. The continuous
extensions allow labels on cell boundaries.

Suppose that $\widetilde I_i$ belongs to $U_-$ and $U_+$ with conditional probabilities at least $p_-$ and $p_+$, respectively. By \eqref{eq:two-config-projection}, the corresponding values of the conditional projection differ by at least $\delta h$. Using two independent copies of that projection gives
\begin{equation}\label{eq:two-config-projection-lower}
 \Var\bigl(\E[Y\mid\mathbf N,I_i]\mid\mathbf N\bigr)
 \ge p_-p_+\delta^2h^2
 \quad\text{on }\{N_i=q\}.
\end{equation}
Conditional independence makes the centered single-cell projections orthogonal. Their sum is the first-order orthogonal projection of $Y-\E[Y\mid\mathbf N]$ onto the sum of the single-cell subspaces. Consequently,
\[
 \Var(Y\mid\mathbf N)
 \ge\sum_i\Var\bigl(\E[Y\mid\mathbf N,I_i]\mid\mathbf N\bigr).
\]

Let $p_i=\mu(C_i)$. Since $\sum_i p_i=1$, at least half the cells satisfy
\[
 ch^2\le p_i\le2h^2.
\]
On these cells the conditional one-point law, after rescaling, dominates $(c/2)m|_D$. We may therefore take
\[
 p_-=(c/2)^q\ell_-,\qquad p_+=(c/2)^q\ell_+.
\]
Moreover, $nh^2\to q$. The binomial probabilities $\Pp(N_i=q)$ are bounded below by a positive constant $q_0=q_0(c,q)$, uniformly over these cells and all sufficiently large $n$: their parameters satisfy $np_i$ in a fixed compact subinterval of $(0,\infty)$ and $p_i\to0$. Taking expectations in \eqref{eq:two-config-projection-lower} and summing yields
\[
 \Var Y\ge \frac{k^2}{2}\,q_0p_-p_+\delta^2h^2
 =\frac{q_0p_-p_+\delta^2}{2}.
\]
This proves \eqref{eq:two-config-variance} uniformly over the indicated probability laws.
\end{proof}

Write $\beta$ for the homogeneous first-order mean vector and
\[
 S_1=\Sigma_P-\frac14\beta\beta^\top
\]
for the homogeneous fixed-count limiting covariance from Theorem~\ref{thm:fixed-law}.

\begin{theorem}[Deterministic null-space criterion]\label{thm:kernel}
For every $t\in\R^r$,
\begin{equation}\label{eq:kernel-characterization}
 t\in\ker S_1
 \quad\Longleftrightarrow\quad
 \sup_{S\in\mathcal A}|t^\top H(S)|<\infty.
\end{equation}
If these conditions hold, the oscillation of $H_t$ at every admissible cardinality is at most $2B(t)$, and
\begin{equation}\label{eq:kernel-uniform-bound}
 |H_t(S)|\le2B(t),\qquad S\in\mathcal A.
\end{equation}
Furthermore, $\beta\perp\ker S_1$. For every sampling law in Theorem~\ref{thm:fixed-law}, its fixed-count limiting covariance has the same null space $\ker S_1$.
\end{theorem}

\begin{proof}
If $t\in\ker S_1$, covariance convergence under the uniform law gives
$\Var(t^\top H_{n,m})\to0$. Proposition~\ref{prop:two-config} therefore rules out two admissible configurations of the same cardinality whose costs differ by more than $2B(t)$. A configuration with all labels at one location has cost zero, by translation invariance and homogeneity. Comparing with this configuration proves \eqref{eq:kernel-uniform-bound}; the empty cost is also zero.

Conversely, suppose that $H_t$ is bounded on $\mathcal A$. Along admissible sample sizes, the centered variables $H_t-\E H_t$ have support in a common bounded interval. Their Gaussian limit from Theorem~\ref{thm:fixed-law} must therefore have zero variance. This proves $t\in\ker S_1$.

For such a direction, divide the uniformly bounded expectations by $\sqrt n$ along admissible counts. The homogeneous mean limit gives $t^\top\beta=0$. For a general law $\mu$ covered by Theorem~\ref{thm:fixed-law}, put
\[
 a(\mu)=\int_D\sqrt{g_\mu}\,dm.
\]
Its fixed-count limiting covariance is
\begin{equation}\label{eq:kernel-general-covariance}
 S_1+\frac{1-a(\mu)^2}{4}\beta\beta^\top.
\end{equation}
Cauchy--Schwarz gives $a(\mu)^2\le\int_Dg_\mu\,dm\le1$. Both summands in \eqref{eq:kernel-general-covariance} are positive semidefinite. A null vector of their sum is therefore a null vector of $S_1$; conversely, every null vector of $S_1$ is orthogonal to $\beta$ and annihilates the second summand. This proves equality of the null spaces.
\end{proof}

The theorem also identifies directions in the range of $S_1$ from
sequences of deterministic configurations. If $b_j\to\infty$ and
\begin{equation}\label{eq:kernel-range-witness}
 \frac{H(S_j)}{b_j}\longrightarrow v,
 \qquad S_j\in\mathcal A,
\end{equation}
then \eqref{eq:kernel-uniform-bound} implies $t^\top v=0$ for every $t\in\ker S_1$, and hence $v\in\operatorname{range}(S_1)$. The same conclusion holds for $\beta$. In particular, a scalar cost with nonzero homogeneous mean coefficient has strictly positive limiting variance for every law in Theorem~\ref{thm:fixed-law}.

\begin{corollary}[Separation of regular and parity null directions]
\label{cor:split-kernel}
Suppose $H=(H_R,H_P)$, where $H_R$ consists of regular costs and $H_P$ consists of parity costs. Let $S_{RR}$ and $S_{PP}$ be the corresponding homogeneous fixed-count principal covariance blocks. Then
\begin{equation}\label{eq:split-kernel}
 \ker S_1=\ker S_{RR}\times\ker S_{PP},
 \qquad
 \operatorname{rank}(S_1)=\operatorname{rank}(S_{RR})+
 \operatorname{rank}(S_{PP}).
\end{equation}
\end{corollary}

\begin{proof}
Suppose that $t_R^\top H_R+t_P^\top H_P$ is bounded on even configurations. Given any configuration $S$, double every label. By repetition invariance, the regular vector on the doubled configuration is $H_R(S)$. By the coincident-pair identity in Lemma~\ref{lem:parity-insert}, every parity cost is zero on the doubled configuration: start from the empty configuration and insert its coincident pairs. Thus $t_R^\top H_R$ is bounded on all configurations. Subtracting this bounded function from the original combination shows that $t_P^\top H_P$ is bounded on even configurations. The converse follows by addition. Apply Theorem~\ref{thm:kernel} to the full vector and separately to its two blocks. The rank identity follows from the dimensions of the null spaces.
\end{proof}

The rank identity is compatible with a nonzero cross-covariance block
between the regular and parity costs.

\section{Tours, perfect matchings, and the Held--Karp relaxation}
\label{sec:applications}

We now verify the boundary comparisons for three optimization problems. A tour must visit all sites in a single cycle; cutting the tour and reconnecting it along a square boundary gives the required comparison. A perfect matching pairs every label with exactly one other label. Cutting the sample into cells can leave odd numbers of labels, so we complete the local configurations at the boundary. For Poisson input, we condition on the cell-count parities when estimating the cutting error, as in Section~\ref{sec:parity}. The Held--Karp relaxation uses fractional edge capacities subject to degree and cut constraints. For this problem, we add capacity along the boundary and then remove the auxiliary labels. These constructions give the centered cutting estimates needed for the joint Gaussian limits below.

\begin{theorem}[Joint optimization limits]\label{thm:applications}
Fix a finite vector $H$ whose coordinates are optimal tour lengths, minimum perfect-matching lengths, or Held--Karp values, each computed using a norm on $\R^2$. Evaluate all coordinates on the same sample. Let $\mu$ be a fixed Borel probability on $D$ with Lebesgue decomposition $\mu=g\,m+\tau$, where $\tau\perp m$ and $g>0$ $m$-almost everywhere, and put
\[
 a=\int_D\sqrt g\,dm.
\]
For each matching coordinate, add one labeled point at $(0,0)$ when the sample size is odd; use the original sample when it is even. Let $\beta$ be the vector of mean coefficients under uniform sampling, and let $S_1$ be the limiting covariance for uniform samples of fixed size. Then
\begin{equation}\label{eq:applications-fixed-limit}
 H_{n,\mu}-\E H_{n,\mu}
 \xrightarrow{W_2}
 \mathcal N_r\left(0,S_1+\frac{1-a^2}{4}\beta\beta^\top\right).
\end{equation}
For a Poisson cloud with intensity measure $\lambda\mu$ and count $N$,
\begin{equation}\label{eq:applications-poisson-limit}
 \begin{pmatrix}
 P_{\lambda,\mu}-\E P_{\lambda,\mu}\\[2pt]
 (N-\lambda)/\sqrt\lambda
 \end{pmatrix}
 \xrightarrow{W_2}
 \mathcal N_{r+1}\left(
 0,
 \begin{pmatrix}
 S_1+\beta\beta^\top/4 & a\beta/2\\
 a\beta^\top/2 & 1
 \end{pmatrix}\right).
\end{equation}
All covariance entries converge, and
\[
 \E H_{n,\mu}-\E P_{n,\mu}\longrightarrow0.
\]
Every coordinate of $\beta$ and every marginal limiting variance is strictly positive. For the vector $(T,M,L)$ of Euclidean tour, matching, and Held--Karp costs, the fixed-count limiting covariance is positive definite.
\end{theorem}

The theorem permits a different norm in each coordinate. Such a list can have deterministic linear relations, for example when the same cost is computed in proportional norms. For the Euclidean vector $(T,M,L)$, we use the null-space criterion to exclude every nonzero linear combination with zero limiting variance.

\subsection{Metric patching for tours and matchings}

Let $\nu$ be a norm on $\R^2$, and let $a_\nu,b_\nu$ be its minimum and maximum on the Euclidean unit circle. Thus
\[
 a_\nu|x|\le\nu(x)\le b_\nu|x|.
\]
Write $T_\nu$ for the minimum closed-tour length and $M_\nu$ for the minimum perfect-matching length. A tour on two labels traverses their common edge twice. A tour on at most one label and a matching on no labels have cost zero. In a matching, labels at coincident locations remain distinct. For a square $Q$ of side $s$, let $p=p_\nu(Q)$ and $d=d_\nu(Q)$ be its perimeter and diameter in the norm $\nu$. Then
\begin{equation}\label{eq:norm-square-bounds}
 p\le4b_\nu s,\qquad d\le\sqrt2b_\nu s.
\end{equation}

\begin{lemma}[Tour and matching comparisons]\label{lem:tour-match}
Let $I$ be a nonempty interior configuration in $Q$, let $O$ be an exterior configuration, and let $z\in\partial Q$. Then
\begin{equation}\label{eq:tour-boundary}
 -4p\le T_\nu(I\sqcup O)-T_\nu(I)-T_\nu(O\sqcup\{z\})\le2d.
\end{equation}
For an even total number of labels, allowing $I$ to be empty, put $e=|I|\bmod2$. Then
\begin{equation}\label{eq:matching-boundary}
 0\le M_\nu(I\sqcup\{z\}^{[e]})+
 M_\nu(O\sqcup\{z\}^{[e]})-M_\nu(I\sqcup O)\le p.
\end{equation}
Both comparisons extend to labels on the boundary assigned to either side. In particular, $T_\nu$ satisfies Definition~\ref{def:regular} with boundary constant $16b_\nu$, and $M_\nu$ satisfies Definition~\ref{def:parity} with boundary constant $4b_\nu$. These bounds hold for every choice of the boundary point $z$.

For a nonempty tour background $A$,
\begin{equation}\label{eq:tour-insert}
 0\le T_\nu(A\sqcup\{x\})-T_\nu(A)
 \le2\min_{u\in A}\nu(x-u).
\end{equation}
For an even matching background $S$,
\begin{equation}\label{eq:matching-pair-insert}
 |M_\nu(S\sqcup\{x,y\})-M_\nu(S)|\le\nu(x-y).
\end{equation}
\end{lemma}

\begin{proof}
Both costs are translation invariant and homogeneous. At fixed cardinality, each is the minimum of finitely many continuous sums of edge lengths, and is therefore continuous even when locations coincide. Removing a label from a tour and shortcutting cannot increase its length. Conversely, attach a new label $x$ to a nearest existing site by a doubled edge and shortcut the resulting Euler circuit. This proves \eqref{eq:tour-insert}. In particular, adding a label at an existing location leaves the tour cost unchanged.

For matching, add the edge $xy$ to an optimum on $S$ to obtain one inequality in \eqref{eq:matching-pair-insert}. For the other, remove $x,y$ from an optimum on $S\sqcup\{x,y\}$. If they are paired together, discard their edge. Otherwise, join their former partners and apply the triangle inequality along the path through $x,y$. Keeping the partner of a moved label also shows that moving one matching label changes the cost by at most its displacement in the norm.

For the upper bound in \eqref{eq:tour-boundary}, join optimal tours on $I$ and $O\sqcup\{z\}$ by a doubled segment from an interior site to $z$, then shortcut. Its additional cost is at most $2d$. If $O$ is empty, the comparison is an identity.

For the lower bound, suppose $O$ is nonempty and first put the sites and edges in general position with respect to the boundary. Cut an optimal global tour at every boundary crossing. Each resulting component meets the boundary. On each side, add a boundary cycle through the crossing labels and the corners, including $z$ on the exterior side. Each graph is now connected, with all odd-degree vertices on its boundary cycle. Pair the odd vertices in cyclic order and add the intervening boundary arcs, at cost at most $p$ on each side. The two cycles and the added arcs cost at most $4p$ in total. Shortcut the resulting Euler circuits to tours on $I$ and $O\sqcup\{z\}$ to obtain the lower bound. This is the boundary-patching construction used in the partition arguments of \citet{BeardwoodHaltonHammersley1959,Karp1977,Jaillet1992}. Applying the comparison to different exterior configurations and subtracting cancels the interior tour cost, as required in the replacement argument of Section~\ref{sec:regular}.

For the first inequality in \eqref{eq:matching-boundary}, take the union of the two completed matchings. When $e=0$, it is a matching of the original configuration. When $e=1$, remove the two edges incident to the completion labels and join their other endpoints. The triangle inequality shows that this operation does not increase the cost.

For the reverse inequality, retain all edges of a global optimum whose endpoints lie on the same side, including exterior--exterior edges that pass through $Q$. Split each of the $r$ interior--exterior edges at its boundary exit, with a separate boundary label for each edge. Their two pieces have exactly the original cost. Since $r\equiv|I|\pmod2$, the exits together with one copy of $z$ when $e=1$ form an even collection. Order them around the boundary. Of the two alternating pairings of consecutive labels, one has total boundary-arc length at most $p/2$. Add these arcs on both sides and shortcut the resulting paths to matching edges. This adds at most $p$ in total and proves \eqref{eq:matching-boundary}. The same boundary-pairing idea appears in \citet[Lemma~3]{DyerFrieze1984}; the cyclic pairing estimate is also used in \citet[Lemma~2.3]{RueschendorfSachs2000}.

Distinct crossing labels may have the same location and are then joined along zero-length boundary portions. Alternatively, perturb the original sites to avoid tangencies, boundary segments, and repeated crossings, prove the comparisons there, and pass to the limit by continuity. The same limiting argument permits labels on a shared face to be assigned to either side. The constructions use additivity of the norm along a straight segment and the triangle inequality, so they apply to every planar norm. Finally, \eqref{eq:norm-square-bounds} gives the stated constants.
\end{proof}

\subsection{Weighted cuts and elimination of auxiliary labels}

For a labeled configuration $V$ with at least two labels, let $E(V)$ be the edges of its complete loopless graph, and let $\delta(A)$ denote the cut of $A\subset V$. Define
\begin{equation}\label{eq:hk-lp}
 L_\nu(V)=\min\left\{
 \sum_{ij\in E(V)}\nu(x_i-x_j)y_{ij}:
 \begin{array}{ll}
 y_{ij}\ge0,&ij\in E(V),\\
 y(\delta(i))=2,&i\in V,\\
 y(\delta(A))\ge2,&\varnothing\ne A\subsetneq V
 \end{array}\right\}.
\end{equation}
Set $L_\nu(V)=0$ on at most one label. For two labels, the unique edge has capacity two. For at least three labels, the degree equations and the cut of a two-vertex set imply $y_{ij}\le1$. Thus \eqref{eq:hk-lp} is the usual subtour relaxation, with a convention that also covers two labels; see \citet{DantzigFulkersonJohnson1954,HeldKarp1970,HeldKarp1971}.

The feasible polytope is nonempty, compact, and independent of the site positions. The optimal value is therefore continuous at fixed cardinality, translation invariant, and homogeneous. Distinct labels may have zero distance when their locations coincide, so the next lemma allows pseudometric costs.

\begin{lemma}[Eliminating auxiliary labels]\label{lem:eliminate}
Let a finite set of labels carry nonnegative pseudometric costs $c_{ij}$. Suppose a nonnegative edge-capacity vector $u$ satisfies
\[
 u(\delta(A))\ge2
 \qquad(\varnothing\ne A\subsetneq V).
\]
For every terminal set $T\subset V$ with $|T|\ge2$, there is a feasible subtour vector on exactly $T$ whose cost is at most $\sum_{ij}c_{ij}u_{ij}$.
\end{lemma}

\begin{proof}
Complete the graph by giving each missing edge initial capacity zero; its cost remains the prescribed pseudometric distance. We use the parsimonious property of \citet[Theorem~1]{GoemansBertsimas1993}. For symmetric requirements $r_{ij}\ge0$, consider the real linear program with cut constraints
\[
 y(\delta(A))\ge\max_{i\in A,\,j\notin A}r_{ij},
 \qquad y\ge0.
\]
Under the triangle inequality for the costs, prescribing
\[
 y(\delta(i))=\max_{j\ne i}r_{ij}
\]
at any designated vertices does not change its minimum value. Here the capacities are real and nonnegative.

Set $r_{ij}=2$ for distinct terminals and $r_{ij}=0$ otherwise, and prescribe the degrees at every vertex. The given vector $u$ is feasible before these equalities are imposed. An optimum satisfying the degree equalities has degree two at every terminal and degree zero at every auxiliary label. Nonnegativity then makes every edge incident to an auxiliary label zero. For every nonempty proper subset $A\subset T$, the original requirement across $A$ is two, and its cut capacity is unchanged by deleting the auxiliary labels. Thus deletion leaves both the degree equations and every subtour cut on $T$, at no larger cost.

To cover arbitrary real pseudometric costs directly from the rational-data proof in \citet{GoemansBertsimas1993}, for distinct labels replace $c_{ij}$ by
\[
 c_{ij}^{(q)}=\frac{1+\lceil qc_{ij}\rceil}{q}.
\]
For distinct $i,j,k$, the triangle inequality for $c$ gives
\[
 1+\lceil qc_{ik}\rceil
 \le1+\lceil qc_{ij}\rceil+\lceil qc_{jk}\rceil
 \le q(c_{ij}^{(q)}+c_{jk}^{(q)}).
\]
Thus these rational costs satisfy the triangle inequality and converge to $c_{ij}$. Apply the preceding argument to each $c^{(q)}$. The terminal subtour polytope is compact, so a convergent subsequence of the resulting terminal vectors gives the required inequality in the limit. Thus the conclusion also holds when distinct labels have zero distance.
\end{proof}

Lemma~\ref{lem:eliminate} removes auxiliary labels while preserving all terminal cut constraints and controlling the cost. The vertex-subset monotonicity theorem of \citet{ShmoysWilliamson1990} proves the corresponding monotonicity of the Held--Karp value using the 1-tree dual. The parsimonious property also applies to the augmented capacity vectors in the boundary construction below.

\begin{lemma}[Held--Karp comparison]\label{lem:hk}
For every nonempty background $A$,
\begin{equation}\label{eq:hk-insertion}
 0\le L_\nu(A\sqcup\{x\})-L_\nu(A)
 \le2\min_{u\in A}\nu(x-u).
\end{equation}
For the interior, exterior, and anchor of Lemma~\ref{lem:tour-match},
\begin{equation}\label{eq:hk-boundary}
 -2p\le L_\nu(I\sqcup O)-L_\nu(I)-L_\nu(O\sqcup\{z\})\le2d.
\end{equation}
Consequently, $L_\nu$ satisfies Definition~\ref{def:regular} with boundary constant $8b_\nu$, for any boundary anchor. Its insertion bound has constant $2b_\nu$, and adding a label at an existing location leaves its value unchanged.
\end{lemma}

\begin{proof}
If $|A|\ge2$, eliminate the new label from an optimal solution on $A\sqcup\{x\}$ using Lemma~\ref{lem:eliminate}. This proves monotonicity; the singleton case follows from nonnegativity. In the other direction, add an edge of capacity two from $x$ to a nearest old label in an optimal solution on $A$. Every nontrivial cut still has capacity at least two. Apply Lemma~\ref{lem:eliminate} with all labels terminal to obtain \eqref{eq:hk-insertion}. The case of a singleton $A$ is exactly the two-label convention. In particular, adding a coincident label leaves the value unchanged.

For the upper bound in \eqref{eq:hk-boundary}, join optimal vectors on $I$ and $O\sqcup\{z\}$ by an edge of capacity two from an interior site to $z$. A nontrivial cut either splits one of the two vertex sets, and already has capacity at least two, or separates the two sets, and crosses the new edge. Eliminate the auxiliary anchor while retaining the original labels. The joining cost is at most $2d$. If $O$ is empty, the comparison is an identity; singleton components can be treated as isolated vertices before the joining edge is added.

For the lower bound, take an optimal global vector. Keep the interior--interior and exterior--exterior edges on their respective sides, and split every interior--exterior edge at a separately labeled boundary exit. Give both pieces the original capacity; their costs add to the original edge cost. On each side add a capacity-one boundary cycle through all exits and the four corners, including $z$ on the exterior side.

We check every nontrivial cut in each augmented graph. If it separates boundary-cycle vertices, the cycle alone crosses it with capacity at least two. Otherwise, choose the side of the cut that contains no cycle vertex. It is a nonempty subset of the original interior labels, or of the original exterior labels, and its incident capacities are exactly those of the same cut in the global vector. Their sum is at least two. Lemma~\ref{lem:eliminate} now produces a subtour vector on $I$ and one on $O\sqcup\{z\}$, at total added cost at most $2p$. When the interior has one label, its required cost is zero and the same inequality follows by discarding its nonnegative augmented cost.

Separate cycle labels are retained even when their positions coincide; intervening zero-cost edges are allowed. Boundary sites are handled by continuity. The capacity-two joining operation, capacity-one boundary cycle, and nearest-neighbor insertion are the constructions of \citet[Propositions~7--8 and Lemma~11]{GoemansBertsimas1991}, here used to verify the interior--exterior comparison. The constants follow from \eqref{eq:norm-square-bounds}.
\end{proof}

\subsection{Positive means and positive definiteness in the Euclidean case}

\begin{proof}[Proof of Theorem~\ref{thm:applications}]
Lemmas~\ref{lem:tour-match} and \ref{lem:hk} verify the boundary comparisons in Definitions~\ref{def:regular} and \ref{def:parity}. Theorem~\ref{thm:fixed-law} then gives \eqref{eq:applications-fixed-limit}, \eqref{eq:applications-poisson-limit}, convergence of the covariance matrices, and convergence to zero of the difference between the fixed-count and Poisson expectations.

To prove that the uniform-sampling mean coefficients are positive, for $n\ge2$ sites let
\[
 R_i=\min_{j\ne i}|x_i-x_j|.
\]
A nonnegative edge vector with weighted degree $d_0$ at every vertex satisfies
\begin{equation}\label{eq:applications-nearest-lower}
 \sum_{i<j}y_{ij}\nu(x_i-x_j)
 \ge\frac{a_\nu}{2}\sum_{i<j}y_{ij}(R_i+R_j)
 =\frac{a_\nu d_0}{2}\sum_iR_i.
\end{equation}
For iid uniform sites and $r=(2\pi n)^{-1/2}$, the conditional union bound gives
\[
 \Pp(R_i>r\mid x_i)\ge1-(n-1)\pi r^2\ge\frac12.
\]
Hence the expected cost in \eqref{eq:applications-nearest-lower} is at least
\[
 \frac{a_\nu d_0}{4\sqrt{2\pi}}\sqrt n.
\]
Use $d_0=2$ for tours and Held--Karp solutions, and $d_0=1$ for perfect matchings at even counts. The homogeneous mean limits in Theorem~\ref{thm:fixed-law} thus have strictly positive coefficients. Theorem~\ref{thm:kernel} implies strictly positive homogeneous marginal variances.

To prove joint positive definiteness, first consider the Euclidean regular pair $(T,L)$. Theorem~1.8 of \citet{FriezePegden2017}, together with the preceding lower bound, gives
\begin{equation}\label{eq:applications-mean-separation}
 0<\beta_L<\beta_T.
\end{equation}
Their relaxation agrees with \eqref{eq:hk-lp}: under the degree equations, the cut inequality is equivalent to $y(E(A))\le|A|-1$, and the bounds $y_{ij}\le1$ follow as above. Their almost sure coefficients equal $\beta_T$ and $\beta_L$, since Theorem~\ref{thm:fixed-law} gives the corresponding costs divided by $\sqrt n$ converging to these coefficients in $L^2$, and limits in probability are unique.

We also need configurations on which the two costs grow at the same rate. Let $S_k$ be the $k\times k$ grid in $[1/4,3/4]^2$, with spacing $\delta_k=1/(2(k-1))$. Every feasible subtour vector has total edge capacity $k^2$, and every edge has length at least $\delta_k$. A tour can visit the grid along successive rows and then return to its starting point. Therefore
\begin{equation}\label{eq:applications-grid-direction}
 k^2\delta_k\le L(S_k)\le T(S_k)
 \le(k^2-1)\delta_k+\sqrt2.
\end{equation}
In particular, both $T(S_k)/(k^2\delta_k)$ and $L(S_k)/(k^2\delta_k)$ tend to one.

Suppose that $uT+vL$ is bounded on all finite labeled configurations in $D$. Evaluate it on $S_k$ and divide by $k^2\delta_k$ to obtain $u+v=0$. Taking expectations under uniform sampling and dividing by $\sqrt n$ gives
$u\beta_T+v\beta_L=0$. By \eqref{eq:applications-mean-separation}, these two equations imply $u=v=0$. Theorem~\ref{thm:kernel} therefore gives rank two for the regular covariance block. Matching has positive limiting variance by the preceding mean lower bound. Corollary~\ref{cor:split-kernel}, applied after grouping $(T,L)$ and $M$, gives rank three for the Euclidean vector $(T,M,L)$.

Finally, for the law $\mu$, the fixed-count covariance is
\[
 S_1+\frac{1-a^2}{4}\beta\beta^\top.
\]
Since $a^2\le1$, the added matrix is positive semidefinite. The marginal variances therefore remain positive, and the covariance of the Euclidean vector remains positive definite.
\end{proof}

In particular, the additive difference $T-L$ between the optimal Euclidean tour length and the Held--Karp lower bound has a nondegenerate centered Gaussian limit. Applying the linear map $(x,y,z)\mapsto x-z$ to Theorem~\ref{thm:applications} gives limiting variance
\[
 (1,0,-1)S_1(1,0,-1)^\top
 +\frac{1-a^2}{4}(\beta_T-\beta_L)^2>0.
\]
This conclusion uses the joint covariance result; separate central limit theorems for the tour and relaxation would not determine the fluctuations of their difference.

\section{Scope of the criterion}\label{sec:scope}

We finish by clarifying the assumptions behind the method. First, we show that the deterministic boundary comparison is equivalent to a four-term comparison of replacement responses. We then give an example in which conditional replacements have small responses but the conditional mean still fluctuates. Finally, we construct a cost satisfying the probabilistic criterion with no uniform deterministic boundary bound.

\subsection{Equivalent deterministic comparisons}

In this subsection $H$ is defined on configurations in $\R^2$ and does not depend on an enclosing square. Disjoint unions retain every label, including labels at coincident locations. Squares are closed, and exterior sites lie strictly outside them.

\begin{proposition}[Deterministic boundary comparison]\label{prop:deterministic-boundary}
Let $H$ be finite and real valued on the finite labeled configurations in $\R^2$. Suppose that it is invariant under permutations and that $H(\{z\})=0$ for every $z$. The following uniform properties are equivalent.
\begin{enumerate}
\item For every axis-parallel square $Q$ of side $b$, a prescribed point $z_Q\in\partial Q$ satisfies
\begin{equation}\label{eq:deterministic-anchor}
 |H(I\sqcup O)-H(I)-H(O\sqcup\{z_Q\})|\le Bb
\end{equation}
for all nonempty finite $I\subset Q^\circ$ and finite $O\subset\R^2\setminus Q$.
\item For every such square, all nonempty finite $I,I'\subset Q^\circ$, and all finite $O,O'\subset\R^2\setminus Q$,
\begin{equation}\label{eq:deterministic-mixed}
 |H(I\sqcup O)-H(I'\sqcup O)-H(I\sqcup O')+H(I'\sqcup O')|\le Cb.
\end{equation}
\end{enumerate}
The first property implies the second with $C=4B$. The second implies the first with $B=C$ for every boundary point. For the reverse implication it suffices to assume \eqref{eq:deterministic-mixed} when $O'=\varnothing$ and $I'$ is a singleton. The equivalence holds without a continuity assumption.
\end{proposition}

\begin{proof}
Apply \eqref{eq:deterministic-anchor} to each of the four terms in \eqref{eq:deterministic-mixed}. The interior costs and the anchored exterior costs cancel, leaving a bound $4Bb$.

Conversely, fix $Q,I,O$ and $z\in\partial Q$. Since $O$ is finite and disjoint from the closed square $Q$, it remains outside the concentric square of side $b+\epsilon$ for every sufficiently small $\epsilon>0$. Both $I$ and $\{z\}$ lie strictly inside that larger square. Apply \eqref{eq:deterministic-mixed} with $I'=\{z\}$ and $O'=\varnothing$. Since $H(\{z\})=0$, its left side is
\[
 |H(I\sqcup O)-H(I)-H(O\sqcup\{z\})|.
\]
It is bounded by $C(b+\epsilon)$ for every sufficiently small $\epsilon$, hence by $Cb$. The configurations in this last inequality do not change with $\epsilon$.
\end{proof}

Translation invariance and degree-one homogeneity imply the singleton normalization used in this proposition.

\subsection{Why the conditional mean must be controlled}

\begin{example}\label{ex:state-fluctuation}
For a finite configuration $S$, set
\[
 H(S)=\1_{\{|S|\text{ odd}\}}\operatorname{diam}(S),\qquad
 \operatorname{diam}(\varnothing)=0.
\]
This cost is translation invariant, degree-one homogeneous, and continuous at each fixed cardinality. If $I,I'$ are nonempty configurations in a square of side $b$ with equal count parity, then for every common exterior $O$,
\begin{equation}\label{eq:state-example-response}
 |H(I\sqcup O)-H(I'\sqcup O)|\le\sqrt2\,b.
\end{equation}
The four-term response to two exteriors is consequently at most $2\sqrt2\,b$. Nevertheless, under homogeneous Poisson sampling on $D$,
\begin{equation}\label{eq:state-example-law}
 H(\mathcal P_\lambda)-\E H(\mathcal P_\lambda)
 \Longrightarrow\sqrt2\,(P-\tfrac12),\qquad P\sim\operatorname{Bernoulli}(\tfrac12).
\end{equation}
For a fixed $k\times k$ grid with $k\ge2$, let $\pi$ record the actual cell-count parities. The cutting remainder satisfies
\begin{equation}\label{eq:state-example-variances}
 \E\Var(R_{\lambda,k}\mid\pi)\longrightarrow0,\qquad
 \Var\bigl(\E[R_{\lambda,k}\mid\pi]\bigr)\longrightarrow1.
\end{equation}
\end{example}

\begin{proof}
The parity indicators in \eqref{eq:state-example-response} agree. An extremal diameter pair wholly in $O$ is unchanged. A pair wholly in the interior has length at most $\sqrt2\,b$. For a cross pair, retain its exterior endpoint and replace its interior endpoint by any point of $I'$; the displacement is at most $\sqrt2\,b$. Interchanging $I$ and $I'$ proves \eqref{eq:state-example-response}. Apply this estimate to each exterior to get the four-term bound.

The diameter of $\mathcal P_\lambda$ tends to $\sqrt2$ in $L^2$: fixed neighborhoods of opposite corners are occupied with probability tending to one, and the diameter is bounded by $\sqrt2$. Also
\[
 \E(-1)^{N_\lambda}=e^{-2\lambda},
\]
so the count parity tends to a fair bit. Boundedness gives the centered limit \eqref{eq:state-example-law}.

Write $\xi_{\lambda,i}$ for the parity in cell $i$, and $\xi_\lambda$ for their sum modulo two. The cell counts are independent, and their parities tend jointly to independent fair bits $P_1,\ldots,P_{k^2}$. Their diameters tend to $\sqrt2/k$ in $L^2$. Thus
\[
 R_{\lambda,k}-\left(\sqrt2\,\xi_\lambda-\frac{\sqrt2}{k}\sum_{i=1}^{k^2}\xi_{\lambda,i}\right)
 \longrightarrow0\quad\text{in }L^2.
\]
The expression in parentheses is $\pi$-measurable. Its variance tends to that of
\[
 \sqrt2\,P-\frac{\sqrt2}{k}\sum_{i=1}^{k^2}P_i,
 \qquad P=P_1\oplus\cdots\oplus P_{k^2}.
\]
As $k^2\ge2$, $P$ is independent of each $P_i$. The variance is therefore
\[
 \frac24+\frac2{k^2}\frac{k^2}4=1.
\]
Conditional expectation is an $L^2$ contraction, giving the first limit in \eqref{eq:state-example-variances}; total variance gives the second.
\end{proof}

In this example, replacements preserve count parity and have nonempty interiors. The example concerns the conditional replacement argument itself, before imposing the completion rule of Definition~\ref{def:parity}. It shows why \eqref{eq:conditional-cut} retains the variance of the conditional mean: the conditional variance can tend to zero while this second term has a positive limit.

\subsection{A probabilistic comparison without a deterministic bound}

\begin{proposition}\label{prop:rare-collinear}
There is a translation-invariant, degree-one homogeneous cost $J$, continuous at every fixed cardinality, which satisfies \eqref{eq:boundary-excess} and tends to zero under homogeneous Poisson sampling in $L^p$ for every fixed $p\ge1$, but
\[
 \sup_{S\subset D}|J(S)|=\infty,
\]
where the supremum is over finite labeled configurations, and $J$ admits no uniform deterministic comparison of the form \eqref{eq:deterministic-anchor}.
\end{proposition}

\begin{proof}
We make the cost large only on configurations lying in a narrow strip. For a nonempty configuration put
\[
 d(S)=\operatorname{diam}(S),\qquad
 w(S)=\min_{|u|=1}\left(\max_{x\in S}u\cdot x-\min_{x\in S}u\cdot x\right).
\]
Choose $0<\epsilon<1/16$ and a continuous $\varphi:[0,\infty)\to[0,1]$ with $\varphi(0)=1$ and $\varphi(t)=0$ for $t\ge\epsilon$. Set
\[
 J(S)=\sqrt{|S|}\,d(S)\varphi\!\left(\frac{w(S)}{d(S)}\right)
 \quad\text{if }|S|\ge3\text{ and }d(S)>0,
\]
and set $J(S)=0$ otherwise. Projection widths vary uniformly continuously when finitely many sites are displaced. Thus $w$ and $d$ are continuous at fixed cardinality; at $d=0$, use $|J(S)|\le\sqrt{|S|}\,d(S)$. The formula is permutation invariant, translation invariant and homogeneous. On $D$,
\[
 |J(S)|\le\sqrt{2|S|}.
\]

The triangle with vertices $(0,0),(1,0),(0,1)$ has minimum width $1/\sqrt2$. Indeed, for a unit vector $(u,v)$, the range of $0,u,v$ is at least $\max\{|u|,|v|\}\ge1/\sqrt2$, with equality at $u=v=1/\sqrt2$. Intersect $D$ with the three balls of radius $1/16$ about these vertices. They have the same positive area $\alpha$. If all three sets contain a sampled point, then
\[
 w(S)\ge\frac1{\sqrt2}-\frac18,\qquad d(S)\le\sqrt2,
\]
and therefore $J(S)=0$. If one specified set is empty under intensity $s$, its probability is $e^{-\alpha s}$ and the number of points outside it is an independent Poisson variable with mean $(1-\alpha)s$. A union bound and Poisson moment estimates give, for each fixed $p\ge1$,
\begin{equation}\label{eq:rare-cost-moments}
 \E|J(\mathcal P_s)|^p\le C_p(1+s^{p/2})e^{-\alpha s},\qquad s\ge1.
\end{equation}
In particular, $q(s)=C(1+s)e^{-\alpha s}$ bounds the second moment.

To check the replacement criterion, use $J$ as the cost in every square. Let $Q$ have side $1/k$, let $B\subset Q$ have side $b$, and put $s=\lambda/k^2$ and $u=\lambda b^2\ge1$. The original and replaced global configurations both have the Poisson law of intensity $\lambda$, while each local cost has the law of $k^{-1}J(\mathcal P_s)$. Applying $|x_1+\cdots+x_4|^2\le4\sum_{i=1}^4|x_i|^2$ to the four terms in $\Delta_{B,Q}J$ gives
\[
 \E|\Delta_{B,Q}J|^2\le8\{q(\lambda)+k^{-2}q(s)\}.
\]
Since $\lambda\ge s\ge u$, division by $b^2=u/\lambda$ yields
\begin{equation}\label{eq:rare-cost-excess}
 b^{-2}\E|\Delta_{B,Q}J|^2
 \le\frac{16}{u}\sup_{v\ge u}vq(v)\longrightarrow0.
\end{equation}
This verifies \eqref{eq:boundary-excess}, uniformly over its scales, for every fixed positive boundary constant.

For $n\ge3$ distinct collinear sites of diameter one in $D$, $J(S)=\sqrt n$, proving unboundedness. Finally, fix a square and any proposed boundary anchor. Choose a line through the anchor and the interior of the square, two interior sites on the line, and one exterior site on it tending to infinity. The interior cost and the exterior cost completed by the anchor both vanish by cardinality. The full three-site cost equals $\sqrt3$ times its diameter and is unbounded. Hence \eqref{eq:deterministic-anchor} fails.
\end{proof}

\appendix
\section{Space-filling curves and prefix parities}
\label{app:curve-prefix}

\begin{proof}[Proof of Lemma~\ref{lem:hilbert-order}]
For a southwest-to-southeast traversal, order the
children southwest, northwest, northeast, southeast, with successive
endpoints
\[
 (0,0),\quad(0,\tfrac12),\quad
 (\tfrac12,\tfrac12),\quad(1,\tfrac12),\quad(1,0).
\]
Repeat this construction by isometries on the four equal parameter
intervals. Nested closed squares define the limit map, and the endpoint
rules agree where consecutive parameter intervals meet. The map sends
each level-$j$ parameter interval onto its corresponding dyadic square:
every point in that square belongs to a nested sequence of smaller
squares. In particular, the map is onto $D$. Symmetries of $D$ give any
prescribed ordered pair of adjacent corners.
If \(4^{-j-1}<|u-v|\le4^{-j}\), the two images belong to at most two
consecutive level-\(j\) squares, proving
\eqref{eq:parity-hilbert}.
The first-preimage map $b$ is Borel, since \(\{x:b(x)\le t\}=\phi([0,t])\) is compact.

By definition, $\phi(b(x))=x$ for every $x\in D$.
\end{proof}

\begin{proof}[Proof of Lemma~\ref{lem:parity-prefix}]
Couple future chains maximally from opposite current states,
keeping them together after their first agreement.  Their probability
of disagreement \(j-i\) steps later is at most \(\rho^{j-i}\).
With the preceding states fixed, changing the \(i\)-th state changes
the conditional expectation of \(f\) by at most
\[
 b_i=\sum_{j\ge i}\rho^{j-i}a_j.
\]
Its centered Doob increment has conditional variance at most \(b_i^2/4\)
and conditional moment generating function at most
\(\exp(t^2b_i^2/8)\), by the bounded-variable exponential estimate of
\citet{Hoeffding1963}.
The \(\ell^2\) triangle inequality gives
\(\sum_i b_i^2\le(1-\rho)^{-2}\sum_i a_i^2\).
Orthogonality of the Doob increments proves the variance estimate;
iteration of the
conditional exponential estimate and optimization of the exponential
Markov bound prove the tail estimate.
\end{proof}

\bibliographystyle{apalike}
\bibliography{boundary_method}

\begin{thebibliography}{}

\bibitem[Aizenman and Wehr, 1990]{AizenmanWehr1990Rounding}
Aizenman, M. and Wehr, J. (1990).
\newblock {Rounding effects of quenched randomness on first-order phase
  transitions}.
\newblock {\em Communications in Mathematical Physics}, 130:489--528.

\bibitem[Avram and Bertsimas, 1993]{AvramBertsimas1993}
Avram, F. and Bertsimas, D. (1993).
\newblock {On central limit theorems in geometrical probability}.
\newblock {\em The Annals of Applied Probability}, 3(4):1033--1046.

\bibitem[Bartholdi and Platzman, 1988]{BartholdiPlatzman1988Spacefilling}
Bartholdi, III, J.~J. and Platzman, L.~K. (1988).
\newblock {Heuristics Based on Spacefilling Curves for Combinatorial Problems
  in Euclidean Space}.
\newblock {\em Management Science}, 34(3):291--305.

\bibitem[Baryshnikov and Yukich, 2005]{BaryshnikovYukich2005}
Baryshnikov, Y. and Yukich, J.~E. (2005).
\newblock {Gaussian limits for random measures in geometric probability}.
\newblock {\em The Annals of Applied Probability}, 15:213--253.

\bibitem[Beardwood et~al., 1959]{BeardwoodHaltonHammersley1959}
Beardwood, J., Halton, J.~H., and Hammersley, J.~M. (1959).
\newblock {The shortest path through many points}.
\newblock {\em Proceedings of the Cambridge Philosophical Society},
  55:299--327.

\bibitem[Boucheron et~al., 2005]{BoucheronBousquetLugosiMassart2005}
Boucheron, S., Bousquet, O., Lugosi, G., and Massart, P. (2005).
\newblock {Moment inequalities for functions of independent random variables}.
\newblock {\em The Annals of Probability}, 33(2):514--560.

\bibitem[Chatterjee, 2008]{Chatterjee2008NormalApproximation}
Chatterjee, S. (2008).
\newblock {A new method of normal approximation}.
\newblock {\em The Annals of Probability}, 36:1584--1610.

\bibitem[Chatterjee, 2014]{Chatterjee2014Survey}
Chatterjee, S. (2014).
\newblock {A short survey of Stein's method}.
\newblock In {\em Proceedings of the International Congress of
  Mathematicians---Seoul 2014}, volume~IV, pages 1--24, Seoul. Kyung Moon Sa.

\bibitem[Chatterjee, 2018]{Chatterjee2018RFIMDecay}
Chatterjee, S. (2018).
\newblock {On the decay of correlations in the random field Ising model}.
\newblock {\em Communications in Mathematical Physics}, 362:253--267.

\bibitem[Chatterjee and Sen, 2017]{ChatterjeeSen2017MST}
Chatterjee, S. and Sen, S. (2017).
\newblock {Minimal spanning trees and Stein's method}.
\newblock {\em The Annals of Applied Probability}, 27:1588--1645.

\bibitem[Chen et~al., 2011]{ChenGoldsteinShao2011}
Chen, L. H.~Y., Goldstein, L., and Shao, Q.-M. (2011).
\newblock {\em {Normal Approximation by Stein's Method}}.
\newblock Probability and Its Applications. Springer, Heidelberg.

\bibitem[Dantzig et~al., 1954]{DantzigFulkersonJohnson1954}
Dantzig, G.~B., Fulkerson, D.~R., and Johnson, S.~M. (1954).
\newblock {Solution of a Large-Scale Traveling-Salesman Problem}.
\newblock {\em Journal of the Operations Research Society of America},
  2:393--410.

\bibitem[Dyer and Frieze, 1984]{DyerFrieze1984}
Dyer, M.~E. and Frieze, A.~M. (1984).
\newblock {A partitioning algorithm for minimum weighted Euclidean matching}.
\newblock {\em Information Processing Letters}, 18(2):59--62.

\bibitem[Efron and Stein, 1981]{EfronStein1981}
Efron, B. and Stein, C. (1981).
\newblock {The jackknife estimate of variance}.
\newblock {\em The Annals of Statistics}, 9(3):586--596.

\bibitem[Feng, 2026]{Feng2026MatchingVariance}
Feng, S. (2026).
\newblock {Variance Bounds for Minimum Perfect Matching Problem by Random Walk
  on Multinomial Distribution}.
\newblock arXiv:2610.14337, version 2.

\bibitem[Frieze and Pegden, 2017]{FriezePegden2017}
Frieze, A. and Pegden, W. (2017).
\newblock {Separating subadditive Euclidean functionals}.
\newblock {\em Random Structures \& Algorithms}, 51:375--403.

\bibitem[Gao and Steele, 1994]{GaoSteele1994Spacefilling}
Gao, J. and Steele, J.~M. (1994).
\newblock {Sums of Squares of Edge Lengths and Spacefilling Curve Heuristics
  for the Traveling Salesman Problem}.
\newblock {\em SIAM Journal on Discrete Mathematics}, 7(2):314--324.

\bibitem[Goemans and Bertsimas, 1991]{GoemansBertsimas1991}
Goemans, M.~X. and Bertsimas, D.~J. (1991).
\newblock {Probabilistic Analysis of the Held and Karp Lower Bound for the
  Euclidean Traveling Salesman Problem}.
\newblock {\em Mathematics of Operations Research}, 16:72--89.

\bibitem[Goemans and Bertsimas, 1993]{GoemansBertsimas1993}
Goemans, M.~X. and Bertsimas, D.~J. (1993).
\newblock {Survivable networks, linear programming relaxations and the
  parsimonious property}.
\newblock {\em Mathematical Programming}, 60:145--166.

\bibitem[Held and Karp, 1970]{HeldKarp1970}
Held, M. and Karp, R.~M. (1970).
\newblock {The Traveling-Salesman Problem and Minimum Spanning Trees}.
\newblock {\em Operations Research}, 18:1138--1162.

\bibitem[Held and Karp, 1971]{HeldKarp1971}
Held, M. and Karp, R.~M. (1971).
\newblock {The traveling-salesman problem and minimum spanning trees: Part II}.
\newblock {\em Mathematical Programming}, 1:6--25.

\bibitem[Hoeffding, 1948]{Hoeffding1948}
Hoeffding, W. (1948).
\newblock {A class of statistics with asymptotically normal distribution}.
\newblock {\em The Annals of Mathematical Statistics}, 19(3):293--325.

\bibitem[Hoeffding, 1963]{Hoeffding1963}
Hoeffding, W. (1963).
\newblock {Probability inequalities for sums of bounded random variables}.
\newblock {\em Journal of the American Statistical Association},
  58(301):13--30.

\bibitem[Holroyd and Liggett, 2001]{HolroydLiggett2001ExtraHead}
Holroyd, A.~E. and Liggett, T.~M. (2001).
\newblock {How to find an extra head: optimal random shifts of Bernoulli and
  Poisson random fields}.
\newblock {\em The Annals of Probability}, 29(4):1405--1425.

\bibitem[Jaillet, 1992]{Jaillet1992}
Jaillet, P. (1992).
\newblock {Rates of Convergence for Quasi-Additive Smooth Euclidean Functionals
  and Application to Combinatorial Optimization Problems}.
\newblock {\em Mathematics of Operations Research}, 17:964--980.

\bibitem[Janson, 2007]{Janson2007Conditioning}
Janson, S. (2007).
\newblock {Monotonicity, asymptotic normality and vertex degrees in random
  graphs}.
\newblock {\em Bernoulli}, 13(4):952--965.

\bibitem[Karp, 1977]{Karp1977}
Karp, R.~M. (1977).
\newblock {Probabilistic Analysis of Partitioning Algorithms for the
  Traveling-Salesman Problem in the Plane}.
\newblock {\em Mathematics of Operations Research}, 2:209--224.

\bibitem[Kesten and Lee, 1996]{KestenLee1996MST}
Kesten, H. and Lee, S. (1996).
\newblock {The central limit theorem for weighted minimal spanning trees on
  random points}.
\newblock {\em The Annals of Applied Probability}, 6:495--527.

\bibitem[Koo and Lee, 2007]{KooLee2007}
Koo, Y. and Lee, S. (2007).
\newblock {Rates of Convergence of Means of Euclidean Functionals}.
\newblock {\em Journal of Theoretical Probability}, 20:821--841.

\bibitem[Lachi{\`e}ze-Rey et~al., 2022]{LachiezeReyPeccatiYang2022}
Lachi{\`e}ze-Rey, R., Peccati, G., and Yang, X. (2022).
\newblock {Quantitative two-scale stabilization on the Poisson space}.
\newblock {\em The Annals of Applied Probability}, 32:3085--3145.

\bibitem[Leskel{\"a}, 2024]{Leskela2024}
Leskel{\"a}, L. (2024).
\newblock {Information divergences and likelihood ratios of Poisson processes
  and point patterns}.
\newblock {\em IEEE Transactions on Information Theory}, 70:9084--9101.

\bibitem[Penrose, 2005]{Penrose2005MultivariateSpatial}
Penrose, M.~D. (2005).
\newblock {Multivariate spatial central limit theorems with applications to
  percolation and spatial graphs}.
\newblock {\em The Annals of Probability}, 33:1945--1991.

\bibitem[Penrose, 2007]{Penrose2007}
Penrose, M.~D. (2007).
\newblock {Gaussian limits for random geometric measures}.
\newblock {\em Electronic Journal of Probability}, 12:989--1035.

\bibitem[Penrose and Yukich, 2001]{PenroseYukich2001Graphs}
Penrose, M.~D. and Yukich, J.~E. (2001).
\newblock {Central limit theorems for some graphs in computational geometry}.
\newblock {\em The Annals of Applied Probability}, 11:1005--1041.

\bibitem[Petrov, 2020]{Petrov2020ParityCoupling}
Petrov, F. (2020).
\newblock {Answer to: Coupling a binomial -- parity conditioning}.
\newblock MathOverflow, \url{https://mathoverflow.net/questions/379062}.

\bibitem[Pinelis, 2020]{Pinelis2020ParityCoupling}
Pinelis, I. (2020).
\newblock {Answer to: Coupling a binomial -- parity conditioning}.
\newblock MathOverflow, \url{https://mathoverflow.net/questions/379062}.

\bibitem[Redmond and Yukich, 1994]{RedmondYukich1994}
Redmond, C. and Yukich, J.~E. (1994).
\newblock {Limit Theorems and Rates of Convergence for Euclidean Functionals}.
\newblock {\em The Annals of Applied Probability}, 4:1057--1073.

\bibitem[Redmond and Yukich, 1996]{RedmondYukich1996}
Redmond, C. and Yukich, J.~E. (1996).
\newblock {Asymptotics for Euclidean functionals with power-weighted edges}.
\newblock {\em Stochastic Processes and their Applications}, 61:289--304.

\bibitem[R{\"u}schendorf and Sachs, 2000]{RueschendorfSachs2000}
R{\"u}schendorf, L. and Sachs, G. (2000).
\newblock {Stochastic analysis of partitioning algorithms for matching
  problems}.
\newblock {\em Journal of Applied Probability}, 37(2):494--503.

\bibitem[Santambrogio, 2015]{Santambrogio2015Transport}
Santambrogio, F. (2015).
\newblock {\em {Optimal Transport for Applied Mathematicians: Calculus of
  Variations, PDEs, and Modeling}}, volume~87 of {\em Progress in Nonlinear
  Differential Equations and Their Applications}.
\newblock Birkh{\"a}user, Cham.

\bibitem[Shmoys and Williamson, 1990]{ShmoysWilliamson1990}
Shmoys, D.~B. and Williamson, D.~P. (1990).
\newblock {Analyzing the Held-Karp TSP bound: A monotonicity property with
  application}.
\newblock {\em Information Processing Letters}, 35:281--285.

\bibitem[Steele, 1981]{Steele1981}
Steele, J.~M. (1981).
\newblock {Subadditive Euclidean Functionals and Nonlinear Growth in Geometric
  Probability}.
\newblock {\em The Annals of Probability}, 9:365--376.

\bibitem[Trinh, 2019]{Trinh2019}
Trinh, K.~D. (2019).
\newblock {On central limit theorems in stochastic geometry for add-one cost
  stabilizing functionals}.
\newblock {\em Electronic Communications in Probability}, 24(76):1--15.

\bibitem[Wehr and Aizenman, 1990]{WehrAizenman1990Extensive}
Wehr, J. and Aizenman, M. (1990).
\newblock {Fluctuations of extensive functions of quenched random couplings}.
\newblock {\em Journal of Statistical Physics}, 60:287--306.

\bibitem[Yukich, 1998]{Yukich1998}
Yukich, J.~E. (1998).
\newblock {\em {Probability Theory of Classical Euclidean Optimization
  Problems}}, volume 1675 of {\em Lecture Notes in Mathematics}.
\newblock Springer, Berlin.

\bibitem[Yukich, 2013]{Yukich2013Survey}
Yukich, J.~E. (2013).
\newblock {Limit Theorems in Discrete Stochastic Geometry}.
\newblock In Spodarev, E., editor, {\em Stochastic Geometry, Spatial Statistics
  and Random Fields: Asymptotic Methods}, volume 2068 of {\em Lecture Notes in
  Mathematics}, pages 239--275. Springer, Heidelberg.

\end{thebibliography}
\end{document}